\documentclass[12pt,reqno,oneside]{amsart}
 \usepackage{verbatim,amsmath,amssymb,cite,xspace}
\usepackage{color,cite,graphicx}

\theoremstyle{plain}
\newtheorem{theorem}{Theorem}[section]
\newtheorem{lemma}{Lemma}[section]
\newtheorem{corollary}{Corollary}[section]
\theoremstyle{definition}

\theoremstyle{remark}

\usepackage{color}

\usepackage{amssymb}
\usepackage{amssymb, latexsym}
\usepackage{amsmath}
\usepackage{amscd}
\usepackage{enumerate}
\usepackage{amsfonts}
\usepackage{graphicx}
\usepackage{epsfig}

\newcommand{\tk}{\tilde{t}_k}
\newcommand{\HL}{\dot{H}^1\times L^2}
\newcommand{\OR}{\overrightarrow}
\newcommand{\E}{\mathcal{E}}
\newcommand{\A}{\mathcal{A}}
\newcommand{\B}{\mathcal{B}}
\newcommand{\Fl}{{\rm Flux}}
\newcommand{\tpsi}{\widetilde{\psi}}
\newcommand{\tv}{\tilde{v}}
\newcommand{\bpsi}{\breve{\psi}}
\newcommand{\tg}{\tilde{g}}
\newcommand{\epsi}{\psi^{\epsilon}}
\numberwithin{equation}{section} 

\begin{document}
\title[]{Asymptotic decomposition for semilinear wave and equivariant wave map equations}
\author{Hao Jia}
\author{Carlos Kenig}
\thanks{CEK was partially supported by NSF grants DMS 0968472 and 1265249}

\maketitle
\noindent
{\bf Abstract.}  In this paper we give a unified proof to the soliton resolution conjecture along a sequence of times, for the semilinear focusing energy critical wave equations in the radial case and two dimensional equivariant wave map equations, including the four dimensional radial Yang Mills equation, without using outer energy type inequalities. Such inequalities have played a crucial role in previous works with similar results. Roughly speaking, we prove that along a sequence of times $t_n\to T_+$ (the maximal time of existence), the solution decouples to a sum of rescaled solitons and a term vanishing in the energy space, plus a free radiation term in the global case or a regular part in the finite time blow up case.
The main difficulty is that in general (for instance for the radial four dimensional Yang Mills case and the radial six dimensional semilinear wave case) we do not have a favorable outer energy inequality for the associated linear wave equations. Our main new input is the simultaneous use of two virial identities. \\

\begin{section}{Introduction}
In this paper we study soliton resolution for radial energy critical semilinear wave equations, $2+1$ dimensional equivariant wave maps and $4+1$ dimensional radial Yang Mills equations. These equations are locally wellposed for initial data in the natural energy space (see \cite{ChriTah,Struwe2,KenigYang,Shatah1,Shatah2} and Section 2 below), and admit nontrivial steady states, called {\it solitons}, which together with scaling symmetry play a prominent role in the singularity formation and long time dynamics of solutions. Indeed, singular solutions in the form of a shrinking soliton plus a residue term have been constructed for the energy critical wave equation in \cite{KSTwave} for $3+1$ dimensional case and in \cite{HillRap} for the $4+1$ dimensional case, for the $2+1$ equivariant dimensional wave map equation in \cite{KSTwavemap,RodSter}, and for the radial four dimensional Yang Mills equation in \cite{KSTmills,RapRod}. We also refer to the recent survey \cite{SchlagICM} for further discussion. \\

The above described behavior is expected to hold for all solutions (possibly with multi solitons). In fact the {\it soliton resolution conjecture} predicts, in this context, that modulo a free radiation in the global case or a regular part in the finite blow up case the solution asymptotically decouples to a sum of rescaled solitons and a term that vanishes as time approaches the maximal time of existence. (In the case of semilinear focusing wave equations, one has to assume in addition that the solutions remain bounded in $\HL(R^d)$, i.e., the solution is Type II.) This conjecture remains largely open, except when one assumes additional size restriction on the solutions, such as closeness to solitons. Recently a breakthrough was made in \cite{DKM} for radial three dimensional energy critical focusing wave equation, where the soliton resolution conjecture was settled. The works \cite{DKM,DKMsmall} introduced a very natural and powerful tool, called ``channels of energy", which among many other things, implies that if a solution is not a soliton, then it will radiate energy to large distances, thus proving the crucial ``rigidity theorem" for compact solutions. The arguments followed from energy inequalities for the associated linear wave equation, and have since been extended to general dimensions \cite{KLIS1,RKS} and applied to many other situations, see for example \cite{Cotemap1,Cotemap2,KLS,KLIS1,KLIS2} for applications to wave maps, \cite{Kenig4dWave,DKMsuper,Casey,JiaLiuXu} for application to semilinear focusing wave equations. It turns out the outer energy \footnote{Other types of ``channels of energy"  have been found in other contexts, for instance in \cite{DKM3}, which are different from those ``channels of energy" in \cite{DKM,DKMsmall} for the associated linear wave equation. To distinguish, we call outer energy inequality specifically the inequality of the type in \cite{DKM,DKMsmall} for the linear wave equations.} inequality is strongest in odd dimensions, and in even dimensions only weaker versions hold, for special initial data. Using some weaker version of the outer energy inequality, the soliton resolution conjecture along a sequence of times was proved for $4+1$ dimensional focusing energy critical wave equation in \cite{Kenig4dWave}, and for $1-$ equivariant wave maps into the sphere in \cite{cotesoliton}. \\

However, for applications to six dimensional energy critical wave equation, $k-$ equivariant wave maps with $k\ge 2$, and $4+1$ dimensional radial Yang Mills equation, the outer energy inequality is not strong enough, consequently the soliton resolution along a sequence of times for them has been open until now. We remark that in \cite{Kenig4dWave,cotesoliton}, the outer energy inequality was only used in the last step, to show that the remainder term (which was already known to vanish in the Strichartz norm) vanishes in the energy norm. The difference between convergence in Strichartz norm and in energy norm is of course very important, especially in the presence of nontrivial solitons in the solution.\\ 

To explain the differences between the equations in \cite{Kenig4dWave,cotesoliton} and ours, let us consider for simplicity the energy critical semilinear wave equation. Suppose $u$ is a radial Type II solution that blows up at $(x,t)=(0,1)$. A key ingredient in \cite{Kenig4dWave} is the fact that there is asymptotically no energy in the self similar region $|x|\sim 1-t$ as $t\to 1-$, which by virial identities implies that the kinetic energy inside the backward lightcone vanishes asymptotically
\begin{equation}\label{eq:firstvirialintro}
\lim_{t\to1-}\frac{1}{1-t}\int_t^1\int_{|x|\leq 1-s}(\partial_su)^2\,dxds=0.
\end{equation}
From (\ref{eq:firstvirialintro}), using arguments in \cite{DKMsmall} one can deduce that along a sequence of times $t_n\to 1$, one has
\begin{equation}\label{eq:averageiszerointro}
\lim_{n\to\infty}\sup_{\tau\in(0,1-t_n)}\frac{1}{\tau}\int_{t_n}^{t_n+\tau}\int_{|x|\leq 1-t}(\partial_tu)^2\,dxdt=0.
\end{equation}
This is already sufficient to show that once we subtract the regular part in the profile decomposition of $\OR{u}(t_n)$ inside the lightcone $|x|\leq 1-t$, all the nontrivial profiles are given by solitons (in this case, the unique radial finite energy steady state, called ground state). Moreover the dispersive term is asymptotically of the form $\sim (w_{0,n},0)$. Up to this point no outer energy inequalities are needed. Then to show $w_{0,n}$ tends to zero in the energy norm, one has to appeal to outer energy inequalities. The argument goes as follows. Suppose this is not true, then the outer energy inequality developed in \cite{RKS} which in $4+1$ dimensions holds exactly for initial data of the form $(f,0)$, implies that there is nontrivial amount of radiation coming from $w_{0,n}$ that would propagate outside of the backward lightcone $|x|\leq 1-t$. This is a contradiction with the known fact that outside the lightcone $|x|\leq 1-t$ the energy of the solution was taken completely by the regular part of the solution. This argument breaks down in $6+1$ dimensions, as in this case the outer energy inequality in \cite{RKS} holds only for initial data of the form $(0,g)$. Similar difficulties were also met in the $k$ equivariant wave maps into the sphere with $k\ge 2$ and radial $4+1$ dimensional Yang Mills equations. \\

To overcome this difficulty, our main new idea is to use a second virial type inequality to show that $w_{0,n}$ vanishes in energy norm. In addition to (\ref{eq:firstvirialintro}), the lack of self similar blow up and virial identities also imply
\begin{equation}
\lim_{t\to1-}\frac{1}{1-t}\int_t^1\int_{|x|\leq 1-s}|\nabla u|^2-|u|^{\frac{2d}{d-2}}\,dxds=0,
\end{equation}
which we will see implies along some sequence of times $t'_n$
\begin{equation}\label{eq:secondvirialintro}
\lim_{n\to\infty}\sup_{\tau\in(0,1-t'_n)}\frac{1}{\tau}\int_{t'_n}^{t'_n+\tau}\int_{|x|\leq 1-t}|\nabla u|^2-|u|^{\frac{2d}{d-2}}\,\,dxdt\leq0.
\end{equation}
(\ref{eq:secondvirialintro}) is not coercive and hence does not seem to be particularly useful. However we observe that it becomes coercive if $u$ is dispersive in the sense that $\|u(\cdot,t)\|_{L^{\frac{2d}{d-2}}}$ becomes small in the limit. This is not quite true due to the presence of solitons, but for the solitons (\ref{eq:secondvirialintro}) holds trivially by a routine integration by parts argument. With the help of some pseudo-orthogonality properties of the profiles, we can indeed carry out this argument and conclude $w_{0,n}\to 0$ in the energy space, {\it provided} we can apply the profile decomposition for $t_n$ to the sequence $t_n'$. The same type of arguments apply to the wave map case as well, though in a more technical way (with the help of a new virial identity). The difficulty is thus that $t'_n$ need not be chosen as the same sequence $t_n$. One of our main new contributions is to show that we can find one sequence of times $t_n\to 1-$, that achieves both (\ref{eq:averageiszerointro}) and (\ref{eq:secondvirialintro}) simultaneously. There are similar arguments in the global existence case. \\

This approach completely avoids the use of outer energy inequalities (which although extremely powerful, in many cases are not available), and provides a unified approach to prove soliton resolution along a sequence of times for many critical wave type equations in the radial case, such as semilinear energy critical focusing wave equations in $R^d$, $d\ge 3$. For simplicity of presentation, we focuse on $6+1$ dimensional energy critical wave equation, general $k$ equivariant $2+1$ dimensional wave map equations, and $4+1$ dimensional radial Yang Mill equation, which seem unapproachable by outer energy inequalities. We also assume $k=1$ or $k=2$ below to get an easy local Cauchy theory and small data scattering result. The extension to higher equivariance classes is straightforward, with the help of the recent works \cite{xiaoyi,Casey} where the more sophisticated Cauchy theory for higher dimensional energy critical wave equation (involving fractional power nonlinearities) has been developed. Below we briefly state our main results, divided in two parts. The first part deals with the $6+1$ dimensional wave equation and the second part deals with wave map equations (in which the radial Yang Mills equation can be incorporated.).\\

In the first part we consider the radial energy critical focusing wave equation in $R^6\times[0,\infty)$
\begin{equation}\label{eq:mainwaveequation}
\left.\begin{array}{rl}
\partial_{tt}u-\Delta u-|u|u&=0\\
\OR{u}(0)&=(u_0,u_1)
\end{array}\right\}\quad (x,t)\in R^6\times [0,\infty).
\end{equation}
Equation (\ref{eq:mainwaveequation}) is invariant under the scaling
\begin{eqnarray}
u(x,t)&\to& u_{\lambda}(x,t)=\lambda^2u(\lambda x,\lambda t);\\
(u_0(x),u_1(x)) &\to& (u_{0\lambda}(x),u_{1\lambda}(x))=(\lambda^2u_0(\lambda x),\lambda^3u_1(\lambda x)),\label{eq:scalinginitial}
\end{eqnarray}
for any $\lambda>0$, in the sense that if $u(x,t)$ is a solution to equation (\ref{eq:mainwaveequation}) with initial data $(u_0,u_1)$ then $u_{\lambda}(x,t)$ is also a solution with initial data $(u_{0\lambda}(x),u_{1\lambda}(x))$. The initial data norm in $\HL(R^6)$ is left invariant under the natural scaling (\ref{eq:scalinginitial}) and this space is thus called a critical space. Equation (\ref{eq:mainwaveequation}) is locally well posed in the energy space. For any initial data $(u_0,u_1)\in\HL(R^6)$, there exists a unique solution $u\in C([0,T_+),\dot{H}^1(R^6))\cap L^2_tL^4_x(R^6\times I)$ for any time interval $I\Subset [0,T_+)$ (see \cite{KenigYang}). Here $T_{+}$ is the maximum time of existence for the solution $u$. Moreover, the energy
\begin{equation}
\E(\OR{u}(t))=\int_{R^6}\frac{(\partial_tu)^2}{2}+\frac{|\nabla u|^2}{2}-\frac{|u|^3}{3}(x,t)\,dx
\end{equation}
is conserved along the evolution for $t\in [0,T_{+})$. Equation (\ref{eq:mainwaveequation}) admits a unique (up to scaling and sign) radial finite energy steady state 
\begin{equation*}
W(x)=\frac{1}{(\frac{1}{24}+|x|^2)^2},
\end{equation*}
called the ground state, satisfying
\begin{equation*}
-\Delta W=W^2.
\end{equation*}
The ground state can be characterized as the minimizer for the best constant in the Sobolev embedding $\dot{H}^1(R^6)\hookrightarrow L^3(R^6)$ (see Talenti \cite{Talenti}), and plays a fundamental role in our analysis. In this paper we consider the extended type II solutions only, i.e., we assume
\begin{equation}
M:=\sup_{t\in[0,T_+)}\|\OR{u}(t)\|_{\HL(R^6)}<\infty, \,\,\,{\rm with}\,\,T_+<\infty,\,\,{\rm or}\,\,T_+=\infty.
\end{equation}
We show that the solution asymptotically decouples into a sum of rescaled ground states and a term vanishing in the limit, plus a regular part in the finite blow up case or a free radiation in the global case, at least along a sequence of times $t_n\to T_+$. More precisely we prove the following theorem.
\begin{theorem}\label{th:maintheoremsemi}
Let $u\in C([0,T_+),\dot{H}^1(R^6))\cap L^2_tL^4_x(R^6\times I)$ for any $I\Subset [0,T_+)$ be the solution to equation (\ref{eq:mainwaveequation}), with radial initial data $(u_0,u_1)\in\HL(R^6)$. Suppose $u$ is an extended Type II solution. Then there exist a sequence of times $t_n\to T_+$, an integer $J_0\ge 0$, signs $\{l_j\}$ with $l_j\in\{\pm 1\}$, and sequences $\lambda_{jn},\,\,j\leq J_0$ with 
\begin{equation}
0<\lambda_{1n}\ll \lambda_{2n}\ll\cdots\ll\lambda_{J_0n}, \footnote{\, For two sequences of positive numbers $a_n,\,b_n$, we write $a_n\ll b_n$ if $\lim\limits_{n\to\infty}\frac{a_n}{b_n}=0$.}\,\,{\rm as}\,\,t_n\to T_+,
\end{equation}
such that the following holds\\

\quad\quad (1) {\it (Global existence case)}\\
\begin{eqnarray*}
&&\lambda_{J_0n}\ll t_n, \,\,\,{\rm and}\\
&&\OR{u}(t_n)=\OR{u}^L(t_n)+\sum_{j=1}^{J_0}l_j\left(\frac{1}{\lambda_{jn}^2}W(\frac{x}{\lambda_{jn}}),0\right)+o_{{\tiny \HL(R^6)}}(1),\,\,{\rm as}\,\,n\to\infty,
\end{eqnarray*}
where the function $\OR{u}^L$ is a radial finite energy solution to the linear wave equation (the free radiation).\\

\quad\quad (2) {\it (Finite time blow up case)}\\
\begin{eqnarray*}
&&J_0\ge 1;\\
&&\lambda_{J_0n}\ll T_+-t_n, \,\,\,{\rm and}\\
&&\OR{u}(t_n)=\OR{v}+\sum_{j=1}^{J_0}l_j\left(\frac{1}{\lambda_{jn}^2}W(\frac{x}{\lambda_{jn}}),0\right)+o_{{\tiny \HL(R^6)}}(1),\,\,{\rm as}\,\,n\to\infty,
\end{eqnarray*}
where $\OR{v}\in \HL(R^6)$.
\end{theorem}

\smallskip
\noindent
{\it Remark.} The theorem implies quantization of energy for $\OR{u}(t)$ modulo free radiation as $t\to\infty$ in the global case, and inside the backward lightcone $|x|\leq T_+-t$ as $t\to T_+$ in the finite time blow up case. The question of whether the soliton resolution holds for all times $t\uparrow T_+$, instead of just along a sequence, is still open. However, if we assume an {\it additional} size restriction, such as 
\begin{equation}
\sup_{t<T_+}\left(\|u(t)\|^2_{\dot{H}^1(r<T_+-t)}+\|\partial_tu(t)\|^2_{L^2(r<T_+-t)}\right)<2\|W\|^2_{\dot{H}^1(R^6)},
\end{equation}
in the Type II blow up case, then $J_0\leq 1$ and the soliton resolution holds for all sequences $t_n\uparrow T_+$.  A similar statement holds for $T_+=\infty.$ See \cite{Kenig4dWave} for more details.

\bigskip
\noindent
In the second part, we consider the equation
\begin{equation}\label{eq:WP}
\partial_{tt}\psi-\partial_{rr}\psi-\frac{\partial_r\psi}{r}+\frac{f(\psi)}{r^2}=0,
\end{equation}
where $f(\rho)=g(\rho)g'(\rho)$ and $\partial_{rr}+\frac{\partial_r}{r}$ is the two dimensional radial Laplacian. We make the following assumptions on $g$.
\begin{eqnarray*}
(A1)&& G(x):=\int_0^x|g(y)|dy\to\infty, \,\,\,{\rm as}\,\,|x|\to\infty,\\
(A2)&& g\in C^3,\,\,\,{\rm define}\,\,\,\mathcal{V}:=\{l\in R:\,g(l)=0\}, \\
&&{\rm then\,\,}\mathcal{V}\,\,{\rm \,\,is\,\,discrete\,\,and\,\,the\,\,cardinality\,\,of\,\,}\mathcal{V},|\mathcal{V}|\ge 2,\\
(A3)&&\forall l\in \mathcal{V}, \,\,g'(l)\in\{-2,-1,1,2\},\,\,{\rm and\,\,if\,\,}g'(l)=\pm1,\,\,{\rm then}\,\,g''(l)=0.
\end{eqnarray*}
Assumption (A1) ensures that finite energy solutions are bounded, (A2) implies that equation (\ref{eq:WP}) admits non-constant steady states (harmonic maps), (A3) is needed for an easy local Cauchy theory and small data scattering results. By the work \cite{xiaoyi,Casey} we can relax assumption (A3) to
\begin{equation*}
(A3)'\quad\,\forall l\in \mathcal{V}, \,\,|g'(l)|\,\,{\rm is\,\,an\,\, integer}\ge 1,\,\,{\rm and\,\,if\,\,}g'(l)=\pm1,\,\,{\rm then}\,\,g''(l)=0.
\end{equation*}
This equation arises in the study of equivariant wave maps into spheres and the radial four dimensional Yang Mills equation. 

The general wave map $U$ from $R^{d+1}$ with Minkowski metric to a Riemannian manifold $M$ with metric $h$ takes the form (in local coordinates)
\begin{eqnarray*}
\partial_{tt}U^a-\Delta U^a&=&\Gamma^a_{b,c}\partial_{\mu}U^b\partial^{\mu}U^c,\\
\OR{U}(0)&=&(U_0,U_1),
\end{eqnarray*}
where $\Gamma^a_{b,c}$ are the Christoffel symbols for the metric $h$. Suppose $M$ is a surface of revolution with polar coordinates $(\rho,\theta)$ and metric $d\rho^2+h^2(\rho)d\theta^2$. Assume the wave map $U$ is $k$ equivariant, given by $(r,\omega,t)\to (\rho,\theta)=(\psi(r,t),k\omega)$, we then obtain equation (\ref{eq:WP}) with $g=kh(\rho)$. In the special case the target manifold being the unit sphere $S^2\subset R^3$, we have
\begin{equation*}
f(\psi)=k^2\sin{\psi}\,\cos{\psi},\,\,\,{\rm with}\,\,g(\psi)=k\sin{\psi},
\end{equation*}
and the solitons (harmonic maps) are given by
\begin{equation*}
Q(r)=2\arctan{\left(r^k\right)}+\ell \pi,\,\,\ell\in Z.
\end{equation*}
We refer the readers to \cite{Struwe2,cotesoliton} for further discussion on the background of this equation.

In the case the radial $4+1$ dimensional Yang Mills system for the $so(4)$-valued gauge potential $A_{\alpha}$ and curvature $F_{\alpha\beta}$
\begin{eqnarray*}
&&F_{\alpha\beta}=\partial_{\alpha}A_{\beta}-\partial_{\beta}A_{\alpha}+[A_{\alpha},\,A_{\beta}],\\
&&\partial_{\beta}F^{\alpha\beta}+[A_{\beta},F^{\alpha\beta}]=0,\,\,\alpha,\,\beta=0,\dots,3,
\end{eqnarray*}
if we make the equivariant reduction
\begin{equation*}
A_{\alpha}^{ij}=\left(\delta^i_{\alpha}x^j-\delta^j_{\alpha}x^i\right)\frac{1-\psi(r,t)}{r^2},
\end{equation*}
we arrive at equation (\ref{eq:WP}) with 
\begin{equation*}
f(\psi)=-2\psi(1-\psi^2),\,\,\,{\rm with}\,\,g(\psi)=\frac{1}{2}(1-\psi^2).
\end{equation*}
The solitons (instantons) are given by 
\begin{equation*}
Q(r)=\frac{1-r^2}{1+r^2}\,\,\,{\rm or}\,\,\,\frac{r^2-1}{1+r^2}\,.
\end{equation*}
We refer to \cite{RapRod} for more discussion of the physical background of this equation.\\

Define the energy and the Hilbert space $H\times L^2$ as follows: for $\OR{\phi}=(\phi_0,\phi_1)$ and $0\leq r_1<r_2\leq \infty$, 
\begin{eqnarray*}
&&\E(\OR{\phi},r_1,r_2):=\int_{r_1}^{r_2}\left(\phi_1^2+(\partial_r\phi_0)^2+\frac{g^2(\phi_0)}{r^2}\right)\,rdr,\\
&&\|\phi_0\|^2_{H([r_1,r_2])}:=\int_{r_1}^{r_2}\left((\partial_r\phi_0)^2+\frac{\phi_0^2}{r^2}\right)\,rdr,\\
&&\|\OR{\phi}\|_{H\times L^2([r_1,r_2])}=\int_{r_1}^{r_2}\left(\phi_1^2+(\partial_r\phi_0)^2+\frac{\phi_0^2}{r^2}\right)\,rdr.
\end{eqnarray*}
By \cite{Shatah2} we know that for each $(\psi_0,\psi_1)$ with $\E(\psi_0,\psi_1):=\int_0^{\infty}\left(\psi_t^2+\psi_r^2+\frac{g^2(\psi)}{r^2}\right)\,rdr<\infty$, there exists a unique solution $\OR{\psi}\in C(I,H\times L^2)$, defined on a maximal interval $I:=[0,T_+)$ which preserves the energy, and satisfies $\psi_0(0,t)\equiv \psi_0(0)$, $\psi(\infty,t)=\psi_0(\infty)$ for each $t\in [0,T_+)$. Moreover, we know that
\begin{equation*}
\frac{f(\psi)}{r^{\frac{2k}{k+2}}}\in L^{\frac{k+2}{k}}_tL^{\frac{2k+4}{k}}_x, 
\end{equation*}
locally, where $k=|g'(\psi_0(0))|$.
Equation (\ref{eq:WP}) admits nontrivial steady states $Q$ of finite energy,
\begin{equation*}
\partial_{rr}Q+\frac{\partial_rQ}{r}=\frac{f(Q)}{r^2}.
\end{equation*}
Thanks to \cite{Cote},  we know $Q$ is monotone, and satisfies either $r\partial_rQ=g(Q)$ or $r\partial_rQ=-g(Q)$. Moreover $\{Q(0),\,Q(\infty)\}=\{l,\,m\}$ with some $l,\,m\in \mathcal{V}$, $l<m$, satisfying $\mathcal{V}\cap (l,m)=\emptyset$. $\E(Q,0)=2(G(m)-G(l))$. To study large solutions, we also need the linearized wave equation around any $l\in\mathcal{V}$:
\begin{equation}\label{eq:linearizedequationl}
\begin{array}{rc}
({\rm LW}_{\ell})&\partial_{tt}\phi-\partial_{rr}\phi-\frac{\partial_r\phi}{r}+\frac{g'(l)^2}{r^2}\phi=0.\\
\end{array}
\end{equation}
Direct calculation shows that if $\phi$ is a solution to $({\rm LW}_{\ell})$, then  $u_L=r^{-|g'(l)|}\phi$ verifies the $2|g'(l)|+2$ dimensional linear wave equation. Hence solutions to $({\rm LW}_{\ell})$ preserve the following energy
\begin{equation}
\|\OR{\phi}(t)\|^2_{\mathcal{H}_{\ell}\times L^2}:=\int_0^{\infty}\left(\phi_t^2+\phi_r^2+\frac{g'(l)^2\phi^2}{r^2}\right)\,rdr=\|\OR{\phi}(0)\|^2_{\mathcal{H}_{\ell}\times L^2}.
\end{equation}
Our main goal is to prove the following result. 
\begin{theorem}\label{th:JiaKenig}
Let $\OR{\psi}(t)$ be a finite energy solution to equation (\ref{eq:WP}). Then there exists a sequence of times $t_n\uparrow T_+$, an integer $J\ge 0$, $J$ sequences of scales $0<r_{Jn}\ll \cdots\ll r_{2n}\ll r_{1n}$ and $J$ harmonic maps $Q_1,\cdots,Q_J$, such that
\begin{equation*}
Q_J(0)=\psi_0(0),\,\,\,Q_{j+1}(\infty)=Q_j(0),\,\,\,{\rm for}\,\,j=1,\cdots,J-1,
\end{equation*}
and such that the following holds.\\
\quad\quad (1)\,\, (Global case) If $T_+=\infty$, denote $l=\psi(\infty)$. Then $Q_1(\infty)=l$, $r_{1n}\ll t_n$ and there exists a radial finite energy solution $\OR{\phi}_L$ to the linear wave equation $(LW_{\ell})$ with $l=\psi_0(\infty)$, such that 
\begin{equation}
\OR{\psi}(t_n)=\sum_{j=1}^J(Q_j(\cdot/r_{jn})-Q_j(\infty),0)+(l,0)+\OR{\phi}_L(t_n)+\OR{b}_n.
\end{equation}
(2)\,\,(Blow up case) If $T_+<\infty$, denote $l=\lim_{t\to T_+}\psi(T_+-t,t)$ (it is well defined). Then $J\ge 1$, $r_{1,n}\ll T_+-t_n$ and there exists a function $\OR{\phi}\in H\times L^2$ of finite energy, such that $Q_1(\infty)=\phi(0)=l$ and 
\begin{equation}
\OR{\psi}(t_n)=\sum_{j=1}^J(Q_j(\cdot/r_{jn})-Q_j(\infty),0)+\OR{\phi}+\OR{b}_n.
\end{equation}
In both cases, $\OR{b}_n=(b_{0,n},b_{1,n})$ vanishes in the following sense
\begin{equation}
\|\OR{b}_n\|_{H\times L^2}\to 0,\,\,{\rm as}\,\,n\to\infty.
\end{equation}
\end{theorem}

\smallskip
\noindent
{\it Remark.} This theorem has been proved by C\^ote\cite{cotesoliton} under the assumptions (A1),(A2) and that $g'(l)\in\{-1,1\}$ for each $l\in\mathcal{V}$, which does not include the cases of $k$ equivariant wave maps for $k\ge 2$ or the radial Yang Mills equation in $R^4$. The theorem implies quantization of energy inside the lightcone $r\leq T_+-t$ in the finite blow up case and after subtracting the radiation term in the global existence case. The question of whether the soliton resolution holds for all times $t\uparrow T_+$, instead of just along a sequence, remains open. However, if one assumes that the topological constraint requires at least $J$ solitons inside the lightcone $|x|\leq T_+-t$, (for example, assuming $Q_j(\infty)$ is monotone in $j$, so that connecting the value $l=\lim_{t\to T_+}\psi(T_+-t,t)$ to the value $\psi(0,t)$ would require at least the energy of the $J$ solitons $Q_j$ while the backward lightcone has only enough energy to contain these $J$ solitons in the limit $t\to T_+$), then the decomposition holds for all time sequences $t_n\uparrow T_+$. We refer the reader to \cite{cotesoliton} for more discussions.

\bigskip

Our paper is organized as follows:  in Section 2  we treat the radial energy critical semilinear wave equation (\ref{eq:mainwaveequation}), in Section 3 we consider the equivariant wave map equation (\ref{eq:WP}), in Section 4 we prove the crucial real analysis lemma needed to choose the sequences of times $t_n$ to achieve (\ref{eq:averageiszerointro}) and (\ref{eq:secondvirialintro}) simultaneously, in the Appendix we review various useful facts on the Cauchy theory for (\ref{eq:mainwaveequation}) and (\ref{eq:WP}), and give an alternative treatment for the profile decomposition adapted to the equivariant wave map equation, introduced by C\^ote\cite{cotesoliton}.
\end{section}

\begin{section}{Part I: Asymptotics for six dimensional radial energy critical semilinear wave equation}
In this section we prove Theorem \ref{th:maintheoremsemi}. In the first subsection we consider the finite time blow up case, and in the second subsection we consider the global existence case. The proofs are based on similar strategies, with some notable differences. 
\begin{subsection}{Asymptotics for $6d$ wave, blow up case}\label{subsec:6dwaveblowup}
Let $u\in C([0,T_+),\dot{H}^1(R^6))\cap L^2_tL^4_x(R^6\times [0,T))$ for any $T<T_+$ be the unique solution to equation (\ref{eq:mainwaveequation}) with radial initial data $(u_0,u_1)\in \HL(R^6)$. By rescaling, without loss of generality, we assume $\OR{u}$ is a Type II blow up solution which blows up at $(x,t)=(0,1)$ with $T_+=1$. Let us denote
\begin{equation}
M:=\sup_{t\in [0,1)}\|\OR{u}(t)\|_{\HL(R^6)}<\infty.
\end{equation}
With slight abuse of notation, we also write $u(r,t)=u(x,t)$ with $|x|=r$. By Cauchy Schwartz inequality we have 
\begin{equation}\label{eq:radialestimate}
|u|(r,t)=\int_r^{\infty}|\partial_{\tau}u(\tau,t)|\,\tau^{\frac{5}{2}}\,\tau^{-\frac{5}{2}}\,d\tau\leq Cr^{-2}\left(\int_{R^6}|\nabla u|^2\,dx\right)^{\frac{1}{2}}.
\end{equation}
Thus $u(x,t)$ is bounded away from the origin $x=0$ when $t\to 1-$. Thus from the assumption that $u$ is Type II, Sobolev inequality and interpolation inequalities, we see $u\in L^2_tL^4_x((R^6\backslash B_r(0))\times [0,1))$ for any $r>0$. Then by the local Cauchy theory and finite speed of propagation we can extend the solution up to time $t=1$ except at the origin and obtain $\nabla_{x,t}u\in C([0,1],L^2(R^6\backslash B_r(0)))$ for any $r>0$. Thus $\OR{u}(1)$ is well defined as a function in $\HL(R^6)$. See also \cite{KenigYang}. We emphasize that we can only say $u\in C([0,1],\dot{H}^1(R^6\backslash B_r(0)))$ for $r>0$, and in a small ball around the origin, this continuity does not hold. Let $\OR{v}$ be the unique solution to equation (\ref{eq:mainwaveequation}) with $\OR{v}(1)=\OR{u}(1)$. We call $v$ the regular part of $u$. By the local Cauchy theory for equation (\ref{eq:mainwaveequation}), $\OR{v}$ is defined at least in a short time interval $[1-\delta,1+\delta]$ for some small $\delta>0$. By a similar estimate to (\ref{eq:radialestimate}) we see
\begin{equation}
\sup_{t\in[t-\delta,1],\,r\in (0,\infty)}r^2|v(r,t)|<\infty.
\end{equation} 
Thus using approximation by smooth functions in $\dot{H}^1$ and the fact that $v\in C([1-\delta,1],\dot{H}^1)$ we see that
\begin{equation}
\lim\limits_{r\to0+}\sup_{t\in[1-\delta,1]}r^2|v(r,t)|=0.
\end{equation}
Since $v\in C([1-\delta,1],\dot{H}^1(R^6\backslash B_r(0)))$ and $\OR{v}(1)=\OR{u}(1)$, by the local Cauchy theory and finite speed of propagation we see that
\begin{equation}
u(x,t)=v(x,t)\,\,\,{\rm for}\,\, |x|\ge r+|t-1|, \,\,t\in [1-\delta,1)\,\,\, {\rm and \,\,\,any}\,\, r>0.
\end{equation}
Thus $u(x,t)=v(x,t)$ for all $(x,t)$ satisfying $|x|\ge |1-t|$ and $t\in [1-\delta,1)$. Hence we have
\begin{equation}\label{eq:tendtozero}
\lim_{t\to 1-}(1-t)^2|u(1-t,t)|=\lim_{t\to 1-}(1-t)^2|v(1-t,t)|=0.
\end{equation}
Hardy's inequality then gives
\begin{equation}\label{eq:hardy}
\def\arraystretch{2.2}
\left.\begin{array}{rl}
            &\limsup\limits_{t\to 1-}\int_{|x|\leq 1-t} \frac{|u|^2(x,t)}{|x|^2}\,dx\\
              &\leq C\,\limsup\limits_{t\to 1-}\left(\int_{|x|\leq 1-t}|\nabla u|^2(x,t)\,dx+\frac{1}{1-t}\int_{|x|=1-t}|u|^2(x,t)\,d\sigma(x)\right)\\
              &\leq CM^2.
         \end{array}\right.
\end{equation}
Let us now estimate $u$ in the interior of the lightcone $|x|<1-t$ for $t\in [1-\delta,1)$. Using
\begin{equation}
r^2u(r,t)=(1-t)^2u(1-t,t)-\int_r^{1-t}\partial_{\tau}\left(\tau^2u(\tau,t)\right)\,d\tau,
\end{equation}
by H\"older inequality, we have
\begin{equation}\label{eq:smallnessnearlightcone1}
r^2|u|(r,t)\leq (1-t)^2|u(1-t,t)|+C\log{\left(\frac{1-t}{r}\right)}\,\left(\int_{r\leq |x|\leq 1-t}|\nabla u|^2+\frac{|u|^2}{|x|^2}(x,t)\,dx\right)^{\frac{1}{2}}.
\end{equation}
Thus if we choose $\lambda_0=\lambda_0(M)<1$ sufficiently close to $1$, then by (\ref{eq:tendtozero}) and (\ref{eq:hardy}) we have 
\begin{equation}\label{eq:smallnessnearlightcone}
\limsup_{t\to 1-}\sup_{1-t\ge r\ge \lambda_0(1-t)}r^2|u(r,t)|\leq CM\log{\frac{1}{\lambda_0}}<\frac{1}{10}.
\end{equation}
This shows that near the boundary of lightcone $|x|=1-t$, $u$ is small in a scale invariant way.\\

Our  first main result in this section is the following theorem, which roughly speaking asserts that asymptotically there is no energy concentration in the self similar region $|x|\sim 1-t$, as $t\to 1-$. Thus energy concentraction in the lightcone $|x|<1-t$ as $t\to 1-$ happens in a much smaller scale. This shows that Type II blow up is quite ``slow". 
\begin{theorem}\label{th:selfsimilar}
Let $u\in C([0,1),\dot{H}^1(R^6))\cap L^2_tL^4_x(R^6\times I)$ for any $I\Subset [0,1)$ be the solution to equation (\ref{eq:mainwaveequation}), with radial initial data $(u_0,u_1)\in\HL(R^6)$. Suppose that $u$ is a Type II blow up solution with blow up point $(x,t)=(0,1)$. Then for any $\lambda\in (0,1)$ we have
\begin{equation}
\lim_{t\to 1-}\int_{\lambda(1-t)\leq |x|\leq 1-t}|\nabla_{x,t}u|^2(x,t)+\frac{u^2}{|x|^2}(x,t)\,dx=0.
\end{equation}
\end{theorem}

\smallskip
\noindent
The proof of this theorem is similar to the case of four dimensional energy critical focusing wave equation \cite{Kenig4dWave} and we follow closely their arguments (which rely on techniques introduced by Christodoulou, Tahvildar-Zadeh, \cite{Chri,ChriTah}, and Shatah, Tahvildar-Zadeh, \cite{Shatah1,Shatah2}), with a slight modification which allows us to eliminate the assumption that the solution is smooth before the blow up time. The key observation in \cite{Kenig4dWave}, in the context of $R^6$, is to study instead $\psi(r,t)=r^2u(r,t)\in C(R^+\times [0,1]\backslash\{(0,1)\})$, which satisfies
\begin{equation}\label{eq:equationforpsi}
\partial_{tt}\psi-\partial_{rr}\psi-\frac{\partial_r\psi}{r}+\frac{4\psi-|\psi|\psi}{r^2}=0,
\end{equation}
and
\begin{equation}
\sup_{t\in (0,1)}\int_0^{\infty}\left(\frac{\psi_t^2}{2}+\frac{\psi_r^2}{2}+\frac{\psi^2}{r^2}\right)rdr<\infty.
\end{equation}
The main advantage of switching to $\psi$ is that when $\psi$ is small, 
\begin{equation*}
\frac{2\psi^2-|\psi|^3/3}{r^2}>0,
\end{equation*}
which implies non-negativity of the energy flux on the boundary of lightcone $\{(r,t)|\,r\leq 1-t\}$. Define
\begin{eqnarray}
&&f(\psi)=4\psi-|\psi|\psi,\\
&&F(\psi)=\int_0^{\psi}f(\rho)\,d\rho=2\psi^2-|\psi|^3/3.
\end{eqnarray}
Let us also set
\begin{eqnarray}
&&e(r,t)=\frac{(\partial_r\psi)^2}{2}+\frac{(\partial_t\psi)^2}{2}+\frac{F(\psi)}{r^2};\\
&&m(r,t)=\partial_t\psi\,\partial_r\psi;\\
&&L(r,t)=-\frac{(\partial_t\psi)^2}{2}+\frac{(\partial_r\psi)^2}{2}+\frac{F(\psi)}{r^2}-\frac{2f(\psi)\partial_r\psi}{r};\\
&&\Fl(t_0;t_1)=\int_{t_0}^{t_1}(1-l)e(1-l,l)-(1-l)m(1-l,l)\,dl.
\end{eqnarray}
It's easy to check the energy
\begin{equation}
\int_0^{\infty}r\,e(r,t)\,dr
\end{equation}
is conserved for equation (\ref{eq:WP}) along the evolution. 

We prove Theorem \ref{th:selfsimilar} as an immediate consequence of the following lemma.
\begin{lemma}\label{lm:noenergyforpsi}
Let $\lambda\in(0,1)$, then
\begin{equation}
\E^{\lambda}_{{\rm ext}}(t):=\int_{\lambda (1-t)}^{1-t}\left(\frac{\psi_t^2}{2}+\frac{\psi_r^2}{2}+\frac{\psi^2}{r^2}\right)(r,t)\,r\,dr\to 0\,\,{\rm as}\,\,t\to 1-.
\end{equation}
\end{lemma}

\smallskip
\noindent
{\it Proof.} We prove this lemma holds for each $\lambda=\lambda_0^j$, $j\ge 1$, by induction. Here $\lambda_0$ is taken from (\ref{eq:smallnessnearlightcone}). For $\tau\in (0,1)$, denote
\begin{equation}
K^{\lambda}_{{\rm ext}}(\tau):=\{(r,t):\,1-t\ge r\ge \lambda (1-t),\,t\in (\tau,1)\}.
\end{equation}
Let us first prove the case $j=1$. By estimate (\ref{eq:smallnessnearlightcone}) there exists $\tau<1$ sufficiently close to $1$ such that 
\begin{equation*}
|\psi(r,t)|\leq \frac{1}{10}\quad {\rm for}\,\,(r,t)\in K^{\lambda}_{{\rm ext}}(\tau).
\end{equation*}
Hence $F(\psi)>0$ on $K^{\lambda}_{{\rm ext}}(\tau)$, and consequently $re\pm rm\ge 0$. Then it's meaningful to define  
\begin{eqnarray}
&&\A^2(r,t)=re+rm;\label{eq:defforAI}\\
&&\B^2(r,t)=re-rm.\label{eq:defforBI}
\end{eqnarray}
Elementary calculations show (in these calculations we can assume $\psi$ to be smooth by suitable approximation arguments)
\begin{eqnarray}
&&\partial_t(re)-\partial_r(rm)=0;\label{eq:someidentities1}\\
&&\partial_t(rm)-\partial_r(re)=L.\label{eq:someidentities2}
\end{eqnarray}
Thus
\begin{eqnarray}
&&\partial_t\A^2-\partial_r\A^2=L;\label{eq:someidentities3}\\
&&\partial_t\B^2+\partial_r\B^2=-L.\label{eq:someidentities4}
\end{eqnarray}

It's more convenient to work in the null coordinates 
\begin{equation}
\xi=t-1-r,\,{\rm and}\,\,\,\eta=t-1+r.
\end{equation}
Then from (\ref{eq:someidentities3}), (\ref{eq:someidentities4}) we see
\begin{equation}\label{eq:nullderivative}
\def\arraystretch{2.2}
\left.\begin{array}{rl}
&\partial_{\xi}\A^2=L/2;\\
&\partial_{\eta}\B^2=-L/2.
\end{array}\right.
\end{equation}
Note that both $\B^2\leq Cr e$ and $|\partial_{\eta}\B^2|\leq |L/2|\leq C e$ are integrable in $t,\,r$ locally away from $r=0$. Hence, by Lemma \ref{lm:restriction}, $\B^2$ can be restricted to $\eta={\rm constant}$, that is, $t+r={\rm constant}$, at least away from $r=0$. Consequently ${\rm Flux}(t_0,t_1)$ is well defined. Moreover by a suitable mollification argument, we have the following local energy flux identity
\begin{equation}
\int_0^{1-t_0}re(r,t_0)\,dr=\int_0^{1-t_1}re(r,t_1)\,dr+\Fl(t_0;t_1)\,\,{\rm for\,\,}0<t_0<t_1<1.
\end{equation}
For $(r,t)\in K^{\lambda}_{{\rm ext}}(\tau)$, $F(\psi)>0$, and thus 
\begin{equation}
 e-m=\frac{1}{2}(\partial_r\psi-\partial_t\psi)^2+\frac{F(\psi)}{r^2}>0.
\end{equation}
It follows that
\begin{equation*}
\int_0^{t_1}re(r,t_1)\,dr
\end{equation*}
 is decreasing and tends to some finite limit as $t_1\to 1-$. Thus
\begin{equation}
\Fl(t_0):=\lim_{t_1\to 1-} \Fl(t_0;t_1)
\end{equation} 
is well defined. Furthermore 
\begin{equation}\label{eq:fluxvanishing}
\lim_{t_0\to 1-}\Fl(t_0)=0.
\end{equation}

We next claim the bound 
\begin{equation}\label{eq:boundL}
L^2\leq C\frac{\A^2\B^2}{r^2},\quad\,\,   {\rm on }\,\,\, K^{\lambda}_{{\rm ext}}(\tau)
\end{equation}
The proof of this claim is exactly the same as in Claim 3.8 of \cite{Kenig4dWave}. We sketch some of the details for completeness. A direct computation yields
\begin{equation}\label{eq:estimateofL}
L^2\leq \frac{1}{2}(\psi_r^2-\psi_t^2)^2+\frac{4}{r^4}F^2(\psi)+\frac{16}{r^2}\psi_r^2f^2(\psi).
\end{equation}
By the bound of $\psi$ on $K^{\lambda}_{{\rm ext}}(\tau)$, we have
\begin{eqnarray*}
&&|f(\psi)|\leq 5|\psi|,\\
&&F(\psi)\ge |\psi|^2.
\end{eqnarray*}
The above inequalities imply
\begin{equation}\label{eq:estimateoff}
|f(\psi)|^2\leq 25F(\psi),\,\,\,{\rm on}\,\,\,K_{{\rm ext}}^{\lambda}(\tau).
\end{equation}
Plugging (\ref{eq:estimateoff}) into (\ref{eq:estimateofL}) we obtain
\begin{equation}
\def\arraystretch{2.2}
\left.\begin{array}{rl}
            L^2 &\leq \frac{1}{2}(\psi_r^2-\psi_t^2)^2+\frac{4}{r^4}F^2(\psi)+\frac{16}{r^2}\psi_r^2f^2(\psi)\\
          & \leq \,C\left(\frac{1}{4}(\psi_r^2-\psi_t^2)^2+\frac{1}{r^4}F^2(\psi)+\frac{1}{r^2}F(\psi)(\psi_r^2+\psi_t^2)\right).
               \end{array}\right.
\end{equation}
On the other hand, 
\begin{equation}
\frac{\A^2\B^2}{r^2}=\frac{1}{4}(\psi_r^2-\psi_t^2)^2+\frac{1}{r^4}F^2(\psi)+\frac{1}{r^2}F(\psi)(\psi_r^2+\psi_t^2),
\end{equation}
which together with the preceding inequality, establishes our claimed inequality (\ref{eq:boundL}). \\
Now we can combine (\ref{eq:nullderivative}) and (\ref{eq:boundL}) to get
\begin{equation}\label{eq:nullderivativeestimate}
|\partial_{\xi}\A|\leq \frac{C}{r}\B,\quad\,\,|\partial_{\eta}\B|\leq \frac{C}{r}\A,\,\,\,{\rm on}\,\,K^{\lambda}_{{\rm ext}}(\tau).
\end{equation}
Fix $\xi_0$ such that $(0,\xi_0)\in K_{{\rm ext}}^{\lambda}(\tau)$ and consider the rectangle
\begin{equation}\label{eq:geometric}
\Gamma(\eta,\xi):=[\eta,0]\times[\xi_0,\xi]\subset K^{\lambda}_{{\rm ext}}(\tau).
\end{equation}
Then elementary geometry shows that on $K^{\lambda}_{{\rm ext}}(\tau)$ there exists $C(\lambda)>1$ such that
\begin{equation}\label{eq:geom}
\frac{|\xi|}{C(\lambda)}\leq r\leq |\xi|,\, \,\,\,{\rm and}\,\,\,\,|\eta|\leq |\xi|.
\end{equation}
Viewing $\A,\,\B$ as functions of $\xi,\,\eta$, combining (\ref{eq:nullderivativeestimate}) and (\ref{eq:geom}), we obtain
\begin{equation}\label{eq:controlofflux}
\def\arraystretch{2.2}
\left.\begin{array}{rl}
    \A(\eta,\xi)&= \A(\eta,\xi_0)+\int_{\xi_0}^{\xi}\partial_u\A(\eta,u)\,du\\
                    &\leq A(\eta,\xi_0)+C(\lambda)\int_{\xi_0}^{\xi}\frac{B(\eta,u)}{|u|}\,du\\
                    &= A(\eta,\xi_0)+C(\lambda)\int_{\xi_0}^{\xi}\frac{B(0,u)}{|u|}\,du-C(\lambda)\int_{\xi_0}^{\xi}\int_{\eta}^0\frac{\partial_v\B(v,u)}{|u|}\,dudv\\
                   &\leq A(\eta,\xi_0)+C(\lambda)\int_{\xi_0}^{\xi}\frac{B(0,u)}{|u|}\,du+C(\lambda)\int_{\xi_0}^{\xi}\int_{\eta}^0\frac{\A(v,u)}{u^2}\,dudv\\
       \end{array}\right.\\
\end{equation}
In the above we adopt the convention that the constant $C(\lambda)$ may change from line to line. We can estimate the 2nd term in the last line of (\ref{eq:controlofflux}) as follows
\begin{equation*}
\def\arraystretch{2.2}
\left.\begin{array}{rl}
\int_{\xi_0}^{\xi}\frac{B(0,u)}{|u|}du&\leq \left(\int_{\xi_0}^{\xi}\B^2(0,u)du\right)^{\frac{1}{2}}\left(\int_{\xi_0}^{\xi}\frac{1}{u^2}du\right)^{\frac{1}{2}}\\
                 &\leq C \sqrt{\frac{\Fl(\xi_0)}{|\xi|}}.
                  \end{array}\right.
\end{equation*}
Take $(\eta_1,\xi_1)$ with $\Gamma(\eta_1,\xi_1)\subset  K^{\lambda}_{{\rm ext}}(\tau)$ and $\xi_1\ge\xi_0$. For $(\eta,\xi)\in \Gamma(\eta_1,\xi_1)$, taking $L^2$ norm of (\ref{eq:controlofflux}) in $\eta$, applying Minkowski inequality and noting that (\ref{eq:geom}) holds on $\Gamma(\eta_1,\xi_1)$, we obtain
\begin{equation}
\def\arraystretch{2.2}
\left.\begin{array}{rl}
     \|\A(\cdot,\xi)\|_{L^2(\eta_1,0)}&\leq \|\A(\cdot,\xi_0)\|_{L^2(\eta_1,0)}+C(\lambda)\,\sqrt{\Fl(\xi_0)}\,+\\
          &\quad\quad\quad\,\,+\,C(\lambda)\int_{\xi_0}^{\xi}\frac{1}{u^2}\,\|\int_{\eta_1}^0\A(v,u)\chi_{(\eta,0)}(v)dv\|_{L^2_{\eta}(\eta_1,0)}\,du\\
                         &{\rm by\,\,Schwartz\,\,inequality}\\
                                     &\leq \|\A(\cdot,\xi_0)\|_{L^2(\eta_1,0)}+C(\lambda)\,\sqrt{\Fl(\xi_0)}+C(\lambda)\int_{\xi_0}^{\xi}\frac{|\eta_1|}{u^2}\,\|\A(\cdot,u)\|_{L^2(\eta_1,0)}\,du.
           \end{array}\right.
\end{equation}


Note that $\xi\leq \xi_1<0$, thus $|\eta_1|\leq |\xi_1|\leq |\xi|$. Then
\begin{equation*}
 \int_{\xi_0}^{\xi}\frac{|\eta_1|}{u^2}\,du\leq \int_{-\infty}^{\xi}\frac{|\eta_1|}{u^2}\,du\leq 1.
\end{equation*}
By Gronwall's inequality, we see that for $\forall\,\xi_0\leq \xi\leq \xi_1$
\begin{equation}
\|\A(\cdot,\xi)\|_{L^2(\eta_1,0)}\leq C(\lambda)\left(\|\A(\cdot,\xi_0)\|_{L^2(\eta_1,0)}+\sqrt{\Fl(\xi_0)}\right).
\end{equation}
By the geometric constraint $\Gamma(\eta_1,\xi_1)\subset  K^{\lambda}_{{\rm ext}}(\tau)$, we have $\frac{1-\lambda}{1+\lambda}\xi_1\leq\eta_1\leq 0$. Hence
\begin{equation}
\lim_{\xi_1\to 0-}\eta_1=0.
\end{equation}
We claim 
\begin{equation}
\lim_{\xi_1\to 0-}\|\A(\cdot,\xi_0)\|_{L^2(\eta_1,0)}=0.
\end{equation}
Indeed, noting that $\A^2\leq Cre$ is integrable in $t,\,r$ (and hence in $\xi,\,\eta$) in a neigborhood of $(\eta,\xi)=(0,\xi_0)$. Moreover,
\begin{equation}
\left|\partial_{\xi}\A^2\right|\leq |L|\leq C e
\end{equation}
is also integrable (in that neighborhood). Thus $\A^2|_{\xi=\xi_0}$ is well defined as an integrable function for $|\eta|$ sufficiently small. The claim then follows. 

Combining the above claim with (\ref{eq:fluxvanishing}) we obtain
\begin{equation}
\lim_{\xi\to0}\|\A(\cdot,\xi)\|_{L^2\left(\frac{1-\lambda}{1+\lambda}\xi,0\right)}=\lim_{\xi_1\to0}\|\A(\cdot,\xi_1)\|_{L^2\left(\frac{1-\lambda}{1+\lambda}\xi_1,0\right)}=0.
\end{equation}
Now we can integrate (\ref{eq:someidentities1}) over the triangle with vortices $(\eta,\xi)$, $(0,\xi)$, and $(0,\eta+\xi)$ with $\eta=\frac{1-\lambda}{1+\lambda}\xi$, to get
\begin{equation}
\def\arraystretch{2.2}
\left.\begin{array}{rl}
     \int_{\lambda (1-t)}^{1-t}re(r,t)\,dr&=-\int_{\eta}^0r(e+m)(\eta',\xi)\,d\eta'-\int_{\eta+\xi}^{\xi}r(e-m)(0,\xi')\,d\xi'\\
       &\leq C\left(\|\A(\cdot,\xi)\|^2_{L^2(\eta,0)}+\Fl(\eta+\xi)\right)\\
       &=I+II.
      \end{array}\right.
\end{equation}
Since both I and II tend to zero as $t\to 1-$, we conclude the left hand side also tends to zero. This, combining with the fact that $F(\psi)\ge \psi^2$ on $K^{\lambda}_{{\rm ext}}(\tau)$, proves the case $\lambda=\lambda_0$ of our theorem.\\
Let us now assume the lemma holds for $\lambda=\lambda_0^j$ for some $j\ge 1$ and prove it holds also for $\lambda=\lambda_0^{j+1}$. By the relation $u=r^2\psi$ and the induction step, we have
\begin{equation}\label{eq:inductivestep}
\lim_{t\to1-}\int_{\lambda_0^j(1-t)\leq |x|\leq (1-t)}\left(\frac{u_t^2}{2}+\frac{|\nabla u|^2}{2}+\frac{u^2}{|x|^2}\right)(x,t)\,dx=0.
\end{equation}
Let us estimate $\psi=r^2u(r,t)$, for $(r,t)$ with $\lambda_0^{j+1}(1-t)\leq r\leq 1-t$. By (\ref{eq:smallnessnearlightcone1}) we see
\begin{equation}
|\psi|(r,t) \leq (1-t)^2|u(1-t,t)|+C\log{\left(\frac{1}{\lambda_0^j}\right)}\,\left(\int_{\lambda_0^j(1-t)\leq |x|\leq 1-t}|\nabla u|^2(x,t)+\frac{|u|^2}{|x|^2}(x,t)\,dx\right)^{\frac{1}{2}}
\end{equation}
for $\lambda_0^j(1-t)\leq r\leq 1-t$. Thus by (\ref{eq:tendtozero}) and (\ref{eq:inductivestep}), we get
\begin{equation}\label{eq:tiny}
\lim_{t\to 1-}\sup_{\lambda_0^j(1-t)\leq r\leq 1-t}|\psi(r,t)|=0.
\end{equation}
By a similar estimate to (\ref{eq:smallnessnearlightcone1}), for $(r,t)$ with $\lambda_0^{j+1}(1-t)\leq r\leq \lambda_0^j(1-t)$, we have
\begin{equation*}
\def\arraystretch{2.2}
\left.\begin{array}{rl}
|\psi(r,t)|=r^2|u|(r,t)&\leq |\psi|(\lambda_0^j(1-t),t)+\\
&\quad\quad\quad \,\,+\,C\log{\left(\frac{\lambda_0^j(1-t)}{r}\right)}\,\left(\int_{r\leq |x|\leq \lambda_0^j(1-t)}|\nabla u|^2+\frac{|u|^2(x,t)}{|x|^2}\,dx\right)^{\frac{1}{2}}\\
            &\leq  |\psi|(\lambda_0^j(1-t),t)+CM\log{\left(\frac{1}{\lambda_0}\right)}\\
         \end{array}\right.
\end{equation*}
Hence if we take $\tau$ sufficiently close to $1$, by the choice of $\lambda_0$, and (\ref{eq:tiny}) we obtain
\begin{equation}
|\psi(r,t)|\leq \frac{1}{10}\,\,\,\,\,\,{\rm on}\,\,\, \, K^{\lambda}_{{\rm ext}}(\tau).
\end{equation}
Thus we can repeat the arguments in the case of $\lambda=\lambda_0$ and conclude the proof. \\

Suppose $\OR{u}$ is a Type II blow up solution to equation (\ref{eq:mainwaveequation}), that blows up at $(x,t)=(0,1)$.  Then we know that there is essentially no energy in the self similar region $|x|\sim 1-t$ as $t$ approaches the blow up time $T_+=1$.  Based on this fact and virial identities, we can show the time average of kinetic energy of $u$ inside the lightcone $|x|\leq 1-t$ tends to zero. More precisely we have the following lemma.
\begin{lemma}\label{lm:timeaverages}
Let $\OR{u}$ be a radial finite energy solution to equation (\ref{eq:mainwaveequation}) on $R^6\times [0,1)$. Suppose $\OR{u}$ is Type II and blows up at spacetime point $(x,t)=(0,1)$. Then we have
\begin{eqnarray}
&&\lim_{t\to 1-}\frac{1}{1-t}\int_t^1\int_{|x|\leq 1-s}(\partial_tu)^2(x,s)\,dxds=0;\\
&&\lim_{t\to 1-}\frac{1}{1-t}\int_t^1\int_{|x|\leq 1-s}\left(|\nabla u|^2-|u|^3\right)(x,s)\,dxds=0.
\end{eqnarray}
\end{lemma}

\smallskip
\noindent
{\it Remark.} Set $a=u-v$ where $v$ is the regular part of $u$. Then from the above lemma, the fact that $v$ is regular near time $t=1$, and $u=v$ for $|x|\ge 1-t$ and $t\in (1-\delta,1)$, we see that
\begin{eqnarray}
&&\lim_{t\to 1-}\frac{1}{1-t}\int_t^1\int_{R^6}(\partial_ta)^2(x,s)\,dxds=0;\label{eq:vanishingtimeaver1}\\
&&\lim_{t\to 1-}\frac{1}{1-t}\int_t^1\int_{R^6}\left(|\nabla a|^2-|a|^3\right)(x,s)\,dxds=0.\label{eq:vanishingtimeaver2}
\end{eqnarray}

\smallskip
\noindent
{\it Proof.} Take a smooth radial function $\phi\in C_c^{\infty}(B_1(0))$ with $\phi\equiv 1$ on $B_{\frac{1}{2}}$. Then direct calculations show
\begin{eqnarray}
&&\frac{d}{dt}\int_{R^6}\partial_tu\,(x\cdot\nabla u+2u)\,\phi(\frac{x}{1-t})\,dx\nonumber\\
&&=-\int_{R^6}u_t^2\,\phi(\frac{x}{1-t})\,dx+O\left(\int_{\frac{1-t}{2}\leq |x|\leq 1-t}|\nabla_{x,t}u|^2+\frac{u^2}{|x|^2}\,dx\right);\label{eq:virial1}\\
{\rm and}&&\nonumber\\
&&\frac{d}{dt}\int_{R^6}\partial_tu\,(x\cdot\nabla u+3u)\,\phi(\frac{x}{1-t})\,dx\nonumber\\
&&=\int_{R^6}\left(-|\nabla u|^2+|u|^3\right)\,\phi(\frac{x}{1-t})\,dx+O\left(\int_{\frac{1-t}{2}\leq |x|\leq 1-t}|\nabla_{x,t}u|^2+\frac{u^2}{|x|^2}\,dx\right).\label{eq:virial2}
\end{eqnarray}
Using Theorem \ref{th:selfsimilar}, we see that
\begin{equation}
 \lim_{t\to 1-}\int_{\frac{1-t}{2}\leq |x|\leq 1-t}|\nabla_{x,t}u|^2+\frac{u^2}{|x|^2}\,dx=0.
\end{equation}
We claim that
\begin{equation}
\limsup_{s\to 1-}\frac{1}{1-s}\int_{|x|\leq 1-s}|\partial_tu|\left(|x||\nabla u|+|u|\right)(x,s)\,dx=0.
\end{equation}
Indeed, for any $\lambda\in (0,1)$, we have
\begin{eqnarray*}
&&\limsup_{s\to 1-}\frac{1}{1-s}\int_{|x|\leq \lambda(1-s)}|\partial_tu|\left(|x||\nabla u|+|u|\right)(x,s)\,dx\\
&&\leq C\limsup_{s\to 1-}\left(\lambda \int_{|x|\leq \lambda (1-s)}\frac{(\partial_tu)^2}{2}+\frac{|\nabla u|^2}{2}\,dx+\right.\\
&&\quad\quad\,\,\left.+\frac{1}{1-s}\left(\int_{|x|\leq\lambda (1-s)}(\partial_tu)^2\,dx\right)^{\frac{1}{2}}\,\left(\int_{|x|\leq\lambda (1-s)}u^2\,dx\right)^{\frac{1}{2}}\right)\\
&&\leq C\limsup_{s\to 1-}\left(\lambda \int_{|x|\leq \lambda (1-s)}\frac{(\partial_tu)^2}{2}+\frac{|\nabla u|^2}{2}\,dx\right.\\
&&\quad\quad\,\,\left.+\lambda \left(\int_{|x|\leq \lambda (1-s)}(\partial_tu)^2\,dx\right)^{\frac{1}{2}} \left(\int_{|x|\leq\lambda (1-s)}u^3\,dx\right)^{\frac{1}{3}}\right)\\
&&\leq C\lambda M^2\to 0, \quad {\rm as}\,\,\lambda\to 0+,
\end{eqnarray*}
and
\begin{eqnarray}
 &&\frac{1}{1-s}\int_{\lambda(1-s)<|x|\leq 1-s}|\partial_tu|\left(|x||\nabla u|+|u|\right)(x,s)\,dx\\
&&\leq C\int_{\lambda (1-s)\leq |x|\leq 1-s}\frac{(\partial_tu)^2}{2}+\frac{|\nabla u|^2}{2}\,dx+ C\int_{\lambda (1-s)\leq |x|\leq 1-s}|\partial_tu|\,\frac{|u|}{|x|}(x,s)\,dx\\
&&\leq C\int_{\lambda (1-s)\leq |x|\leq 1-s}\frac{(\partial_tu)^2}{2}+\frac{|\nabla u|^2}{2}+\frac{u^2}{|x|^2}\,dx\to 0, \quad {\rm as}\,\,s\to 1- 
\end{eqnarray}
by Theorem \ref{th:selfsimilar}. From the above inequalities, the claim follows straightforwardly.

Integrating the identities (\ref{eq:virial1}), (\ref{eq:virial2}) over the time interval $(t,s)$ and taking time average in the limit $s\to 1-$, we obtain
\begin{eqnarray*}
&&\frac{1}{1-t}\int_{R^6}\partial_tu\,(x\cdot\nabla u+2u)\phi(\frac{x}{1-t})\,dx=\frac{1}{1-t}\int_t^1\int_{R^6}u_t^2\,\phi(\frac{x}{1-s})\,dxds+o(1),\\
{\rm and}&&\\
&&\frac{1}{1-t}\int_{R^6}\partial_tu\,(x\cdot\nabla u+3u)\phi(\frac{x}{1-t})\,dx\\
&&\quad\quad\quad=\frac{1}{1-t}\int_t^1\int_{R^6}\left(-|\nabla u|^2+|u|^3\right)\phi(\frac{x}{1-s})\,dxds+o(1),
\end{eqnarray*}
as $t\to 1-$. Sending $t\to 1-$, we get 
\begin{eqnarray}
&&\lim_{t\to 1-}\frac{1}{1-t}\int_t^1\int_{R^6}u_t^2\,\phi(\frac{x}{1-s})\,dxds=0;\\
&&\lim_{t\to 1-}\frac{1}{1-t}\int_t^1\int_{R^6}\left(-|\nabla u|^2+|u|^3\right)\phi(\frac{x}{1-s})\,dxds=0.
\end{eqnarray}
This, combined with Theorem \ref{th:selfsimilar}, completes the proof. \\

To use (\ref{eq:vanishingtimeaver1}), (\ref{eq:vanishingtimeaver2}) we need the following real analysis lemma, which is a variant of Corollary 5.3 in \cite{DKMsmall}.

\begin{lemma}\label{lm:realanalysislemma111}
Let $a$ be as above. Then there exists a positive monotone function $\lambda(t)$ defined on $(0,1)$, with $\lim\limits_{t\to1-}\lambda(t)=\infty$, and a sequence of times $1>\tilde{t}_k>0$ approaching $1$ as $k$ goes to infinity, such that 
\begin{equation}\label{eq:importantbound}
\limsup_{k\to\infty}\sup_{\tau\in(0,1-\tilde{t}_k)}\frac{1}{\tau}\int_{\tilde{t}_k}^{\tk+\tau}\left(\lambda(s)\int_{R^6}(\partial_sa)^2(x,s)\,dx+\int_{R^6}|\nabla a|^2-|a|^3(x,s)\,dx\right)ds\leq 0.
\end{equation}
\end{lemma}

We postpone the proof of this lemma to Section \ref{sec:realanalysislemma} below. \\

Now we are ready to follow arguments in \cite{DKMsmall} to prove our main theorem in this section, which is the finite time blow up case of Theorem \ref{th:maintheoremsemi}.

\begin{theorem}\label{th:maintheorem}
Let $u\in C([0,T_+),\dot{H}^1(R^6))\cap L^2_tL^4_x(R^6\times I)$ for any $I\Subset [0,T_+)$ be the solution to equation (\ref{eq:mainwaveequation}), with radial initial data $(u_0,u_1)\in\HL(R^6)$. We assume without loss of generality $T_+=1$. Suppose $u$ is a Type II blow up solution with blow up point $(x,t)=(0,1)$. Then there exist a sequence of times $t_n\to 1-$, an integer $J_0\ge 1$, signs $\{l_j\}$ with $l_j\in\{\pm 1\}$, and sequences $\{\lambda_{jn}\}$ for $1\leq j\leq J_0$ with 
\begin{equation}
0<\lambda_{1n}\ll \lambda_{2n}\ll\cdots\ll\lambda_{J_0n}\ll 1-t_n, \,\,{\rm as}\,\,t_n\to 1-,
\end{equation}
such that
\begin{equation}
\OR{u}(t_n)=\OR{v}(1)+\sum_{j=1}^{J_0}l_j\left(\frac{1}{\lambda_{jn}^2}W(\frac{x}{\lambda_{jn}}),0\right)+o_{{\tiny \HL(R^6)}}(1),\,\,{\rm as}\,\,n\to\infty,
\end{equation}
where $\OR{v}$ is the regular part of $u$ and
\begin{equation}
W(x)=\frac{1}{(\frac{1}{24}+|x|^2)^2}
\end{equation}
 is the unique (up to scaling and sign) radial positive finite energy solution to
\begin{equation}
-\Delta W=|W|W.
\end{equation}
\end{theorem}

\smallskip
\noindent
{\it Proof.} We follow the arguments in \cite{DKMsmall} and \cite{Kenig4dWave}. With the new ingredient Lemma \ref{lm:realanalysislemma111}, we no longer need outer energy inequalities. Let $v$ be the regular part of $u$, and set $a=u-v$. By (\ref{eq:vanishingtimeaver1}), (\ref{eq:vanishingtimeaver2}) and Lemma \ref{lm:realanalysislemma111} there exists a sequence of times $t_n\to 1-$ and a positive function $\lambda(t)$ with $\lim_{t\to 1-}\lambda(t)=\infty$, such that
\begin{equation}\label{eq:importantbound1}
\limsup_{n\to\infty}\sup_{\tau\in(0,1-t_n)}\frac{1}{\tau}\int_{t_n}^{t_n+\tau}\left(\lambda(s)\int_{R^6}(\partial_sa)^2(x,s)\,dx+\int_{R^6}|\nabla a|^2-|a|^3(x,s)\,dx\right)ds\leq 0.
\end{equation}
Thanks to
\begin{equation}\label{eq:uniformbound1}
\limsup_{t\to 1-}\left(\|\nabla_{x,t}a\|_{L^2(R^6)}(t)+\|a\|_{L^3(R^6)}(t)\right)\leq CM<\infty,
\end{equation}
and the fact $\lim\limits_{t\to1-}\lambda(t)=\infty$, we get
\begin{equation}\label{eq:uniformbound2}
\lim_{n\to\infty}\sup_{\tau\in(0,1-t_n)}\frac{1}{\tau}\int_{t_n}^{t_n+\tau}\int_{R^6}(\partial_sa)^2(x,s)\,dxds=0.
\end{equation}
Letting $\tau\to0+$ in the above inequality, we also have
\begin{equation}\label{eq:limitzero}
\lim_{n\to\infty}\int_{R^6}(\partial_ta)^2(x,t_n)\,dx=0.
\end{equation}
Next we use the profile decompositions introduced in the context of wave equations by Bahouri and Gerard \cite{BaGe}, (see also \cite{Bulut} in the case of $R^{1+6}$). By passing to a subsequence of $\{t_n\}$ (which we still denote as $\{t_n\}$), we can assume $\OR{a}(t_n)$ has the following profile decomposition
\begin{equation}
\OR{a}(t_n)=\sum_{j=1}^J\frac{1}{\lambda_{jn}^2}U_L^j(\frac{x}{\lambda_{jn}},\frac{-t_{jn}}{\lambda_{jn}})+(w^J_{0,n},w^J_{1,n}),
\end{equation}
where $U_L^j(x,t)$ are radial finite energy solutions to the six dimensional linear wave equation. In addition, let $w^J_{n}(x,t)$ be the radial  solution to the linear wave equation with initial datat $(w^J_{0,n},w^J_{1,n})$, then $w^J_{n}$, $\lambda_{jn}$ and $t_{jn}$ satisfy the pseudo-orthogonality properties
\begin{eqnarray*}
&&t_{jn}\equiv 0\,\,{\rm for\,\,all\,\,}n,\,\,{\rm or}\,\,\lim_{n\to\infty}\frac{t_{jn}}{\lambda_{jn}}\in\{\pm\infty\};\\
&&\lim_{n\to\infty}\left(\frac{\lambda_{jn}}{\lambda_{j'n}}+\frac{\lambda_{j'n}}{\lambda_{jn}}+\frac{|t_{jn}-t_{j'n}|}{\lambda_{jn}}\right)=\infty, \,\,\,{\rm for\,\,all}\,\,j\neq j';\\
&&{\rm write}\,\,w^J_{n}(x,t)=\frac{1}{\lambda_{jn}^2}\widetilde{w}^j_{Jn}(\frac{x}{\lambda_{jn}},\,\frac{t-t_{jn}}{\lambda_{jn}}),\,\,{\rm then\,\,we\,\,have}\,\,\\
&&\widetilde{w}^j_{Jn}(\cdot,t)\rightharpoonup 0,\,\,{\rm as}\,\,n\to\infty, \,\,{\rm in}\,\,\HL(R^6)\,\,{\rm for}\,\,t\in R\,\,{\rm and\,\,}1\leq j\leq J.
\end{eqnarray*}
Moreover, the dispersive part $w^J_{n}$ satisfies the bound
\begin{equation}
\lim_{J\to\infty}\limsup_{n\to\infty}\|w^J_{n}\|_{L^{\infty}_tL^3_x(R^6\times[0,\infty))\cap L^2_tL^4_x(R^3\times [0,\infty))}=0.
\end{equation}
With the help of estimates (\ref{eq:uniformbound2}) and  (\ref{eq:limitzero}), we can use the same arguments as in (Section 5, \cite{DKMsmall} and the Erratum, \cite{DKMprofilesErratum}) to show that in the profile decomposition of $\OR{a}(t_n)$, all the nontrivial profiles are rescaled ground states, namely for some signs $l_j\in\{\pm 1\}$,
\begin{equation}\label{eq:firstdecom}
\OR{a}(t_n)=\sum_{j=1}^{J_0}l_j\left(\frac{1}{\lambda_{jk}^2}W(\frac{x}{\lambda_{jn}}),0\right)+(w^{J_0}_{0,n},\,w^{J_0}_{1,n}).
\end{equation}
$(w^{J_0}_{0,n},\,w^{J_0}_{1,n})$ is the dispersive part that satisfies 
\begin{equation}\label{eq:auxbound}
\lim_{n\to\infty}\left(\|w^{J_0}_{1,n}\|_{L^2(R^6)}+\|w^{J_0}_{0,n}\|_{L^3(R^6)}\right)=0.
\end{equation}
Note the bound (\ref{eq:importantbound1}) also implies
\begin{equation}
\limsup_{n\to\infty}\sup_{\tau\in(0,1-t_n)}\frac{1}{\tau}\int_{t_n}^{t_n+\tau}\int_{R^6}|\nabla a|^2-|a|^3(x,s)\,dxds\leq 0,
\end{equation}
and hence by taking $\tau\to0+$ 
\begin{equation}\label{eq:secondvirialvanishingblowupcase}
\limsup_{n\to\infty}\int_{R^6}|\nabla a|^2-|a|^3(x,t_n)\,dx\leq 0.
\end{equation}
Since $W$ is a steady state, by integration by parts we get
\begin{equation}
\int_{R^6}|\nabla W|^2-W^3\,dx=0.
\end{equation}
Thus by pseudo-orthogonality of the profiles in the profile decomposition and (\ref{eq:secondvirialvanishingblowupcase}), \footnote{Indeed, by the properties of the parameters, the profiles are concentrated in different {\it scales}, and have asymptotically vanishing interactions.} we see immediately that
\begin{equation}
\limsup_{n\to\infty}\int_{R^6}|\nabla w^{J_0}_{0,n}|^2-|w^{J_0}_{0,n}|^3(x)\,dx\leq 0.
\end{equation}
Hence, from the bound (\ref{eq:auxbound}) we obtain
\begin{equation}
\limsup_{n\to\infty}\int_{R^6}|\nabla w^{J_0}_{0,n}|^2\,dx\leq 0.
\end{equation}
Combining this with the bound (\ref{eq:auxbound}) we see
\begin{equation}
\lim_{n\to\infty}\|(w^{J_0}_{0,n},\,w^{J_0}_{1,n})\|_{\HL(R^6)}=0.
\end{equation}
Thus the decomposition (\ref{eq:firstdecom}) can be upgraded to
\begin{equation}
\OR{a}(t_n)=\sum_{j=1}^{J_0}l_j\left(\frac{1}{\lambda_{jn}^2}W(\frac{x}{\lambda_{jn}}),0\right)+o_{\HL(R^6)}(1),\,\,\,{\rm as}\,\,n\to\infty. 
\end{equation}
Note that $\OR{u}(t)=\OR{v}(t)+\OR{a}(t)$ with $v\in C([1-\delta,1],\dot{H}^1)$, hence
\begin{equation}
\OR{u}(t_n)=\OR{v}(1)+\sum_{j=1}^{J_0}l_j\left(\frac{1}{\lambda_{jn}^2}W(\frac{x}{\lambda_{jn}}),0\right)+o_{\HL(R^6)}(1),\,\,\,{\rm as}\,\,n\to\infty. 
\end{equation}
This completes the proof of our main theorem.

\end{subsection}

\begin{subsection}{Profiles for $6d$ wave, global existence case}
In this section we consider the case when the solution exists globally. More precisely, suppose 
\begin{equation*}
u\in C([0,\infty), \dot{H}^1(R^6))\cap L^2_tL^4_x([0,T)\times R^6)\,\, {\rm for\,\, any}\,\, T<\infty
\end{equation*}
is a radial solution to equation $(\ref{eq:mainwaveequation})$. Assume in addition that the solution is Type II, in the sense that
\begin{equation}
M:=\sup_{t\in [0,\infty)}\|\OR{u}(t)\|_{\HL(R^6)}<\infty.
\end{equation}
Our goal is to prove the global existence case of Theorem \ref{th:maintheoremsemi}.
\begin{theorem}\label{th:maintheorem2}
Suppose that $u\in C([0,\infty), \dot{H}^1(R^6))\cap L^2_tL^4_x([0,T)\times R^6)$ for any $T<\infty$, is a type II radial solution to equation $(\ref{eq:mainwaveequation})$. Then there exist a radial finite energy solution $u^L$ to the six dimensional linear wave equation, a sequence of times $t_n\to\infty$, an integer $J\ge 0$, sequences $\{\lambda_{jn},\,j\leq J\}$ with 
\begin{equation}
0<\lambda_{1n}\ll \lambda_{2n}\ll \dots\ll\lambda_{Jn}\ll t_n, \,\,\,{\rm as}\,\, n\to \infty,
\end{equation}
and signs $l_j\in\{\pm 1\}$, such that 
\begin{equation}
\OR{u}(t_n)=\sum_{j=1}^Jl_j\left(\frac{1}{\lambda_{jn}^2}W(\frac{x}{\lambda_{jn}}),0\right)+\OR{u}^L+o_{\HL(R^6)}(1),\,\,\,{\rm as}\,\,n\to\infty.
\end{equation}
\end{theorem}

The proof again follows closely that in \cite{Kenig4dWave}, with the new real analysis lemma as our main new ingredient, which allows us to rely on virial type identities only, thus overcoming the difficulty that there is no favorable outer energy inequality for solutions to the six dimensional linear wave equation. \\

We begin by recalling that the behavior of $u$ near the lightcone $||x|-|t||=O(1)$ is given by free radiation. More precisely we have
\begin{lemma}\label{lm:freeradiation}
Let $u$ be as above, then there exists a radial solution $u^L$ to linear wave equation with $\OR{u}^L(0)\in \HL(R^6)$, such that for any $A\in R$, 
\begin{equation}
\lim_{t\to\infty}\int_{|x|\ge t-A}\,\left|\nabla_{x,t}(u-u^L)\right|^2(x,t)\,dx=0.
\end{equation}
\end{lemma}
The proof of the lemma is exactly the same as the proof for Proposition 4.1 in \cite{Kenig4dWave}.\\

Next we show that there is essentially no energy in the self similar region. 
\begin{theorem}\label{th:selfsimilarenergyII}
Let $\lambda\in(0,1)$. Then
\begin{equation}
\lim_{A\to\infty}\,\limsup_{t\to\infty}\int_{\lambda t\leq |x|\leq t-A}\,\left(|\nabla_{t,x}u|^2+\frac{|u|^2}{|x|^2}\right)(x,t)\,dx=0.
\end{equation}
\end{theorem}

\bigskip
From this theorem, we can obtain the following simple consequence.
\begin{corollary}\label{cor:ext} 
Let $\lambda\in(0,1)$. Then
\begin{equation}
\lim_{t\to\infty}\|\nabla_{t,x}(u-u^L)\|_{L^2(|x|\ge\lambda t)}=0
\end{equation}
and
\begin{equation}
\lim_{t\to\infty}\|r^2u(r,t)\|_{L^{\infty}(r\ge \lambda t)}=0.
\end{equation}
\end{corollary}

\medskip
\noindent
{\it Proof of the Corollary.} From Theorem \ref{th:selfsimilarenergyII} and Lemma \ref{lm:energypartition} on the concentration property of solutions to the linear wave equation, we have 
\begin{eqnarray*}
&&\limsup_{t\to\infty}\|\nabla_{x,t}(u-u^L)\|_{L^2(|x|\in [\lambda t,t-A])}\\
&&\leq \limsup_{t\to\infty}\left(\|\nabla_{x,t}u\|_{L^2(|x|\in [\lambda t,t-A])}+\|\nabla_{x,t}u^L\|_{L^2(|x|\in [\lambda t,t-A])}\right)\to 0, \quad {\rm as}\,\,A\to\infty.
\end{eqnarray*}
By Lemma \ref{lm:freeradiation}, we also have for any $A\in R$
\begin{equation*}
\lim_{t\to\infty}\|\nabla_{t,x}(u-u^L)\|_{L^2(|x|\ge t-A)}=0.
\end{equation*}
Combining the two inequalities above, we get 
\begin{equation*}
\lim_{t\to\infty}\|\nabla_{t,x}(u-u^L)\|_{L^2(|x|\ge \lambda t)}=0.
\end{equation*}
Note that 
\begin{eqnarray}
\left|r^2(u-u^L)(r,t)\right|&\leq &\left|r^2\int_r^{\infty}\tau^{-\frac{5}{2}}\tau^\frac{5}{2}\partial_{\tau}(u-u^L)(\tau,t)\,d\tau\right|\nonumber\\
&\leq&C\left(\int_{|x|\ge r}|\nabla(u-u^L)|^2(x,t)\,dx\right)^{\frac{1}{2}}.\label{eq:pointwiseestimate}
\end{eqnarray}
Thus,
\begin{eqnarray*}
&&\limsup_{t\to\infty}\,\sup_{r\ge \lambda t}\,|r^2(u-u^L)(r,t)|\\
&&\leq C\limsup_{t\to\infty}\left(\int_{|x|\ge \lambda t}|\nabla(u-u^L)|^2\right)^{\frac{1}{2}}\\
&&=0.
\end{eqnarray*}
As $u^L$ is a linear solution, by dispersive estimates and approximation by smooth solutions with compact support, we see that
\begin{equation}\label{eq:dispersivesmall}
\lim_{t\to\infty}\|r^2u^L(\cdot,t)\|_{L^{\infty}(R^6)}=0.
\end{equation}
Thus
\begin{equation*}
\lim_{t\to\infty}\|r^2u(r,t)\|_{L^{\infty}(r\ge \lambda t)}=0.
\end{equation*}
This completes the proof of the Corollary.\\

As in the finite time blow up case, to prove Theorem \ref{th:selfsimilarenergyII}, we consider instead $\psi= r^2u(r,t)$, which satisfies
\begin{equation}\label{eq:equationforpsi}
\partial_{tt}\psi-\partial_{rr}\psi-\frac{\partial_r\psi}{r}+\frac{4\psi-|\psi|\psi}{r^2}=0,
\end{equation}
with
\begin{equation*}
\sup_{t\in[0,\infty)}\int_0^{\infty}\left(\frac{\psi_t^2}{2}+\frac{\psi_r^2}{2}+\frac{\psi^2}{r^2}\right)\,rdr<\infty
\end{equation*}
Theorem \ref{th:selfsimilarenergyII} is an immediate consequence of the following lemma.
\begin{lemma}\label{lm:noselfsimilarenergyII}
Let $\lambda\in(0,1)$, then
\begin{equation}
\lim_{A\to\infty}\,\limsup_{t\to\infty}\int_{\lambda t}^{t-A}\,\left(\frac{\psi_t^2}{2}+\frac{\psi_r^2}{2}+\frac{\psi^2}{r^2}\right)(r,t)\,rdr=0.
\end{equation}
\end{lemma}

\smallskip
\noindent
{\it Proof.} Let us recall the following notations:
\begin{eqnarray*}
&&f(\psi)=4\psi-|\psi|\psi;\\
&&F(\psi)=\int_0^{\psi}f=2\psi^2-\frac{|\psi|^3}{3};\\
&&e(r,t)=\frac{\psi_r^2}{2}+\frac{\psi_r^2}{2}+\frac{F(\psi)}{r^2};\\
&&m(r,t)=\psi_t\psi_r;\\
&&L(r,t)=-\frac{\psi_t^2}{2}+\frac{\psi_r^2}{2}+\frac{F(\psi)}{r^2}-\frac{2f(\psi)\psi_r}{r},
\end{eqnarray*}
and the following identities
\begin{eqnarray}
&&\partial_t(re)-\partial_r(rm)=0;\label{eq:id1}\\
&&\partial_t(rm)-\partial_r(re)=L.\label{eq:id2}
\end{eqnarray}
Note that by Lemma \ref{lm:freeradiation}
\begin{equation*}
\lim_{t\to\infty}\int_{|x|\ge t-A}\,|\nabla_{x,t}(u-u^L)|^2(x,t)\,dx=0, \quad{\rm for\,\,any}\,\,A\in R.
\end{equation*}
Hence, by an estimate similar to (\ref{eq:pointwiseestimate}), we see
\begin{equation}
\lim_{t\to\infty}\sup_{r\ge t-A}\,|r^2(u-u^L)(r,t)|=0,
\end{equation}
thus by (\ref{eq:dispersivesmall})
\begin{equation}
\lim_{t\to\infty}\sup_{r\ge t-A}|r^2u(r,t)|=0,\,\,\forall A\in R.
\end{equation}
Using
\begin{equation}
r^2u(r,t)=t^2u(t,t)-\int_r^t\partial_{\tau}\left(\tau^2u(\tau,t)\right)d\tau,
\end{equation}
by an application of H\"older inequality, 
\begin{equation}\label{eq:estimateonu}
r^2|u(r,t)|\leq t^2|u(t,t)|+C\log{\left(\frac{t}{r}\right)}\,\cdot\left(\int_{r\leq |x|\leq t}|\nabla u|^2+\frac{|u|^2}{|x|^2}(x,t)\,dx\right)^{\frac{1}{2}}.
\end{equation}
Hence, if we choose $\lambda_0=\lambda_0(M)<1$ sufficiently close to $1$, we have by inequality (\ref{eq:estimateonu})
\begin{equation}\label{eq:smallnearlightcone11}
\limsup_{t\to\infty}\sup_{\lambda t\leq r\leq t}\,r^2|u(r,t)|\leq CM\log{\frac{1}{\lambda_0}}<\frac{1}{10}.
\end{equation}
We prove our lemma holds for each $\lambda=\lambda_0^j$ by induction. We first consider the case $j=1$ and thus $\lambda=\lambda_0$.
Define
\begin{equation}
C_{\lambda}(A,T)=\{(r,t)|\,t\ge T,\,\lambda t\leq r\leq t-A\}.
\end{equation}
If we choose $T$ sufficiently large, then by (\ref{eq:smallnearlightcone11})
\begin{equation}
|\psi|=|r^2u(r,t)|\leq \frac{1}{10},\,\,{\rm on}\,\, C_{\lambda}(A,T).
\end{equation}
Hence $F(\psi)\ge \psi^2$, and consequently $re\pm rm\ge 0$ on $C_{\lambda}(A,T)$. Thus we can define 
\begin{eqnarray}
&&\A^2=re-rm,\\
&&\B^2=re+rm.
\end{eqnarray}
Note that here the definitions differ from the definitions (\ref{eq:defforAI}) and (\ref{eq:defforBI}) for Type II blow up case, to preserve the same geometric interpretations, as the lightcone is ``outgoing" in this case, instead of ``incoming" in the finite blow up case. By (\ref{eq:id1}) and (\ref{eq:id2}), we then have
\begin{eqnarray}
&&\partial_t\A^2+\partial_r\A^2=-L;\\
&&\partial_t\B^2-\partial_r\B^2=L.
\end{eqnarray}
As before we introduce (with corresponding modifications adapted to the new geometries)
\begin{equation}
\xi=t+r,\quad\, {\rm and}\,\,\,\eta=t-r.
\end{equation}
Then in the new coordinates $\xi,\,\eta$ we have
\begin{eqnarray}
&&\partial_{\xi}\A^2=-\frac{L}{2};\label{eq:identity1II}\\
&&\partial_{\eta}\B^2=\frac{L}{2}.\label{eq:identity2II}
\end{eqnarray}
We claim 
\begin{equation}\label{eq:controlofL}
L^2\leq C\frac{\A^2\B^2}{r^2}\quad {\rm on}\,\,C_{\lambda}(A,T).
\end{equation}
The proof of this claim follows exactly the same argument (which is just some algebra) as in the Type II blow up case, hence we omit it. Thus combining (\ref{eq:identity1II}), (\ref{eq:identity2II}) and (\ref{eq:controlofL}) we see
\begin{equation}\label{eq:keyinequalitiesII}
|\partial_{\xi}\A|\leq C\frac{\B}{r},\quad {\rm and}\,\,\,|\partial_{\eta}\B|\leq C\frac{\A}{r}, \quad {\rm on}\,\,C_{\lambda}(A,T).
\end{equation}
Since $\B^2$ and $\partial_{\eta}\B^2$ are both locally integrable in $t,\,r$ (and hence in $\xi,\,\eta$), by Lemma \ref{lm:restriction} we see $\B^2$ can be restricted to hyperplane of the form $\eta={\rm constant}$. Moreover, by approximation with smooth solutions, we have the following energy flux identity for $T<t_1<t_2$
\begin{eqnarray*}
&&\int_{t_1-A}^{\infty}\,\left(\frac{\psi_t^2}{2}+\frac{\psi_r^2}{2}+\frac{F(\psi)}{r^2}\right)(r,t_1)\,rdr\\
&&=\quad \int_{t_2-A}^{\infty}\,\left(\frac{\psi_t^2}{2}+\frac{\psi_r^2}{2}+\frac{F(\psi)}{r^2}\right)(r,t_2)\,rdr+\\
&&\quad\quad+\int_{t_1}^{t_2}\left(\frac{\psi_t^2}{2}+\frac{\psi_r^2}{2}+\frac{F(\psi)}{r^2}+\psi_t\psi_r\right)(l-A,l)\,(l-A)dl.
\end{eqnarray*}
Since on $C_{\lambda}(A,T)$, $F(\psi)\ge \psi^2$, we see
\begin{equation}
\int_{t_2-A}^{\infty}\,\left(\frac{\psi_t^2}{2}+\frac{\psi_r^2}{2}+\frac{F(\psi)}{r^2}\right)(r,t)\,rdr
\end{equation}
is decreasing in $t_2$ and tends to a finite limit as $t_2\to\infty$. We can therefore conclude
\begin{equation}\label{eq:estimateonthefluxII}
\lim_{t_1\to\infty}\int_{t_1}^{\infty}\left(\frac{\psi_t^2}{2}+\frac{\psi_r^2}{2}+\frac{F(\psi)}{r^2}+\psi_t\psi_r\right)(l-A,l)\,(l-A)dl=0.
\end{equation}
Fix $\left(\frac{1-\lambda}{1+\lambda}\xi_0,\xi_0\right)\in C_{\lambda}(A,T)$. Take $(\eta,\xi)\in C_{\lambda}(A,T)$ and $\xi_1>\xi>\xi_0$. By elementary geometry, $r<\xi<C(\lambda)r$, and $A\leq \eta\leq C(\lambda)r$. By the inequalities (\ref{eq:keyinequalitiesII}), we have
\begin{eqnarray*}
\A(\eta,\xi)&=&\A(\eta,\xi_1)-\int_{\xi}^{\xi_1}\partial_u\A(\eta,u)\,du\\
&\leq& \A(\eta,\xi_1)+C(\lambda)\int_{\xi}^{\xi_1}\frac{\B(\eta,u)}{u}\,du\\
&=&\A(\eta,\xi_1)+C(\lambda)\int_{\xi}^{\xi_1}\frac{\B(A,u)}{u}\,du+C(\lambda)\int_{\xi}^{\xi_1}\int_A^{\eta}\frac{\partial_v\B(v,u)}{u}\,dudv\\
&\leq& \A(\eta,\xi_1)+C(\lambda)\int_{\xi}^{\xi_1}\frac{\B(A,u)}{u}\,du+C(\lambda)\int_{\xi}^{\xi_1}\int_A^{\eta}\frac{\A(v,u)}{u^2}\,dudv.
\end{eqnarray*}
We follow similar arguments as in the Type II blow up case to estimate $\A$. The second term in the last line above can be estimated as
\begin{eqnarray*}
\int_{\xi}^{\xi_1}\frac{\B(A,u)}{u}\,du&\leq& \left(\int_{\xi}^{\xi_1}\B^2(A,u)\,du\right)^{\frac{1}{2}}\left(\int_{\xi}^{\xi_1}\frac{1}{u^2}\,du\right)^{\frac{1}{2}}\\
&\leq& C\frac{1}{|\xi|^{\frac{1}{2}}}\Fl(\xi,A)^{\frac{1}{2}},
\end{eqnarray*}
where we have set
\begin{equation}
\Fl(\xi,A):=\int_{\xi}^{\infty}\B^2(A,u)\,du.
\end{equation}
By (\ref{eq:estimateonthefluxII}), we know $\limsup_{\xi\to\infty}\Fl(\xi,A)=0$. From the geometric constraint $(\eta,\xi)\in C_{\lambda}(A,T)$, we see $A\leq \eta\leq \frac{1-\lambda}{1+\lambda}\xi$. Taking $L^2$ norm in $\eta$ and applying Minkowski inequality and H\"older inequality, we obtain
\begin{eqnarray*}
\|\A(\cdot,\xi)\|_{L^2(A,\frac{1-\lambda}{1+\lambda}\xi_0)}&\leq& \|\A(\cdot,\xi_1)\|_{L^2(A,\frac{1-\lambda}{1+\lambda}\xi_0)}+C(\lambda) F(\xi_0,A)^{\frac{1}{2}}+\\
&&\quad+C(\lambda)\int_{\xi}^{\xi_1}\left\|\int_A^{\eta}A(v,u)dv\right\|_{L^2\left(\eta\in(A,\frac{1-\lambda}{1+\lambda}\xi_0)\right)}\,\frac{du}{u^2}\\
&\leq& \|\A(\cdot,\xi_1)\|_{L^2(A,\frac{1-\lambda}{1+\lambda}\xi_0)}+C(\lambda) F(\xi_0,A)^{\frac{1}{2}}+\\
&&\quad+C(\lambda) \int_{\xi}^{\xi_1}\left\|\int_A^{\frac{1-\lambda}{1+\lambda}\xi}A(v,u)\chi_{(A,\eta)}(v)dv\right\|_{L^2\left(\eta\in(A,\frac{1-\lambda}{1+\lambda}\xi_0)\right)}\,\frac{du}{u^2}\\
&\leq& \|\A(\cdot,\xi_1)\|_{L^2(A,\frac{1-\lambda}{1+\lambda}\xi_0)}+C(\lambda) F(\xi_0,A)^{\frac{1}{2}}+\\
&&\quad+C(\lambda)\int_{\xi}^{\xi_1}\|\A(\cdot,u)\|_{L^2(A,\frac{1-\lambda}{1+\lambda}\xi_0)}\,\cdot\|\eta^{\frac{1}{2}}\|_{L^2\left(A,\frac{1-\lambda}{1+\lambda}\xi_0\right)}\,\frac{du}{u^2}\\
&\leq& \|\A(\cdot,\xi_1)\|_{L^2(A,\frac{1-\lambda}{1+\lambda}\xi_0)}+C(\lambda) F(\xi_0,A)^{\frac{1}{2}}+\\
&&\quad+C(\lambda)\int_{\xi}^{\xi_1}\|\A(\cdot,u)\|_{L^2\left(A,\frac{1-\lambda}{1+\lambda}\xi_0\right)}\frac{\xi_0}{u^2}\,du\\
\end{eqnarray*}
By Gronwall's inequality, we get 
\begin{equation*}
\|\A(\cdot,\xi_0)\|_{L^2(A,\frac{1-\lambda}{1+\lambda}\xi_0)}\leq C(\lambda)\left(\Fl(\xi_0,A)^{\frac{1}{2}}+\|\A(\cdot,\xi_1)\|_{L^2\left(A,\frac{1-\lambda}{1+\lambda}\xi_0\right)}\right).
\end{equation*}
Re-labelling $\xi_0$ as $\xi$, we have
\begin{equation}\label{eq:estimateonA}
\|\A(\cdot,\xi)\|_{L^2(A,\frac{1-\lambda}{1+\lambda}\xi)}\leq C(\lambda)\left(\Fl(\xi,A)^{\frac{1}{2}}+\|\A(\cdot,\xi_1)\|_{L^2\left(A,\frac{1-\lambda}{1+\lambda}\xi\right)}\right).
\end{equation}
We already know $\Fl(\xi,A)\to 0$ as $\xi\to\infty$. Let us estimate $\|\A(\cdot,\xi_1)\|_{L^2\left(A,\frac{1-\lambda}{1+\lambda}\xi\right)}$, in the limit $\xi_1\to\infty$. For any $\epsilon>0$, by Lemma \ref{lm:energypartition} and Lemma \ref{lm:freeradiation}, we can choose $A$ sufficiently large such that for any $L>A$, we have
\begin{equation}\label{eq:inequalityhahahaha}
\limsup_{t\to\infty}\int_{t-L\leq r\leq t-A}re\,dr\leq \epsilon.
\end{equation}
For fixed $\xi$, take the rectangle (in $\tilde{\xi}$ and $\tilde{\eta}$)
\begin{equation}
K(A,\xi_1)=\{(\tilde{\eta},\tilde{\xi})|\,A\leq\tilde{\eta}\leq \frac{1-\lambda}{1+\lambda}\xi,\,\xi_1-1\leq \tilde{\xi}\leq\xi_1+1\}.
\end{equation}
By (\ref{eq:inequalityhahahaha}), we can choose $A$ sufficiently large so that 
\begin{equation}
\limsup_{\xi_1\to\infty}\int_{K(A,\xi_1)}re\,drdt\leq \epsilon.
\end{equation}
Since $\A^2\leq Cre$, and $|\partial_{\xi}\A^2|\leq CL\leq Ce$, we have
\begin{equation}
\limsup_{\xi\to\infty}\int_{K(\xi_1)}\A^2\,d\xi' d\eta'\leq C\epsilon,\,\,\,{\rm and}\,\,\,\lim_{\xi_1\to\infty}\int_{K(\xi_1)}|\partial_{\xi}\A^2|\,d\xi' d\eta'=0.
\end{equation}
Hence we see from Lemma \ref{lm:restriction}
\begin{equation}
\limsup_{\xi_1\to\infty}\int_A^{\frac{1-\lambda}{1+\lambda}\xi}\,\A^2(\eta,\xi_1)\,d\eta\leq C\epsilon.
\end{equation}
Taking $\xi_1\to\infty$ in (\ref{eq:estimateonA}), we get
\begin{equation*}
\|\A(\cdot,\xi)\|_{L^2(A,\frac{1-\lambda}{1+\lambda}\xi)}\leq C(\lambda)\left(\Fl(A,\xi)^{\frac{1}{2}}+\epsilon\right),
\end{equation*}
and consequently
\begin{equation}\label{eq:vanishingoffluxtheotherside}
\limsup_{\xi\to\infty}\|\A(\cdot,\xi)\|_{L^2(A,\frac{1-\lambda}{1+\lambda}\xi)}\leq C(\lambda)\epsilon.
\end{equation}
Now integrate the identity $\partial_t(re)-\partial_r(rm)=0$ over the triangle with vortices $(\eta,\xi)=((1-\lambda)t,\,\,(1+\lambda)t)$, $(A,2t-A)$ and $(A,(1+\lambda)t)$. We obtain
\begin{eqnarray*}
&&\int_{\lambda t}^{t-A}\left(\frac{\psi_r^2}{2}+\frac{\psi_t^2}{2}+\frac{F(\psi)}{r^2}\right)(r,t)\,rdr\\
&&\quad=\int_{(1+\lambda)t}^{2t-A}\B^2(A,u)\,du+\int_A^{(1-\lambda)t}\A^2(\eta,(1+\lambda)t)\,d\eta\\
&&\quad \leq \int_{(1+\lambda)t}^{2t-A}\B^2(A,u)\,du+\|\A(\cdot,t+\lambda t)\|^2_{L^2(A,(1-\lambda)t)}\leq C(\lambda)\epsilon,
\end{eqnarray*}
for sufficiently large $t$, by (\ref{eq:vanishingoffluxtheotherside}). Since $\epsilon$ can be made arbitrarily small as long as we choose $A$ correspondingly large, the Lemma is proved in the case $\lambda=\lambda_0$. 

Now suppose the lemma is true for $\lambda=\lambda_0^j$, let us prove the lemma holds also for $\lambda=\lambda_0^{j+1}$. By the induction step and $\psi=r^2u$, we see
\begin{equation}
\lim_{A\to\infty}\,\limsup_{t\to\infty}\,\int_{\lambda_0^jt\leq |x|\leq t-A}|\nabla u|^2+\frac{|u|^2}{|x|^2}(x,t)\,dx=0.
\end{equation}
Setting $a=u-u^L$, then by Lemma \ref{lm:energypartition} and Lemma \ref{lm:freeradiation} we have
\begin{equation}
\limsup_{t\to\infty}\,\int_{|x|\ge\lambda_0^jt}|\nabla a|^2+\frac{|a|^2}{|x|^2}(x,t)\,dx=0.
\end{equation}
Then by estimates similar to estimate (\ref{eq:smallnessnearlightcone1}) we obtain
\begin{equation}
\sup_{r\in(\lambda_0^jt,t)}|r^2a(r,t)|\leq t^2|a(t,t)|+C\log{\left(\frac{1}{\lambda_0^j}\right)}\,\cdot\left(\int_{\lambda_0^jt\leq |x|\leq t}|\nabla a|^2+\frac{|a|^2}{|x|^2}(x,t)\,dx\right)^{\frac{1}{2}}\to 0,
\end{equation}
as $t\to\infty$. Hence for $\lambda t\leq r<\lambda_0^jt$ and $t$ large, we have
\begin{eqnarray*}
r^2|a(r,t)|&\leq& (\lambda_0^jt)^2|a(\lambda_0^jt,t)|+C\log{\left(\frac{1}{\lambda_0}\right)}\,\cdot\left(\int_{\lambda_0^{j+1}t\leq |x|\leq \lambda_0^jt}|\nabla a|^2+\frac{|a|^2}{|x|^2}(x,t)\,dx\right)^{\frac{1}{2}}\\
&\leq& (\lambda_0^jt)^2|a(\lambda_0^jt,t)|+CM\log{\frac{1}{\lambda_0}}.
\end{eqnarray*}
Hence 
\begin{equation}
\limsup_{t\to\infty}\sup_{\lambda_0^{j+1}t\leq r\leq \lambda_0^j t}r^2|a(r,t)|\leq CM\log{\frac{1}{\lambda_0}}<\frac{1}{10}.
\end{equation}
By the property of $u^L$, we then obtain
\begin{equation}
\limsup_{t\to\infty}\sup_{\lambda_0^{j+1}t\leq r\leq t}r^2|u(r,t)|\leq CM\log{\frac{1}{\lambda_0}}<\frac{1}{10}.
\end{equation}
Thus we can use exactly the same arguments as in the case $\lambda=\lambda_0$ to prove the Lemma for $\lambda=\lambda_0^{j+1}$, and finishes the induction.\\

Next we apply virial identities to $u$ to obtain some global estimates. We have the following result.
\begin{lemma}\label{lm:averagekinetic}
Let $u$ be as above and $u^L$ be the radiation term. Set $a=u-u^L$, then
\begin{eqnarray}
&&\lim_{t\to\infty}\frac{1}{t}\int_0^t\int_{R^6}(\partial_sa)^2(x,s)\,dxds=0;\label{eq:timeaverage1II}\\
&&\lim_{t\to\infty}\frac{1}{t}\int_0^t\int_{R^6}|\nabla a|^2-|a|^3(x,s)\,dxds=0.\label{eq:timeaverage2II}
\end{eqnarray}
\end{lemma}

\smallskip
\noindent
{\it Proof.} Take smooth cutoff function $\phi$ with supp\,$\phi\Subset B_{\frac{3}{4}}(0)$, $\phi|_{B_{\frac{1}{2}}}\equiv 1$. Note
\begin{equation}\label{eq:virial1II}
\frac{d}{dt}\int_{R^6}u_t\,\cdot (x\cdot\nabla u+2u)\phi(\frac{x}{t})\,dx=\int_{R^6}u_t^2(x,t)\phi(\frac{x}{t})\,dx+O\left(\int_{\frac{t}{2}\leq |x|\leq \frac{3}{4}t}|\nabla_{x,t} u|^2+\frac{|u|^2}{|x|^2}\,dx\right).
\end{equation}
By Theorem \ref{th:selfsimilarenergyII}, and similar arguments as in the Type II blow up case, we get 
\begin{equation}\label{eq:boundarytermstendtozero}
\lim_{t\to\infty}\frac{1}{t}\int_{R^6}u_t\,\cdot (x\cdot\nabla u+2u)\phi(\frac{x}{t})\,dx=0.
\end{equation}
By Theorem \ref{th:selfsimilarenergyII}, we have
\begin{equation*}
\int_{\frac{t}{2}\leq |x|\leq \frac{3}{4}t}|\nabla_{x,t} u|^2+\frac{|u|^2}{|x|^2}\,dx\to 0,\,\,{\rm as}\,\,t\to\infty,
\end{equation*}
thus
\begin{eqnarray*}
&&\limsup_{t\to \infty}\frac{1}{t}\int_0^t\int_{\frac{s}{2}\leq |x|\leq \frac{3}{4}s}|\nabla_{x,t} u|^2+\frac{|u|^2}{|x|^2}\,dxds\\
&&=\limsup_{t\to\infty}\frac{1}{t}\left(\int_0^{\epsilon t}+\int_{\epsilon t}^t\right)\int_{\frac{s}{2}\leq |x|\leq \frac{3}{4}s}|\nabla_{x,t} u|^2+\frac{|u|^2}{|x|^2}\,dx\,ds\\
&&\leq C\epsilon,\,\,{\rm for \,\,any\,\,}\epsilon>0.
\end{eqnarray*}
As a consequence we see that
\begin{equation}\label{eq:secondaveragetendstozero}
\lim_{t\to \infty}\frac{1}{t}\int_0^t\int_{\frac{s}{2}\leq |x|\leq \frac{3}{4}s}|\nabla_{x,t} u|^2+\frac{|u|^2}{|x|^2}\,dxds=0.
\end{equation}
Hence by integrating (\ref{eq:virial1II}) on $(0,t)$ and taking average, using (\ref{eq:boundarytermstendtozero}) and (\ref{eq:secondaveragetendstozero}), we obtain
\begin{eqnarray*}
\lim_{t\to\infty}\frac{1}{t}\int_0^t\int_{|x|\leq\frac{s}{2}}u_t^2(x,s)\,dxds=0.
\end{eqnarray*}
Similarly, by integrating 
\begin{eqnarray*}
&&\frac{d}{dt}\int_{R^6}u_t\,\cdot (x\cdot\nabla u+3u)\phi(\frac{x}{t})\,dx\\
&&=\int_{R^6}\left(-|\nabla u|^2+|u|^3\right)\phi(\frac{x}{t})\,dx+O\left(\int_{\frac{t}{2}\leq |x|\leq \frac{3}{4}t}|\nabla u|^2+\frac{|u|^2}{|x|^2}\,dx\right),
\end{eqnarray*}
we get
\begin{equation}
\lim_{t\to\infty}\frac{1}{t}\int_0^t\int_{|x|\leq\frac{s}{2}}|\nabla u|^2-|u|^3(x,s)\,dxds=0.
\end{equation}
Then by Lemma \ref{lm:energypartition} and Corollary \ref{cor:ext}, we conclude (\ref{eq:timeaverage1II}) and (\ref{eq:timeaverage2II}) hold. The lemma is proved.\\

The following real analysis lemma plays an important role in our analysis, which allows us to use the two bounds  (\ref{eq:timeaverage1II}), (\ref{eq:timeaverage2II}) simultaneously.
\begin{lemma}\label{lm:keylemmaII}
Let $u$ be a type II global solution to equation (\ref{eq:mainwaveequation}), and $u^L$ be the radiation term from Lemma \ref{lm:freeradiation}. Set $a=u-u^L$, then from Lemma \ref{lm:averagekinetic} we know that $a$ verifies
\begin{eqnarray}
&&\lim_{t\to\infty}\frac{1}{t}\int_0^t\int_{R^3}(\partial_ta)^2(x,s)\,dxds=0;\\
&&\lim_{t\to\infty}\frac{1}{t}\int_0^t\int_{R^3}|\nabla a|^2-|a|^3(x,s)\,dxds=0.
\end{eqnarray}
Then there exists an increasing function $\lambda(t)$ with $\lim\limits_{t\to\infty}\lambda(t)=\infty$, and a sequence of times $\tilde{t}_k\to\infty$ as $k\to\infty$, such that
\begin{equation}\label{eq:vanishingaverageII}
\limsup_{k\to\infty}\sup_{\tau\in (0,\frac{\tk}{4})}\frac{1}{\tau}\int_{\tk}^{\tk+\tau}\left(\lambda(t)\int_{R^6}(\partial_sa)^2(x,t)\,dx+\int_{R^6}|\nabla a|^2-|a|^3\,dx\right)ds\leq 0.
\end{equation}
\end{lemma}
We postpone the proof of this lemma to Section 4, which is similar in spirit to the proof of Corollary 5.3 of \cite{DKMsmall}, with some necessary modifications.\\

Now we are ready to give\\
\noindent
{\it Proof of main Theorem \ref{th:maintheorem2}}. Let $t_n\to\infty$ be the sequence give by Lemma \ref{lm:keylemmaII}, such that
\begin{equation}\label{eq:maintoolII}
\limsup_{n\to\infty}\sup_{\tau\in(0,\frac{t_n}{4})}\,\frac{1}{\tau}\int_{t_n}^{t_n+\tau}\left(\lambda(s)\int_{R^6}(\partial_sa)^2\,dx+\int_{R^6}|\nabla a|^2-|a|^3\,dx\right)ds\leq 0.
\end{equation}
Since $\lim\limits_{s\to\infty}\lambda(s)=\infty$, and $\sup\limits_{s\in(0,\infty)}\int_{R^6}|\nabla a|^2-|a|^3\,dx\leq C(M)$, we obtain
\begin{equation}\label{eq:anothertimeaverage}
\lim_{n\to\infty}\sup_{\tau\in(0,\frac{t_n}{4})}\frac{1}{\tau}\int_{t_n}^{t_n+\tau}\int_{R^6}(\partial_sa)^2\,dxds=0.
\end{equation}
Following the arguments \footnote{Strictly speaking, in \cite{DKMsmall}, it is assumed that (\ref{eq:anothertimeaverage}) holds for all $\tau>0$, not only for $\tau\in (0,\frac{t_n}{4})$. However a close inspection of the arguments in \cite{DKMsmall}  shows that (\ref{eq:anothertimeaverage}) already suffices. The key point is that all the profiles in the profile decomposition of $\OR{u}(t_n)$ have length scale $O(t_n)$, hence only time scales $O(t_n)$ are relevant.} in \cite{DKMsmall} (see also the Erratum in \cite{DKMprofilesErratum}), we conclude for a subsequence of $t_n$ (which we still denote as $t_n$ for ease of notation), that $\OR{u}(t_n)$ has the following profile decompositions
\begin{equation}\label{eq:decomp1II}
\OR{u}(t_n)=\sum_{j=1}^{J_0}l_j\left(\frac{1}{\lambda_{jn}^2}W(\frac{x}{\lambda_{jn}}),0\right)+\OR{u}^L(t_n)+(w^{J_0}_{0,n},\,w^{J_0}_{1,n}),
\end{equation} 
where $l_j\in\{\pm 1\}$, $J_0\ge 0$, $\OR{u}^L(t_n)$ is the linear solution from Lemma \ref{lm:freeradiation}, and
\begin{equation}
0\ll\lambda_{1n}\ll\lambda_{2n}\ll\cdots\ll\lambda_{J_0n}\ll t_n.
\end{equation}
Moreover, the dispersive term $(w^{J_0}_{0,n},\,w^{J_0}_{1,n})$ is uniformly bounded in $\HL$ and satisfies the following ``pseudo-orthogonality" condition:
\begin{equation}\label{eq:pseudoorthogonality}
(\lambda^2_{jn}w^{J_0}_{0,n}(\lambda_{jn}x),\lambda_{jn}^3w^{J_0}_{1,n}(\lambda_{jn}x))\rightharpoonup 0, \,\,\,{\rm as}\,\,n\to\infty,\,\,\,{\rm for\,\,each\,\,}j\leq J_0,
\end{equation}
and
\begin{equation}\label{eq:dispersivedecay}
\lim_{n\to\infty}\left(\|w^{J_0}_{0,n}\|_{L^3_x(R^6)}+\|w^{J_0}_{1,n}\|_{L^2(R^6)}\right)=0.
\end{equation}
Note (\ref{eq:maintoolII}) also implies
\begin{equation}
\limsup_{n\to\infty}\,\sup_{\tau\in (0,\frac{t_n}{4})}\,\frac{1}{\tau}\int_{t_n}^{t_n+\tau}\int_{R^6}|\nabla a|^2-|a|^3\,dxds\leq 0,
\end{equation}
and consequently by sending $\tau\to 0+$
\begin{equation}
\lim_{n\to\infty}\int_{R^6}|\nabla (u-u^L)|^2-|u-u^L|^3\,dx\leq 0.
\end{equation}
Hence, by the orthogonality of the profiles, and the fact 
\begin{equation}
\int_{R^6}|\nabla W|^2-|W|^3\,dx=0
\end{equation}
which follows easily from the fact $W$ is a steady state and an integration by parts argument,
we see
\begin{equation}
\limsup_{n\to\infty}\int_{R^6}|\nabla w^{J_0}_{0,n}|^2-|w^{J_0}_{0,n}|^3(x)\,dx\leq 0.
\end{equation}
By (\ref{eq:dispersivedecay}), we have
\begin{equation}
\lim_{n\to\infty}\int_{R^6}|w^{J_0}_{0,n}|^3(x)\,dx=0.
\end{equation}
Therefore
\begin{equation}
\lim_{n\to\infty}\int_{R^6}|\nabla w^{J_0}_{0,n}|^2(x)\,dx=0.
\end{equation}
Thus, combined with (\ref{eq:dispersivedecay}), the decompostion (\ref{eq:decomp1II}) can be upgraded to
\begin{equation}
\OR{u}(t_n)=\sum_{j=1}^Jl_j\left(\frac{1}{\lambda_{jn}^2}W(\frac{x}{\lambda_{jn}}),0\right)+\OR{u}^L(t_n)+o_{\HL(R^6)}(1),\,\,\,{\rm as}\,\,\,n\to\infty.
\end{equation} 
The theorem is proved.

\end{subsection}

\end{section}

\begin{section}{Asymptotics for two dimensional equivariant wave maps and radial Yang-Mills equation}
In this section we study the profiles for the two dimensional equivariant wave map and radial four dimensional Yang Mills  equation
\begin{equation}\label{eq:generalequation}
\partial_{tt}\psi-\partial_{rr}\psi-\frac{\partial_r\psi}{r}+\frac{f(\psi)}{r^2}=0,
\end{equation}
with $(\psi,\partial_t\psi)|_{t=0}=(\psi_0,\psi_1)$, where $f=gg'$.
We recall the following assumptions on $f(\rho)=g(\rho)g'(\rho)$:\\
\begin{eqnarray*}
(A1)&& G(x):=\int_0^x|g(y)|dy\to\infty, \,\,\,{\rm as}\,\,|x|\to\infty,\\
(A2)&& g\in C^3,\,\,\,{\rm and}\,\,\,\mathcal{V}:=\{l\in R:\,g(l)=0\} \,\,{\rm is\,\,nonempty\,\,and\,\,discrete}\\
&&{\rm\,\,with\,\,cardinality}\,|\mathcal{V}|\ge 2,\\
(A3)&&\forall l\in \mathcal{V}, \,\,g'(l)\in\{-2,-1,1,2\},\,\,{\rm and\,\,if\,\,}g'(l)=\pm1,\,\,{\rm then}\,\,g''(l)=0.
\end{eqnarray*}
Assumption (A1) ensures that finite energy solutions are bounded, (A2) implies that equation (\ref{eq:generalequation}) admits non-constant harmonic maps, (A3) is needed mainly for an easy local Cauchy theory and small data scattering results. By the work \cite{xiaoyi,Casey} we can relax assumption (A3) to
\begin{equation*}
(A3)'\quad\,\forall l\in \mathcal{V}, \,\,|g'(l)|\,\,{\rm is\,\,an\,\, integer}\ge 1,\,\,{\rm and\,\,if\,\,}g'(l)=\pm1,\,\,{\rm then}\,\,g''(l)=0.
\end{equation*}
Let us briefly review some results for equation (\ref{eq:generalequation}).
Define the energy and the Hilbert space $H\times L^2$ as follows: for $\OR{\phi}=(\phi_0,\phi_1)$ and $0\leq r_1<r_2\leq \infty$, 
\begin{eqnarray*}
&&\E(\OR{\phi},r_1,r_2):=\int_{r_1}^{r_2}\left(\phi_1^2+(\partial_r\phi_0)^2+\frac{g^2(\phi_0)}{r^2}\right)\,rdr,\\
&&\|\phi_0\|^2_{H([r_1,r_2])}:=\int_{r_1}^{r_2}\left((\partial_r\phi_0)^2+\frac{\phi_0^2}{r^2}\right)\,rdr,\\
&&\|\OR{\phi}\|_{H\times L^2([r_1,r_2])}=\int_{r_1}^{r_2}\left(\phi_1^2+(\partial_r\phi_0)^2+\frac{\phi_0^2}{r^2}\right)\,rdr.
\end{eqnarray*}
Also denote $\E(\OR{\phi})=\E(\OR{\phi},0,\infty)$ and $\mathcal{H}\times L^2=\{\OR{\phi}:\,\E(\OR{\phi})<\infty\}$.
The finiteness of the energy of $\OR{\psi}$ implies that $\psi_0$ is locally H\"older continuous, and moreover $\psi_0$ has limits both as $r\to\infty$ and $r\to 0$ (See  Lemma \ref{lm:boundonpsi} in the Appendix).

By \cite{Shatah2} we know that for each $(\psi_0,\psi_1)$ with $\E(\psi_0,\psi_1):=\int_0^{\infty}\left(\psi_1^2+(\partial_r\psi_0)^2+\frac{g^2(\psi_0)}{r^2}\right)\,rdr<\infty$, there exists a unique solution $\OR{\psi}\in C(I,\mathcal{H}\times L^2)$, defined on a maximal interval $I:=[0,T_+)$ which preserves the energy, and satisfies $\psi_0(0,t)\equiv \psi_0(0)$, $\psi(\infty,t)\equiv \psi_0(\infty)$ for each $t\in [0,T_+)$. Equation (\ref{eq:generalequation}) admits a non-constant harmonic map $Q$ of finite energy, \footnote{For Yang Mills equations, steady states are more often called {\it instantons}.}
\begin{equation*}
\partial_{rr}Q+\frac{\partial_rQ}{r}=\frac{f(Q)}{r^2}.
\end{equation*}
Thanks to \cite{Cote},  we know that $Q$ is monotone, and satisfies either $r\partial_rQ=g(Q)$ or $r\partial_rQ=-g(Q)$. Moreover $\{Q(0),\,Q(\infty)\}=\{l,\,m\}$ for some $l,\,m\in \mathcal{V}$, $l<m$, satisfying $\mathcal{V}\cap (l,m)=\emptyset$. $\E(Q,0)=2(G(m)-G(l))$. The linearized wave equation around any $l\in\mathcal{V}$ is:
\begin{equation}\label{eq:linearizedequationl}
\begin{array}{rc}
({\rm LW}_{\ell})&\partial_{tt}\phi-\partial_{rr}\phi-\frac{\partial_r\phi}{r}+\frac{g'(l)^2}{r^2}\phi=0.\\
\end{array}
\end{equation}

Direct calculation shows that if $\phi$ is a solution to $({\rm LW}_{\ell})$, then  $u_L=r^{-|g'(l)|}\phi$ verifies the radial $2|g'(l)|+2$ dimensional linear wave equation. Hence solutions to $({\rm LW}_{\ell})$ preserve the following energy
\begin{equation}
\|\OR{\phi}(t)\|^2_{\mathcal{H}_l\times L^2}:=\int_0^{\infty}\left(\phi_t^2+\phi_r^2+\frac{g'(l)^2\phi^2}{r^2}\right)\,rdr
\end{equation}

Suppose we have a finite energy solution $\psi$ to equation (\ref{eq:generalequation}) with (without loss of generality by rescaling) $T_+=1$. 
As in the semilinear wave equation case, we need the fact that there is asymptotically no energy in the self similar scale $r\sim 1-t$ as $t\to 1-$. 
\begin{lemma}\label{lm:noselfsimilarenergymap}
$\forall \lambda\in (0,1)$, we have
\begin{equation}\label{eq:noselfsimilarenergymap}
\lim_{t\to 1-}\int_{\lambda(1-t)}^{1-t}\left(\psi_t^2+\psi_r^2+\frac{g^2(\psi)}{r^2}\right)\,rdr=0.
\end{equation}
\end{lemma}
This important fact, which plays a crucial role in our argument, is proved in Lemma 2.2 of Shatah and Tahvildar-Zadeh \cite{Shatah1}.\\

Using virial identities, Lemma \ref{lm:noselfsimilarenergymap} implies
\begin{lemma}
Let $\psi$ be as above, then
\begin{equation}\label{eq:virial1general}
\lim_{t\to 1-}\frac{1}{1-t}\int_t^1\int_{0<r<1-s}(\partial_t\psi)^2\,rdrds=0,
\end{equation}
and
\begin{equation}\label{eq:virial2general}
\lim_{t\to1-}\frac{1}{1-t}\int_t^1\int_{0<r<1-s}\left(f'(\psi)\psi_r^2+\frac{f^2(\psi)}{r^2}\right)\,rdrds=0.
\end{equation}
\end{lemma}

\smallskip
\noindent
{\it Proof.} The first virial identity is contained in Corollary 2.3 of \cite{Shatah1}. Both virial identities follow from relatively straightforward calculations, let us only prove the second identity, which is new, and refer the reader to Corollary 2.3 in \cite{Shatah1} and Lemma 8 in \cite{KenigYang} for details for the first virial identity. Since $\psi$ is bounded, we have
\begin{equation*}
f^2(\psi)=g^2(\psi)\left(g'(\psi)\right)^2\leq C(\E(\psi_0,\psi_1),g)g^2(\psi).
\end{equation*}
Take $\eta\in C_c^{\infty}([0,1))$ with $\eta|_{[0,\frac{1}{2}]}\equiv 1$. Hence, 
\begin{eqnarray*}
(\ast):&=&\frac{d}{dt}\int_0^{\infty}\psi_tf(\psi)\,\eta\left(\frac{r}{1-t}\right)\,rdr\\
&&=\int_0^{\infty}\psi_{tt}f(\psi)\eta\left(\frac{r}{1-t}\right)\,rdr+\int_0^{\infty}\psi_t^2f'(\psi)\,\eta\left(\frac{r}{1-t}\right)\,rdr\\
&&\quad\quad+\int_0^{\infty}\psi_tf(\psi)\,\eta'\left(\frac{r}{1-t}\right)\,\frac{r}{(1-t)^2}\,rdr\\
&&=\int_0^{\infty}\psi_t^2f'(\psi)\,\eta\left(\frac{r}{1-t}\right)\,rdr+\int_0^{\infty}\left(\partial_{rr}\psi+\frac{\partial_r\psi}{r}-\frac{f(\psi)}{r^2}\right)f(\psi)\,\eta\left(\frac{r}{1-t}\right)\,rdr\\
&&\quad\quad+O\left(\int_{\frac{1-t}{2}}^{1-t}\left(\psi_t^2+\psi_r^2+\frac{g^2(\psi)}{r^2}\right)\,rdr\right)\\
&&=\int_0^{\infty}\psi_t^2f'(\psi)\,\eta\left(\frac{r}{1-t}\right)\,rdr-\int_0^{\infty}\left(f'(\psi)\psi_r^2+\frac{f^2(\psi)}{r^2}\right)\,\eta\left(\frac{r}{1-t}\right)\,rdr+\\
&&\quad\quad+O\left(\int_{\frac{1-t}{2}}^{1-t}\left(\psi_t^2+\psi_r^2+\frac{g^2(\psi)}{r^2}\right)\,rdr\right)\\
&&=\int_0^{\infty}\psi_t^2f'(\psi)\,\eta\left(\frac{r}{1-t}\right)\,rdr-\int_0^{\infty}\left(f'(\psi)\psi_r^2+\frac{f^2(\psi)}{r^2}\right)\,\eta\left(\frac{r}{1-t}\right)\,rdr+o(1),
\end{eqnarray*}
as $t\to 1-$. Note by H\"older inequality and Lemma \ref{lm:noselfsimilarenergymap}, as $t\to 1-$, 
\begin{equation}\label{eq:smallerthanlength}
\left|\int_0^{\infty}\psi_tf(\psi)\,\eta\left(\frac{r}{1-t}\right)\,rdr\right|=\left|\int_0^{\infty}r\psi_t\frac{f(\psi)}{r}\,\eta\left(\frac{r}{1-t}\right)\,rdr\right|=o(1-t).
\end{equation}
Integrate $(\ast)$ over time interval $(t,1)$ and take average, and then take limit $t\to 1-$. By (\ref{eq:smallerthanlength}), Lemma \ref{lm:noselfsimilarenergymap} and (\ref{eq:virial1general}), the lemma follows.\\

We note that for steady states, we have
\begin{lemma}\label{lm:harnomicvirial}
Let $Q$ be a finite energy harmonic map, then
\begin{equation}\label{eq:xxx}
\int_0^{\infty}\left(f'(Q)Q_r^2+\frac{f^2(Q)}{r^2}\right)\,rdr=0.
\end{equation}
\end{lemma}
We first remark that the above integral is absolutely convergent, since the boundedness of $Q$ implies $f(Q)$ is bounded and we have the inequality 
\begin{equation*}
f^2(Q)=g^2(Q)(g'(Q))^2\lesssim g^2(Q).
\end{equation*}
Multiplying equation (\ref{eq:generalequation}) with $f(Q)$ and integrating by parts, then (\ref{eq:xxx}) follows. Hence if $\lambda_n\ll 1-t_n$, then
\begin{equation}
\lim_{n\to\infty}\int_0^{1-t_n}\left(f'(Q(\frac{r}{\lambda_n}))\left(\frac{1}{\lambda_n}\partial_rQ(\frac{r}{\lambda_n})\right)^2+\frac{f^2(Q(\frac{r}{\lambda_n}))}{r^2}\right)\,rdr=0.
\end{equation}

\medskip
\noindent
Since $\psi$ is H\"older continuous in $(r,t)\in (0,\infty)\times (0,1]$, $\OR{\psi}(\cdot,1)$ is well defined with finite energy. Let $\phi(r,t)$ be the local in time solution to equation (\ref{eq:generalequation}) in $t\in[1-\delta,1]$ (for some small $\delta>0$), with initial data $\OR{\phi}(1)=\OR{\psi}(1)$. $\OR{\phi}$ is called the ``regular part" of $\OR{\psi}$, see \cite{cotesoliton} for full details. Set $\OR{a}=\OR{\psi}-\OR{\phi}+(\phi(0,1),0)$. The extra term $(\phi(0,1),0)$ is used to ensure that $a(0,t)\equiv \psi(0,t)$. We observe
\begin{lemma}\label{lm:fromphitoa}
\begin{equation*}
\lim_{t\to 1-}\frac{1}{1-t}\int_t^{1}\int_0^{\infty}(\partial_ta)^2\,rdrds=0,
\end{equation*}
and
\begin{equation*}
\lim_{t\to1-}\frac{1}{1-t}\int_t^1\int_0^{\infty}\left(f'(a)a_r^2+\frac{f^2(a)}{r^2}\right)\,rdrds=0.
\end{equation*}
\end{lemma}

\smallskip
\noindent
{\it Proof.} By finite speed of propagation, we get $\phi(r,t)=\psi(r,t)$, for $r\ge 1-t,\,t\in[1-\delta,1)$. Thus $a(r,t)=\phi(0,1)$ for $r\ge 1-t,\,t\in[1-\delta,1]$ with $g(\phi(0,1))=0$. Therefore we only need to show
\begin{eqnarray}
&&\lim_{t\to 1-}\frac{1}{1-t}\int_t^1\int_0^{1-s}(\partial_ta)^2\,rdrds=0,\label{eq:eqno1}\\
&&\lim_{t\to1-}\frac{1}{1-t}\int_t^1\int_0^{1-s}\left(f'(a)a_r^2+\frac{f^2(a)}{r^2}\right)\,rdrds=0.\label{eq:eqno2}
\end{eqnarray}
Since $\phi$ is a regular solution for $t\in[1-\delta,1]$, we have $\sup\limits_{r\leq 1-t}|\phi(r,t)-\phi(r,1)|\to 0$, as $t\uparrow 1$, and
\begin{equation}\label{eq:regularpartzero}
\lim_{t\to 1-} \int_0^{1-t}\left((\partial_t\phi)^2+(\partial_r\phi)^2+\frac{g^2(\phi)}{r^2}\right)(r,t)\,rdr=0.
\end{equation}
(\ref{eq:eqno1}) then follows immediately from (\ref{eq:virial1general}). We note, due to the non-vanishing of $g'(\phi(0,1))$, that
\begin{equation}
\int_0^{1-t}\frac{|\phi-\phi(0,t)|^2}{r^2}\,rdr\lesssim \int_0^{1-t}\frac{g^2(\phi)}{r^2}\,rdr\to 0, \,\,\,{\rm as}\,\,t\to 1-.
\end{equation}
To prove (\ref{eq:eqno2}), we also need the following elementary inequality
\begin{equation}
|f^2(x+y)-f^2(x)-2f'(x)f(x)y|\lesssim_{\|f\|_{C^2}, M} |y|^2, \,\,\,{\rm for}\,\,x, y\in [-M,M].
\end{equation}
Hence by H\"older inequality and (\ref{eq:regularpartzero})
\begin{eqnarray}
&&\left|\int_{r\leq 1-t}\frac{f^2(a)-f^2(\psi)}{r^2}\,rdr\right|\nonumber\\
&&\lesssim \int_{r\leq 1-t}\frac{|f(\psi)||\phi-\phi(0)|}{r^2}\,rdr+\int_{r\leq 1-t}\frac{|\phi-\phi(0)|^2}{r^2}\,rdr\nonumber\\
&&\lesssim \left(\int_{r\leq 1-t}\frac{f^2(\psi)}{r^2}\,rdr\right)^{\frac{1}{2}}\left( \int_{r\leq 1-t}\frac{|\phi-\phi(0)|^2}{r^2}\,rdr\right)^{\frac{1}{2}}+\int_{r\leq 1-t}\frac{|\phi-\phi(0)|^2}{r^2}\,rdr\nonumber\\
&&\to 0,\,\,\,{\rm as}\,\,t\to 0.\label{eq:hahainequality}
\end{eqnarray}
In the above inequalities, we have used the fact that $\phi(0,t)$ is constant for all $t$ in the time of existence. The above estimate, combined with the fact that $f'(\psi)\to f'(a)$ as $t\to 1-$ and (\ref{eq:virial2general}), implies (\ref{eq:eqno2}).\\

 Applying the real analysis lemma in Section 4, we obtain
\begin{lemma}\label{lm:thetrickgeneral}
There exists a sequence of times $t_n\uparrow 1$ such that 
\begin{eqnarray}
i)&&\sup_{\sigma\in(0,1-t_n)}\frac{1}{\sigma}\int_{t_n}^{t_n+\sigma}\int_0^{\infty}(\partial_ta)^2\,rdrds=o_n(1),\label{eq:key1}\\
ii)&&\limsup_{n\to\infty}\int_0^{\infty}\left(f'(a)a_r^2+\frac{f^2(a)}{r^2}\right)(r,t_n)\,rdr\leq 0.\label{eq:key2}
\end{eqnarray}
\end{lemma}


\medskip
\noindent
Let us also recall some results needed below in the global existence case ($T_+=\infty$).
\begin{lemma}\label{lm:noenergyglobal}
Suppose $\psi$ is a globally defined finite energy solution to equation (\ref{eq:generalequation}), and $\psi_0(\infty)=\ell\in\mathcal{V}$. Then there exists a radial finite energy solution $\phi^L$ to $({\rm LW}_{\ell})$, such that
\begin{equation}
\lim_{t\to\infty}\int_{t-A}^{\infty}\left(|\partial_r(\psi-\phi^L)|^2+|\partial_t(\psi-\phi^L)|^2+\frac{|\psi-\ell-\phi^L|^2}{r^2}\right)\,rdr=0,\,\,\,{\rm for\,\,any\,\,}A>0.
\end{equation}
Moreover, there is asymptotically no energy in the self similar region as $t\to\infty$, more precisely
for all $\lambda\in (0,1)$
\begin{equation}\label{eq:noenergyglobal}
\lim_{A\to\infty}\lim_{t\to\infty}\E(\OR{\psi}(t);\lambda t,t-A)=0.
\end{equation}
\end{lemma}
We refer to Proposition 5.1 in \cite{cotesoliton} and Proposition 2.1 in \cite{Cotemap2} for proofs.\\

Set $\OR{a}=\OR{\psi}-\OR{\phi}^L$, as an immediate consequence of Lemma \ref{lm:noenergyglobal} and Lemma \ref{lm:energypartition}, we obtain
\begin{lemma}\label{lm:psitoaglobal}
For any $\lambda\in(0,1)$,
\begin{equation}\label{eq:noenergya}
\lim_{t\to\infty}\int_{\lambda t}^{\infty}\left(|\partial_ra|^2+|\partial_ta|^2+\frac{|a-\ell|^2}{r^2}\right)\,rdr=0.
\end{equation}\label{lm:noenergya}
\end{lemma}

Lemma \ref{lm:noenergyglobal}, combined with virial identies, implies 
\begin{lemma}\label{lm:virialglobalgeneral}
Let $\psi$ be as above, then
\begin{eqnarray*}
i)&&\limsup_{T\to\infty}\frac{1}{T}\int_0^{T}\int_0^{\frac{T}{2}}|\partial_t\psi|^2\,rdrdt=0\\
ii)&&\limsup_{T\to\infty}\frac{1}{T}\int_0^{T}\int_0^{\frac{T}{2}}\left(f'(\psi)\psi_r^2+\frac{f^2(\psi)}{r^2}\right)\,rdrdt= 0.
\end{eqnarray*}
\end{lemma}

\smallskip
\noindent
{\it Proof.} The proof follows from similar calculations as in the finite time blow up case, this time based on (\ref{eq:noenergyglobal}), we omit the routine details.\\

%
Combining Lemma \ref{lm:psitoaglobal}, Lemma \ref{lm:virialglobalgeneral} and estimates similar to (\ref{eq:hahainequality}), we get
\begin{lemma}\label{lm:averagetozeroglobalmap}
\begin{eqnarray}
&&\limsup_{T\to\infty}\frac{1}{T}\int_0^T\int_0^{\infty}(\partial_ra)^2\,rdrdt=0,\\
&&\limsup_{T\to\infty}\frac{1}{T}\int_0^T\int_0^{\infty}\left(f'(a)a_r^2+\frac{f^2(a)}{r^2}\right)\,rdrdt=0.
\end{eqnarray}
\end{lemma}

Applying the real analysis lemma, we get
\begin{lemma}\label{lm:thetrickgeneralglobal}
There exists a sequence of times $t_n\uparrow \infty$, such that 
\begin{eqnarray*}
i)&&\sup_{\sigma\in(0,\frac{t_n}{4})}\frac{1}{\sigma}\int_{t_n}^{t_n+\sigma}\int_0^{\infty}(\partial_ta)^2\,rdrds=o_n(1),\\
ii)&&\limsup_{n\to\infty}\int_0^{\infty}\left(f'(a)a_r^2(r,t_n)+\frac{f^2(a)}{r^2}(r,t_n)\right)\,rdr\leq 0.
\end{eqnarray*}
\end{lemma}

In \cite{cotesoliton}, the following theorem was proved (using only the first virial type estimate from Lemma \ref{lm:thetrickgeneral} and Lemma \ref{lm:thetrickgeneralglobal} for the sequence $t_n\uparrow\infty$).
\begin{theorem}\label{th:cote}
Let $\OR{\psi}(t)$ be a finite energy solution to equation (\ref{eq:generalequation}). Then there exists a sequence of times $t_n\uparrow T_+$, an integer $J\ge 0$, $J$ sequences of scales $0<r_{Jn}\ll \cdots\ll r_{2n}\ll r_{1n}$ and $J$ harmonic maps $Q_1,\cdots,Q_J$ such that
\begin{equation*}
Q_J(0)=\psi(0),\,\,\,Q_{j+1}(\infty)=Q_j(0),\,\,\,{\rm for}\,\,j=1,\cdots,J-1,
\end{equation*}
and such that the following holds.\\
\quad\quad (1)\,\, (Global case) If $T_+=\infty$, denote $l=\psi(\infty)$. Then $Q_1(\infty)=l$, $r_{1n}\ll t_n$ and there exists a radial finite energy solution $\OR{\phi}^L$ to the linear wave equation $(LW_{\ell})$ such that 
\begin{equation}
\OR{\psi}(t_n)=\sum_{j=1}^J(Q_j(\cdot/r_{jn})-Q_j(\infty),0)+(l,0)+\OR{\phi}^L(t_n)+\OR{b}_n.
\end{equation}
(2)\,\,(Blow up case) If $T_+<\infty$, denote $l=\lim\limits_{t\to T_+}\psi(T_+-t,t)$ (it is well defined). Then $J\ge 1$, $\lambda_{1,n}\ll T_+-t_n$ and there exists a function $\OR{\phi}$ with $\E(\OR{\phi})<\infty$ such that $Q_1(\infty)=\phi(0)=l$ and 
\begin{equation}
\OR{\psi}(t_n)=\sum_{j=1}^J(Q_j(\cdot/r_{jn})-Q_j(\infty),0)+\OR{\phi}+\OR{b}_n.
\end{equation}
In both cases, $\OR{b}_n=(b_{0,n},b_{1,n})$ satisfies $\|b_n\|_{H\times L^2}=O(1)$ and vanishes in the following sense
\begin{equation}\label{eq:heheinequality}
\|\partial_rb_{0,n}\|_{L^2(r\leq r_{Jn},\,rdr)}+\|b_{0,n}\|_{L^{\infty}}+\|b_{1,n}\|_{L^2(R^+,\,rdr)}\to 0\,\,{\rm as}\,\,n\to\infty,
\end{equation}
and moreover, for any sequence $\lambda_n>0$ and $A>1$
\begin{equation}
\|b_{0,n}\|_{H(\frac{\lambda_n}{A}\leq r\leq A\lambda_n)}\to 0,\,\,{\rm as}\,\,n\to\infty.
\end{equation}
If we further assume\\

\medskip
\noindent
{\sl (A) for all $l\in \mathcal{V}$, $g'(l)\in \{-1,1\}$},\\

\medskip
\noindent
then 
\begin{equation}\label{eq:additionalconvergence}
\|(b_{0,n},b_{1,n})\|_{H\times L^2}\to 0,\,\,{\rm as}\,\,n\to\infty.
\end{equation}
\end{theorem}

\smallskip
\noindent
{\it Remark.} Let $\phi$ be the regular part of $\psi$ in the finite time blow up case, then $\OR{\phi}$ in the theorem is $\OR{\phi}(\cdot,1)$. We note that in the case of $k$ equivariant wave maps, $g(\rho)=k\sin{\rho}$, the assumption ({\sl A}) holds only if $k=1$. In the case of radial Yang Mills equation $g(\rho)=\frac{1}{2}(1-\rho^2)$, and ({\sl A}) fails. 
Hence in these cases, (\ref{eq:additionalconvergence}) is not proved in \cite{cotesoliton}. The main difficulty is the lack of a favorable outer energy inequality for the associated linear equations, which has played an important part in previous results of this type. \\

We shall prove (\ref{eq:additionalconvergence}) without assuming ({\sl A}), namely the following theorem.
\begin{theorem}\label{th:jk}
Theorem \ref{th:cote} holds without assuming {\sl (A)}.
\end{theorem}

\medskip
\noindent
{\it Proof of Theorem \ref{th:jk}.} We first consider the finite time blow up case. Assume without loss of generality $T_+=1$. Take the sequence $t_n\uparrow 1$ from Lemma \ref{lm:thetrickgeneral}. For this sequence of times, Theorem \ref{th:cote} holds, \footnote{In fact, the weaker Theorem \ref{th:cotesec2}, proved in the Appendix, is already sufficient for us.} thus there exist $J$ sequences of scales ($J\ge 1$)
\begin{equation*}
0<r_{Jn}\ll\cdots\ll r_{1,n}\ll 1-t_n,
\end{equation*}
$J$ finite energy non-constant steady states $Q_j$, and a function $\OR{\phi}$ with $\E(\OR{\phi})<\infty$, such that $Q_1(\infty)=\phi(0)=l$, $Q_J(0)=\psi_0(0)$, $Q_j(0)=Q_{j+1}(\infty)$ for $1\leq j\leq J-1$, and 
\begin{eqnarray*}
\OR{\psi}(t_n)&=&\sum_{j=1}^J(Q_j(\cdot/r_{jn})-Q_j(\infty),0)+\OR{\phi}+\OR{b}_n\\
&=&\OR{\phi}-(\phi(0),0)+(Q_1(\frac{r}{r_{1n}}),0)+\sum_{j=2}^J(Q_j(\frac{r}{r_{jn}})-Q_j(\infty),0)+(b_{0,n},b_{1,n}).
\end{eqnarray*}
where $\OR{b}_n=(b_{0,n},b_{1,n})$ vanishes in the following sense
\begin{equation}\label{eq:smallb}
\|b_{0,n}\|_{L^{\infty}}+\|b_{1,n}\|_{L^2(R^+,\,rdr)}\to 0\,\,{\rm as}\,\,n\to\infty,
\end{equation}
and moreover, for any sequence $\lambda_n>0$ and $A>1$
\begin{equation}\label{eq:smallb2}
\|b_{0,n}\|_{H(\frac{\lambda_n}{A}\leq r\leq A\lambda_n)}\to 0,\,\,{\rm as}\,\,n\to\infty.
\end{equation}

We note that $\OR{\phi}$ is simply $\OR{\phi}(\cdot,1)$, the regular part of $\psi$ at time $t=1$. $\phi$ satisfies
\begin{equation}
\int_0^{1-t_n}\left((\partial_r\phi)^2+(\partial_t\phi)^2+\frac{g^2(\phi)}{r^2}\right)\,rdr=o_n(1),\,\,{\rm as}\,\,n\to\infty.
\end{equation}
$Q_j$ satisfies
\begin{equation}
\int_0^{1-t_n}\left(\left(\frac{1}{r_{jn}}\partial_rQ_j(\frac{r}{r_{jn}})\right)^2+\frac{g^2(Q_j(\frac{r}{r_{jn}}))}{r^2}\right)\,rdr=o_n(1),\,\,{\rm as}\,\,n\to\infty.
\end{equation}
Our main additional tool is the following estimate, which is a consequence of $ii)$ in Lemma \ref{lm:thetrickgeneral} and the choice of the time sequence $t_n$,
\begin{equation}\label{eq:keyvirialgeneral}
\limsup_{n\to\infty}\int_0^{\infty}\left(f'(a)(\partial_ra)^2+\frac{f^2(a)}{r^2}\right)(r,t_n)\,rdr\leq 0,
\end{equation}
where 
\begin{equation}
\OR{a}(r,t_n)=\OR{\psi}(r,t_n)-\OR{\phi}+(\phi(0),0)=(Q_1(\frac{r}{r_{1n}}),0)+\sum_{j=2}^J(Q_j(\frac{r}{r_{jn}})-Q_j(\infty),0)+(b_{0,n},b_{1,n}).
\end{equation}

The key point in the proof of Theorem \ref{th:jk} is to study the expression (\ref{eq:keyvirialgeneral}) in various length scales. Let us note the following non-vanishing property of $f'$ which is important for us below:
\begin{equation}\label{eq:nonvanishingoff}
f'(l)=(g'(l))^2>0,\,\,{\rm for\,\,all\,\,}l\in\mathcal{V}.
\end{equation}
We make the following preliminary observation 
\begin{equation}
\forall L>1,\, 1\leq k\leq J,\,\, a(r_{kn}r,t_n)\to Q_k, \,\,\,{\rm uniformly\,\,in}\,\,r\in [\frac{1}{L},L].
\end{equation}
The above convergence follows immediately from
\begin{eqnarray*}
a(r,t_n)&=&Q_1(\frac{r}{r_{1n}})+\sum_{j=2}^J(Q_j(\frac{r}{r_{jn}})-Q_j(\infty))+b_{0,n}\\
&&=Q_k(\frac{r}{r_{kn}})+\sum_{j=k+1}^J(Q_j(\frac{r}{r_{jn}})-Q_j(\infty))+\sum_{j=1}^{k-1}(Q_j(\frac{r}{r_{jn}})-Q_j(0))+b_{0,n},
\end{eqnarray*}
and the fact that $\|b_{0,n}\|_{L^{\infty}}\to 0$ as $n\to\infty$.
This asymptotics, combined with (\ref{eq:smallb}) and (\ref{eq:smallb2}), implies that we can choose $\beta_{jn}\to\infty$, such that
\begin{eqnarray*}
&&\beta_{j+1,n}r_{j+1,n}\ll \frac{r_{jn}}{\beta_{jn}}, \,\,\,{\rm for\,\,}1\leq j\leq J-1;\\
&&\sup_{r\in[\frac{1}{\beta_{kn}},\,\beta_{kn}]}\left|a(r_{kn}r,t_n)- Q_k\right|\to 0,\,\,{\rm for\,\,}1\leq k\leq J;\\
&&\|\partial_rb_{0,n}\|_{L^2\left(\frac{r_{jn}}{\beta_{jn}}\leq r\leq \beta_{jn}r_{jn},\,rdr\right)}+\left\|\frac{b_{0,n}}{r^2}\right\|_{L^2\left(\frac{r_{jn}}{\beta_{jn}}\leq r\leq \beta_{jn}r_{jn},\,rdr\right)}\to 0, \,\,{\rm as}\,\,n\to\infty.
\end{eqnarray*}

Let us first study the expression $\int f'(a)(\partial_ra)^2\,rdr$ for $r$ in different scales.\\

{\it Case 1},  $r\ge \beta_{1n}r_{1n}$, then
\begin{equation}
a(r,t_n)=Q_1(\infty)+o_n(1).
\end{equation}
Hence 
\begin{eqnarray*}
&&\int_{\beta_{1n}r_{1n}}^{\infty}f'(a)(\partial_ra)^2(r,t_n)\,rdr\\
&&=\int_{\beta_{1n}r_{1n}}^{\infty}f'(Q_1(\infty))\left(\sum_{j=1}^J\frac{1}{r_{jn}}\partial_rQ_j(\frac{r}{r_{jn}})+\partial_rb_{0,n}\right)^2\,rdr+o_n(1)\\
&&=\int_{\beta_{1n}r_{1n}}^{\infty}f'(Q_1(\infty))(\partial_rb_{0,n})^2\,rdr+o_n(1).
\end{eqnarray*}

\bigskip
\noindent
{\it Case 2.}  $\beta_{kn}^{-1}r_{kn}\leq r\leq \beta_{kn}r_{kn}$. Then
\begin{equation*}
a(r,t_n)=Q_k(\frac{r}{r_{kn}})+o_n(1).
\end{equation*}
Consequently,
\begin{eqnarray*}
&&\int_{\beta_{kn}^{-1}r_{kn}}^{\beta_{kn}r_{kn}}f'(a)(\partial_ra)^2\,rdr\\
&&=\int_{\beta_{kn}^{-1}r_{kn}}^{\beta_{kn}r_{kn}}f'(Q_k(\frac{r}{r_{kn}}))\left(\sum_{j=1}^J\frac{1}{r_{jn}}\partial_rQ(\frac{r}{r_{jn}})+\partial_rb_{0,n}\right)^2\,rdr+o_n(1)\\
&&=\int_{\beta_{kn}^{-1}r_{kn}}^{\beta_{kn}r_{kn}}f'(Q_k(\frac{r}{r_{kn}}))\left(\frac{1}{r_{kn}}\partial_rQ_k(\frac{r}{r_{jn}})\right)^2\,rdr+o_n(1)\\
&&=\int_0^{\infty}f'(Q_k)(\partial_rQ_k)^2\,rdr+o_n(1).
\end{eqnarray*}

\bigskip
\noindent
{\it Case 3.} $\beta_{k+1,n}r_{k+1,n}\leq r\leq \frac{r_{kn}}{\beta_{kn}}$ for $1\leq k\leq J-1$. Then
\begin{equation*}
a(r,t_n)=Q_k(0)+o_n(1)=Q_{k+1}(\infty)+o_n(1).
\end{equation*}
Hence,
\begin{eqnarray*}
&&\int_{\beta_{k+1,n}r_{k+1,n}}^{\frac{r_{kn}}{\beta_{kn}}}f'(a)(\partial_ra)^2\,rdr\\
&&=\int_{\beta_{k+1,n}r_{k+1,n}}^{\frac{r_{kn}}{\beta_{kn}}}f'(Q_k(0))\left(\sum_{j=1}^J\frac{1}{r_{jn}}\partial_rQ(\frac{r}{r_{jn}})+\partial_rb_{0,n}\right)^2\,rdr+o_n(1)\\
&&=\int_{\beta_{k+1,n}r_{k+1,n}}^{\frac{r_{kn}}{\beta_{kn}}}f'(Q_k(0))(\partial_rb_{0,n})^2\,rdr+o_n(1).
\end{eqnarray*}

\bigskip
\noindent
{\it Case 4.} $r\leq \frac{r_{Jn}}{\beta_{Jn}}$. Then 
\begin{equation*}
a(r,t_n)=Q_J(0)+o_n(1)=\psi_0(0)+o_n(1).
\end{equation*}
Consequently,
\begin{eqnarray*}
&&\int_0^{\beta^{-1}_{Jn}r_{Jn}}f'(a)(\partial_ra)^2(r,t_n)\,rdr\\
&&=\int_0^{\beta^{-1}_{Jn}r_{Jn}}f'(Q_J(0))\left(\sum_{j=1}^J\frac{1}{r_{jn}}\partial_rQ(\frac{r}{r_{jn}})+\partial_rb_{0,n}\right)^2\,rdr\\
&&=\int_0^{\beta^{-1}_{Jn}r_{Jn}}f'(Q_J(0))(\partial_rb_{0,n})^2\,rdr+o_n(1).
\end{eqnarray*}

Thus in summary, we have 
\begin{eqnarray}
&&\int_0^{\infty}f'(a)(\partial_ra)^2\,rdr\nonumber\\
&&=\sum_{j=1}^J\int_0^{\infty}f'(Q_j)(\partial_rQ_j)^2\,rdr+\int_0^{\beta^{-1}_{Jn}r_{Jn}}f'(Q_J(0))(\partial_rb_{0,n})^2\,rdr+\nonumber\\
&&\quad+\sum_{j=1}^{J-1}\int_{\beta_{j+1,n}r_{j+1,n}}^{\frac{r_{jn}}{\beta_{jn}}}f'(Q_j(0))(\partial_rb_{0,n})^2\,rdr+\int_{\beta_{1n}r_{1n}}^{\infty}f'(Q_1(\infty))(\partial_rb_{0,n})^2\,rdr.\label{eq:firstterm}
\end{eqnarray}

Next we consider the term $\int \frac{f^2(a)}{r^2}\,rdr$ for $r$ in different scales, which is slightly trickier due to the nonlinear nature of $f$. For this, we need the following  lemma.
\begin{lemma}
For $x,\,y$ with $|x|+|y|\leq M$, $l\in \mathcal{V}$, we have
\begin{equation}\label{eq:estimateonf}
|f^2(l+x+y)-f^2(l+x)-f^2(l+y)|\leq C(\|f\|_{C^2([0,2M])},M)|x||y|.
\end{equation}
\end{lemma}
The proof follows from Calculus.\\

 We again divide into several cases.\\

{\it Case 1.}  $r\ge \beta_{1n}r_{1n}$, then
\begin{equation*}
a(r,t_n)=Q_1(\infty)+o_n(1).
\end{equation*}
We have
\begin{eqnarray*}
&&\int_{\beta_{1n}r_{1n}}^{\infty}\frac{f^2(a)}{r^2}\,rdr\\
&&=\int_{\beta_{1n}r_{1n}}^{\infty}\frac{f^2(Q_1(\infty)+\sum_{j=1}^J(Q_j(\frac{r}{r_{jn}})-Q_j(\infty))+b_{0,n})}{r^2}\,rdr\\
&&\,\\
&&\quad\quad{\rm by\,\,\,(\ref{eq:estimateonf}})\\
&&\,\\
&&=\int_{\beta_{1n}r_{1n}}^{\infty}\frac{f^2(Q_1(\infty)+b_{0n})}{r^2}\,rdr+\int_{\beta_{1n}r_{1n}}^{\infty}\frac{f^2(Q_1(\infty)+\sum_{j=1}^J(Q_j(\frac{r}{r_{jn}})-Q_j(\infty)))}{r^2}\,rdr+\\
&&\quad\quad+O\left(\int_{\beta_{1n}r_{1n}}^{\infty}\frac{|b_{0,n}||\sum_{j=1}^J(Q_j(\frac{r}{r_{jn}})-Q_j(\infty))|}{r^2}\,rdr\right).
\end{eqnarray*}
Note for $n$ large
\begin{eqnarray*}
&&\int_{\beta_{1n}r_{1n}}^{\infty}\frac{|\sum_{j=1}^J(Q_j(\frac{r}{r_{jn}})-Q_j(\infty))|^2}{r^2}\,rdr\\
&&\lesssim\int_{\beta_{1n}r_{1n}}^{\infty}\frac{\sum_{j=1}^J(Q_j(\frac{r}{r_{jn}})-Q_j(\infty))^2}{r^2}\,rdr\\
&&\lesssim\int_{\beta_{1n}r_{1n}}^{\infty}\frac{\sum_{j=1}^Jf^2(Q_j(\frac{r}{r_{jn}}))}{r^2}\,rdr=o_n(1),\,\,\,({\rm by}\,\,f^2\lesssim g^2).
\end{eqnarray*}
We obtain 
\begin{equation}
\int_{\beta_{1n}r_{1n}}^{\infty}\frac{f^2(a)}{r^2}\,rdr\gtrsim \int_{\beta_{1n}r_{1n}}^{\infty}\frac{f^2(Q_1(\infty)+b_{0,n})}{r^2}\,rdr+o_n(1)\gtrsim \int_{\beta_{1n}r_{1n}}^{\infty}\frac{b^2_{0,n}}{r^2}\,rdr+o_n(1).
\end{equation}
In the above inequality we have used the nonvanishing property (\ref{eq:nonvanishingoff}).

\bigskip
\noindent
{\it Case 2.} $\beta_{kn}^{-1}r_{kn}\leq r\leq \beta_{kn}r_{kn}$. Then
\begin{eqnarray*}
a(r,t_n)&=&Q_1(\frac{r}{r_{1n}})+\sum_{j=2}^J(Q_j(\frac{r}{r_{jn}})-Q_j(\infty))+b_{0,n}\\
&&=Q_k(\frac{r}{r_{kn}})+\sum_{j=k+1}^J(Q_j(\frac{r}{r_{jn}})-Q_j(\infty))+\sum_{j=1}^{k-1}(Q_j(\frac{r}{r_{jn}})-Q_j(0))+b_{0,n}.\\
\end{eqnarray*}
We note that for $j\ge k+1$
\begin{equation*}
\int_{\beta_{kn}^{-1}r_{kn}}^{\beta_{kn}r_{kn}}\frac{|Q_j(\frac{r}{r_{jn}})-Q_j(\infty)|^2}{r^2}\,rdr\lesssim \int_{\beta_{kn}^{-1}r_{kn}}^{\beta_{kn}r_{kn}}\frac{f^2(Q_j(\frac{r}{r_{jn}}))}{r^2}\,rdr=o_n(1),
\end{equation*}
and for $j\leq k-1$
\begin{equation*}
\int_{\beta_{kn}^{-1}r_{kn}}^{\beta_{kn}r_{kn}}\frac{|Q_j(\frac{r}{r_{jn}})-Q_j(0)|^2}{r^2}\,rdr\lesssim\int_{\beta_{kn}^{-1}r_{kn}}^{\beta_{kn}r_{kn}}\frac{f^2(Q_j(\frac{r}{r_{jn}}))}{r^2}\,rdr=o_n(1).
\end{equation*}
Let us recall that 
\begin{equation}\label{eq:tinyintegralb}
\int_{\beta_{kn}^{-1}r_{kn}}^{\beta_{kn}r_{kn}}\frac{b^2_{0,n}}{r^2}\,rdr=o_n(1).
\end{equation}

Define
\begin{equation}
\varphi_n=\sum_{j=k+1}^J(Q_j(\frac{r}{r_{jn}})-Q_j(\infty))+\sum_{j=1}^{k-1}(Q_j(\frac{r}{r_{jn}})-Q_j(0)).
\end{equation}
Then
\begin{eqnarray}
&&\int_{\beta_{kn}^{-1}r_{kn}}^{\beta_{kn}r_{kn}}\frac{\varphi_n^2}{r^2}\,rdr\nonumber\\
&&\lesssim \sum_{j=k+1}^J\int_{\beta_{kn}^{-1}r_{kn}}^{\beta_{kn}r_{kn}}\frac{(Q_j(\frac{r}{r_{jn}})-Q_j(\infty))^2}{r^2}\,rdr+\sum_{j=1}^{k-1}\int_{\beta_{kn}^{-1}r_{kn}}^{\beta_{kn}r_{kn}}\frac{(Q_j(\frac{r}{r_{jn}})-Q_j(0))^2}{r^2}\,rdr\nonumber\\
&&\lesssim \sum_{j=k+1}^J\int_{\beta_{kn}^{-1}r_{kn}}^{\beta_{kn}r_{kn}}\frac{f^2(Q_j(\frac{r}{r_{jn}}))}{r^2}\,rdr+\sum_{j=1}^{k-1}\int_{\beta_{kn}^{-1}r_{kn}}^{\beta_{kn}r_{kn}}\frac{f^2(Q_j(\frac{r}{r_{jn}}))}{r^2}\,rdr\label{eq:estimatevarphi}\\ 
&&=o_n(1).\nonumber
\end{eqnarray}
Therefore,
\begin{eqnarray*}
(\ast):&=&\int_{\beta_{kn}^{-1}r_{kn}}^{\beta_{kn}r_{kn}}\frac{f^2(a)}{r^2}\,rdr\\
&&=\int_{\beta_{kn}^{-1}r_{kn}}^{r_{kn}}\frac{f^2(Q_k(\frac{r}{r_{kn}})+\varphi_n+b_{0,n})}{r^2}\,rdr+\\
&&\quad+\int_{r_{kn}}^{\beta_{kn}r_{kn}}\frac{f^2(Q_k(\frac{r}{r_{kn}})+\varphi_n+b_{0,n})}{r^2}\,rdr\\
&&=\int_{\beta_{kn}^{-1}r_{kn}}^{r_{kn}}\frac{f^2(Q_k(0)+Q_k(\frac{r}{r_{kn}})-Q_k(0)+\varphi_n+b_{0,n})}{r^2}\,rdr+\\
&&\quad+\int_{r_{kn}}^{\beta_{kn}r_{kn}}\frac{f^2(Q_k(\infty)+Q_k(\frac{r}{r_{kn}})-Q_k(\infty)+\varphi_n+b_{0,n})}{r^2}\,rdr\\
\end{eqnarray*}
by (\ref{eq:estimateonf}) we get
\begin{eqnarray*}
(\ast)&=&\int_{\beta_{kn}^{-1}r_{kn}}^{r_{kn}}\frac{f^2(Q_k(0)+Q_k(\frac{r}{r_{kn}})-Q_k(0))}{r^2}\,rdr+\int_{\beta_{kn}^{-1}r_{kn}}^{r_{kn}}\frac{f^2(Q_k(0)+\varphi_n+b_{0,n})}{r^2}\,rdr\\
&&+\int_{r_{kn}}^{\beta_{kn}r_{kn}}\frac{f^2(Q_k(\infty)+Q_k(\frac{r}{r_{kn}})-Q_k(\infty))}{r^2}\,rdr\\
&&\quad+\int_{r_{kn}}^{\beta_{kn}r_{kn}}\frac{f^2(Q_k(\infty)+\varphi_n+b_{0,n})}{r^2}\,rdr\\
&&\quad+O\left(\int_{\beta_{kn}^{-1}r_{kn}}^{r_{kn}}\frac{|Q_k(\frac{r}{r_{kn}})-Q_k(0)||\varphi_n+b_{0,n}|}{r^2}\,rdr\right)\\
&&\quad\quad+ O\left(\int_{r_{kn}}^{\beta_{kn}r_{kn}}\frac{|Q_k(\frac{r}{r_{kn}})-Q_k(\infty)||\varphi_n+b_{0,n}|}{r^2}\,rdr\right)\\
\end{eqnarray*}
Hence by (\ref{eq:tinyintegralb}), (\ref{eq:estimatevarphi}), the preceding inequality and H\"older inequality, we see
\begin{eqnarray*}
(\ast)&=&\int_{\beta_{kn}^{-1}r_{kn}}^{\beta_{kn}r_{kn}}\frac{f^2(Q_k(\frac{r}{r_{kn}}))}{r^2}\,rdr+o_n(1)\\
&&=\int_0^{\infty}\frac{f^2(Q_k)}{r^2}\,rdr+o_n(1).
\end{eqnarray*}

\bigskip
\noindent
{\it Case 3.}  $\beta_{k+1,n}r_{k+1,n}\leq r\leq \frac{r_{kn}}{\beta_{kn}}$, $1\leq k\leq J-1$. We can re-write, using $Q_{j+1}(\infty)=Q_j(0)$, 
\begin{eqnarray*}
&&Q_1(\frac{r}{r_{1n}})+\sum_{j=2}^J(Q_j(\frac{r}{r_{jn}})-Q_j(\infty))+b_{0,n}\\
&&=\sum_{j=k+1}^J(Q_j(\frac{r}{r_{jn}})-Q_j(\infty))+Q_1(\frac{r}{r_{1n}})+\sum_{j=2}^kQ_j(\frac{r}{r_{n}})-\sum_{j=1}^{k-1}Q_j(0)+b_{0,n}\\
&&=\sum_{j=k+1}^J(Q_j(\frac{r}{r_{jn}})-Q_j(\infty))+Q_k(\frac{r}{r_{kn}})+\sum_{j=1}^{k-1}(Q_j(\frac{r}{r_{jn}})-Q_j(0))+b_{0,n}.
\end{eqnarray*}
We note for $j\ge k+1$,
\begin{eqnarray*}
&&\int_{\beta_{k+1,n}r_{k+1,n}}^{\frac{r_{kn}}{\beta_{kn}}}\frac{|Q_j(\frac{r}{r_{jn}})-Q_j(\infty)|^2}{r^2}\,rdr\\
&&\lesssim \int_{\beta_{k+1,n}r_{k+1,n}}^{\frac{r_{kn}}{\beta_{kn}}}\frac{f^2(Q_j(\frac{r}{\beta_{jn}}))}{r^2}\,rdr=o_n(1),
\end{eqnarray*}
and for $j\leq k-1$
\begin{eqnarray*}
&&\int_{\beta_{k+1,n}r_{k+1,n}}^{\frac{r_{kn}}{\beta_{kn}}}\frac{|Q_j(\frac{r}{r_{jn}})-Q_j(0)|^2}{r^2}\,rdr\\
&&\lesssim \int_{\beta_{k+1,n}r_{k+1,n}}^{\frac{r_{kn}}{\beta_{kn}}}\frac{f^2(Q_j(\frac{r}{\beta_{jn}}))}{r^2}\,rdr=o_n(1).
\end{eqnarray*}
Set
\begin{equation}
\varphi_n=\sum_{j=k+1}^J(Q_j(\frac{r}{r_{jn}})-Q_j(\infty))+\sum_{j=1}^{k-1}(Q_j(\frac{r}{r_{jn}})-Q_j(0)).
\end{equation}
Arguing exactly as in (\ref{eq:estimatevarphi}), we get
\begin{equation}
\int_{\beta_{k+1,n}r_{k+1,n}}^{\frac{r_{kn}}{\beta_{kn}}}\frac{\varphi_n^2+(Q_k(\frac{r}{r_{kn}})-Q_k(0))^2}{r^2}\,rdr=o_n(1).
\end{equation}
Therefore
\begin{eqnarray*}
&&\int_{\beta_{k+1,n}r_{k+1,n}}^{\frac{r_{kn}}{\beta_{kn}}}\frac{f^2(a)}{r^2}\,rdr\\
&&\gtrsim\int_{\beta_{k+1,n}r_{k+1,n}}^{\frac{r_{kn}}{\beta_{kn}}}\frac{(a(r,t_n)-Q_k(0))^2}{r^2}\,rdr\\
&&\gtrsim \int_{\beta_{k+1,n}r_{k+1,n}}^{\frac{r_{kn}}{\beta_{kn}}}\frac{|\varphi_n+(Q_k(\frac{r}{r_{kn}})-Q_k(0))+b_{0,n}|^2}{r^2}\,rdr\\
&&\gtrsim \int_{\beta_{k+1,n}r_{k+1,n}}^{\frac{r_{kn}}{\beta_{kn}}}\frac{b^2_{0,n}}{r^2}\,rdr+o_n(1).
\end{eqnarray*}

\bigskip
\noindent
{\it Case 4.}  $0<r\leq \frac{r_{Jn}}{\beta_{Jn}}$. Then
\begin{eqnarray*}
a(r,t_n)&=&Q_1(\frac{r}{r_{1n}})+\sum_{j=2}^J(Q_j(\frac{r}{r_{jn}})-Q_j(\infty))+b_{0,n}\\
&&=Q_J(\frac{r}{r_{Jn}})+\sum_{j=1}^{J-1}(Q_j(\frac{r}{r_{jn}})-Q_j(0))+b_{0,n}\\
&&=Q_J(0)+o_n(1).
\end{eqnarray*}
Hence we have for each $j\leq J$
\begin{equation*}
\int_0^{\frac{r_{Jn}}{\beta_{Jn}}}\frac{|Q_j(\frac{r}{r_{jn}})-Q_j(0)|^2}{r^2}\,rdr\lesssim \int_0^{\frac{r_{Jn}}{\beta_{Jn}}}\frac{f^2(Q_j(\frac{r}{r_{jn}}))}{r^2}\,rdr=o_n(1).
\end{equation*}
Consequently, 
\begin{eqnarray*}
&&\int_0^{\frac{r_{Jn}}{\beta_{Jn}}}\frac{f^2(a(r,t_n))}{r^2}\,rdr\\
&&\gtrsim\int_0^{\frac{r_{Jn}}{\beta_{Jn}}}\frac{(a(r,t_n)-Q_J(0))^2}{r^2}\,rdr\\
&&\gtrsim\int_0^{\frac{r_{Jn}}{\beta_{Jn}}}\frac{|b_{0,n}+\sum_{j=1}^{J}(Q_j(\frac{r}{r_{jn}})-Q_j(0))|^2}{r^2}\,rdr\\
&&\gtrsim\int_0^{\frac{r_{Jn}}{\beta_{Jn}}}\frac{|b_{0,n}|^2}{r^2}\,rdr+O\left(\int_0^{\frac{r_{Jn}}{\beta_{Jn}}}\frac{\sum_{j=1}^{J}(Q_j(\frac{r}{r_{jn}})-Q_j(0))^2}{r^2}\,rdr\right)\\
&&\gtrsim\int_0^{\frac{r_{Jn}}{\beta_{Jn}}}\frac{|b_{0,n}|^2}{r^2}\,rdr+o_n(1).
\end{eqnarray*}

Summarizing all cases, we get for some small $c>0$,
\begin{eqnarray}
&&\int_0^{\infty}\frac{f^2(a(r,t_n))}{r^2}\,rdr\nonumber\\
&&\ge c\int_{\beta_{1n}r_{1n}}^{\infty}\frac{b^2_{0,n}}{r^2}\,rdr+c\sum\limits_{k=1}^{J-1}\int_{\beta_{k+1,n}r_{k+1,n}}^{\frac{r_{kn}}{\beta_{kn}}}\frac{b^2_{0,n}}{r^2}\,rdr+c\int_0^{\frac{r_{Jn}}{\beta_{Jn}}}\frac{|b_{0,n}|^2}{r^2}\,rdr+\nonumber\\
&&\quad\quad\quad+ \sum\limits_{k=1}^J\int_0^{\infty}\frac{f^2(Q_k)}{r^2}\,rdr+o_n(1)\nonumber\\
&&\ge c\int_0^{\infty}\frac{b^2_{0,n}}{r^2}\,rdr+\sum_{k=1}^J\int_0^{\infty}\frac{f^2(Q_k)}{r^2}\,rdr+o_n(1).\label{eq:secondterm}
\end{eqnarray}

Hence, combining (\ref{eq:keyvirialgeneral}), (\ref{eq:firstterm}) and (\ref{eq:secondterm}) we conclude
\begin{eqnarray*}
0&\ge&\limsup_{n\to 0}\int_0^{\infty}\left(f'(a)(\partial_ra)^2+\frac{f^2(a)}{r^2}\right)(r,t_n)\,rdr\\
&\ge& \limsup_{n\to\infty}\left(\sum_{k=1}^{J}\int_0^{\infty}\left(f'(Q_k)(\partial_rQ_k)^2+\frac{f^2(Q_k)}{r^2}\right)\,rdr+c\int_0^{\infty}\left((\partial_rb_{0,n})^2+\frac{b_{0,n}^2}{r^2}\right)\,rdr\right)\\
&\gtrsim&\limsup_{n\to\infty}\int_0^{\infty}(\partial_rb_{0,n})^2+\frac{b_{0,n}^2}{r^2}\,rdr.
\end{eqnarray*}
Thus 
\begin{equation}
\lim_{n\to\infty}\|\partial_rb_{0,n}\|_{L^2((0,\infty),\,rdr)}+\left\|\frac{b_{0,n}}{r^2}\right\|_{L^2((0,\infty),\,rdr)}=0
\end{equation}
and the theorem in the finite time blow up case is proved.\\

In the global case set
\begin{equation}
\OR{a}=\OR{\psi}-\OR{\phi}^L.
\end{equation}
By Lemma \ref{lm:thetrickgeneralglobal} there exists a sequence of times $t_n\uparrow \infty$, such that 
\begin{eqnarray*}
i)&&\sup_{\sigma\in(0,\frac{t_n}{4})}\frac{1}{\sigma}\int_{t_n}^{t_n+\sigma}\int_0^{\infty}(\partial_ta)^2\,rdrds=o_n(1),\\
ii)&&\limsup_{n\to\infty}\int_0^{\infty}\left(f'(a)a_r^2(r,t_n)+\frac{f^2(a)}{r^2}(r,t_n)\right)\,rdr\leq 0.
\end{eqnarray*}
then we can use the decomposition from Theorem \ref{th:cote} and exactly the same arguments as in the finite blow up case to show $\|b_{0,n}\|_{H(0,\infty)}\to 0$, which finishes the proof.
\end{section}

\begin{section}{The real analysis lemma}\label{sec:realanalysislemma}
Let $u$ be a radial type II blow up solution to the six dimensional energy critical wave equation, that blows up at spacetime point $(x,t)=(0,1)$. Let $v$ be the regular part of the solution and set $a=u-v$. Then we know that $a$ has essentially no energy in the self similar region. Consequently by using virial identities, we can conclude the following two pieces of information:
\begin{eqnarray}
&&\lim_{t\to 1-}\frac{1}{1-t}\int_t^1\int_{R^6} (\partial_sa)^2(x,s)dxds=0;\\
&&\lim_{t\to 1-}\frac{1}{1-t}\int_t^1\int_{R^6} |\nabla a|^2-|a|^3dxds=0.
\end{eqnarray}
Moreover, we know
\begin{equation}\label{eq:uniformbound}
\sup_{t\in(0,1)}\left(\|\nabla_{x,t}a\|_{L^2(R^6)}(t)+\|a\|_{L^3(R^6)}(t)\right)\leq M<\infty,
\end{equation}
since $u$ is a type II blow up solution. 
Similarly, if $\psi$ is a finite energy equivariant wave map that blows up at $(r,t)=(0,1)$. Let $\phi$ be the regular part of $\psi$, and set $a=\psi-\phi$. Then we have
\begin{eqnarray}
&&\lim_{t\to 1-}\frac{1}{1-t}\int_t^1\int_{0}^{\infty} (\partial_sa)^2(r,s)\,rdrds=0;\\
&&\lim_{t\to 1-}\frac{1}{1-t}\int_t^1\int_0^{\infty} \left(f'(a)(\partial_ra)^2+\frac{f^2(a)}{r^2}\right)\,rdrds=0.
\end{eqnarray}
Set
\begin{eqnarray*}
&&g(t)=\int_{R^6} (\partial_ta)^2(x,t)\,dx,\\
&&h(t)=\int_{R^6} |\nabla a|^2-|a|^3\,dx,
\end{eqnarray*}
in the case of semilinear wave equation, and set 
\begin{eqnarray*}
&&g(t)=\int_{0}^{\infty} (\partial_ta)^2(r,t)\,rdr,\\
&&h(t)=\int_0^{\infty} \left(f'(a)(\partial_ra)^2+\frac{f^2(a)}{r^2}\right)\,rdr,
\end{eqnarray*}
in the case of equivariant wave maps.
Then in both cases, we have $g\ge 0$, $g,\,h$ bounded on $(0,1)$, moreover,
\begin{eqnarray}
&&\lim_{t\to 1-}\frac{1}{1-t}\int_t^1g(s)\,ds=0,\label{eq:eqforg}\\
&&\lim_{t\to 1-}\frac{1}{1-t}\int_t^1h(s)\,ds=0.\label{eq:eqforh}
\end{eqnarray}
Our goal is to prove the following statement.
\begin{lemma}\label{lm:realanalysislemma}
Suppose $g\ge 0$, $g,\,h$ are bounded on $(0,1)$, and satisfy (\ref{eq:eqforg}), (\ref{eq:eqforh}). Then there exists a positive monotone function $\lambda(t)$ defined on $(0,1)$, with $\lim\limits_{t\to1-}\lambda(t)=\infty$, and a sequence of times $\tilde{t}_k>0$ approaching $1$ as $k$ goes to infinity, such that 
\begin{equation}\label{eq:importantbound}
\limsup_{k\to\infty}\sup_{\tau\in(0,1-\tilde{t}_k)}\frac{1}{\tau}\int_{\tilde{t}_k}^{\tk+\tau}\left(\lambda(s)g(s)+h(s)\right)ds\leq 0.
\end{equation}
\end{lemma}

\smallskip
\noindent
{\it Remark.} This lemma implies that along the time sequence $\tk\uparrow 1$, 
\begin{eqnarray}
&&\lim_{k\to\infty}\sup_{\tau\in (0,1-\tk)}\frac{1}{\tau}\int_{\tilde{t}_k}^{\tk+\tau}g(s)\,ds=0,\\
&&\lim_{k\to\infty}\sup_{\tau\in (0,1-\tk)}\frac{1}{\tau}\int_{\tilde{t}_k}^{\tk+\tau}h(s)\,ds\leq 0.
\end{eqnarray}

To prove Lemma \ref{lm:realanalysislemma}, we first establish the following preliminary step.
\begin{lemma}\label{lm:boostestimate}
Let $g$ be as above. Then there exists a positive monotone function $\lambda(t)$ defined on $(0,1)$, with $\lim\limits_{t\to1-}\lambda(t)=\infty$ and an increasing sequence $t_k\in (0,1)$ approaching $1$, such that
\begin{equation}
\lim_{k\to\infty}\frac{1}{1-t_k}\int_{t_k}^1\lambda(s)g(s)\,ds=0.
\end{equation}

\end{lemma}

\smallskip
\noindent
{\it Proof.} Since $\lim\limits_{t\to 1-}\frac{1}{1-t}\int_t^1g(s)\,ds=0$, we can choose an increasing sequence $t_k$ with the following properties:
\begin{eqnarray}
&&1-t_{k+1}\leq \frac{1-t_k}{4}, \,\,{\rm for\,\,all\,\,}k\ge 1;\label{eq:no1}\\
&&\frac{1}{1-t_{k+1}}\int_{t_{k+1}}^1g(s)\,ds\leq \frac{1}{4(1-t_k)}\int_{t_k}^1g(s)\,ds.
\end{eqnarray}
Define $\lambda(t)$ as
\begin{equation}
\lambda(t)=\sum_{k=1}^{\infty}2^k\chi_{[t_k,t_{k+1})}(t), \,\,{\rm for\,\,}t\in(0,1).
\end{equation}
Since $t_k$ is increasing, there is only one term in the sum for each $t\in(0,1)$. It's easy to check that $\lambda(t)$ is positive and monotone with $\lim\limits_{t\to1-}\lambda(t)=\infty$.
\end{section}
Let us estimate $\frac{1}{1-t_j}\int_{t_j}^1\lambda(s)g(s)\,ds$. We have
\begin{eqnarray*}
&&\frac{1}{1-t_j}\int_{t_j}^1\lambda(s)g(s)\,ds\\
&&\quad= \sum_{k\ge j}\frac{1}{1-t_j}\int_{t_k}^{t_{k+1}}\lambda(s)g(s)\,ds\\
&&\quad\,\,\leq \sum_{k\ge j}2^k\frac{1}{1-t_j}\int_{t_k}^{t_{k+1}}g(s)\,ds\\
&&\quad\,\,\,\leq \sum_{k\ge j}2^k\frac{1-t_k}{1-t_j}\frac{1}{1-t_k}\int_{t_k}^1g(s)\,ds\\
&&\quad\,\,\,\,\leq \sum_{k\ge j}2^k\,4^{j-k}\,4^{1-k}\frac{1}{1-t_1}\int_{t_1}^1g(s)\,ds\\
&&\quad\quad\leq C2^{-j}\frac{1}{1-t_1}\int_{t_1}^1g(s)\,ds\to0+,\,\,{\rm as\,\,}j\to\infty.
\end{eqnarray*}
Thus the lemma is proved.\\

\medskip
\noindent
Now we can turn to \\
{\it Proof of Lemma \ref{lm:realanalysislemma}.} Take $\lambda(t)$ and $t_k$ from Lemma \ref{lm:boostestimate}. Let us denote
\begin{equation}
\tg(t)=\lambda(t)g(t).
\end{equation}
We know $\tg(t)\ge 0$ and $\lim\limits_{k\to\infty}\frac{1}{1-t_k}\int_{t_k}^1\tg(s)ds=0$. Also $|h(t)|\leq C(M)$, $\lim\limits_{t\to1-}\frac{1}{1-t}\int_t^1h(s)ds=0$. Let us consider the sequence
\begin{equation}
\frac{1}{1-t_k}\int_{t_k}^1(\tg(s)+h(s))ds, \,\,k\ge 1.
\end{equation}
We distinguish two cases.\\
{\it Case 1.} There exists a subsequence of $t_k$ (which we still denote as $t_k$ for ease of notation) such that $\frac{1}{1-t_k}\int_{t_k}^1(\tg(s)+h(s))ds>0$ for each $k$. In this case we will follow the arguments in \cite{DKMsmall} to find $\tk$. Since $\frac{1}{s-t_k}\int_{t_k}^s(\tg(\tau)+h(\tau))d\tau$ is continuous in $s\in[\frac{t_k+1}{2},1]$, there exists $\tk\in[\frac{t_k+1}{2},1]$ such that
\begin{equation}
\frac{1}{\tk-t_k}\int_{t_k}^{\tk}(\tg(\tau)+h(\tau))d\tau=\max_{s\in[\frac{t_k+1}{2},1]}\frac{1}{s-t_k}\int_{t_k}^s(\tg(\tau)+h(\tau))d\tau.
\end{equation}
Take $j>k$ sufficiently large such that $t_j>\frac{t_k+1}{2}$ and
\begin{equation}
\frac{1}{1-t_j}\int_{t_j}^1(\tg(\tau)+h(\tau))ds<\frac{1}{2(1-t_k)}\int_{t_k}^1(\tg(\tau)+h(\tau))d\tau,
\end{equation}
then it follows immediately that
\begin{equation}
\frac{1}{t_j-t_k}\int_{t_k}^{t_j}(\tg(\tau)+h(\tau))d\tau>\frac{1}{1-t_k}\int_{t_k}^1(\tg(\tau)+h(\tau))d\tau.
\end{equation}
Thus $\tk<1$. Since $\max\limits_{s\in[\frac{t_k+1}{2},1]}\frac{1}{s-t_k}\int_{t_k}^s(\tg(\tau)+h(\tau))d\tau$ is achieved at $s=\tk$, a moment's reflection shows
\begin{equation}\label{eq:controlaverage}
\sup_{\sigma\in(0,1-\tilde{t}_k)}\frac{1}{\sigma}\int_{\tilde{t}_k}^{\tk+\sigma}(\tg(\tau)+h(\tau))d\tau\leq \frac{1}{\tk-t_k}\int_{t_k}^{\tk}(\tg(\tau)+h(\tau))d\tau.
\end{equation}
Let us now estimate the right hand side from above. Since $\tk\ge\frac{t_k+1}{2}$ and $\tg(t)\ge 0$, we have
\begin{equation}\label{eq:controlofaverageforg}
\frac{1}{\tk-t_k}\int_{t_k}^{\tk}\tg(\tau)d\tau\leq \frac{2}{1-t_k}\int_{t_k}^1\tg(\tau)d\tau\,\to \,0,\,\,{\rm as}\,\,k\to\infty.
\end{equation}
We claim 
\begin{equation}\label{eq:controlaveragesforh}
\left|\frac{1}{\tk-t_k}\int_{t_k}^{\tk}h(s)ds\right|\,\leq \,3\left|\frac{1}{1-t_k}\int_{t_k}^1h(s)ds\right|+3\left|\frac{1}{1-\tk}\int_{\tk}^1h(s)ds\right|.
\end{equation}
Otherwise, recalling $\tk\ge \frac{t_k+1}{2}$, we would have
\begin{equation}
\left|\int_{\tk}^1h(s)ds\right|< \frac{1-\tk}{3(\tk-t_k)}\left|\int_{t_k}^{\tk}h(s)ds\right|\leq \frac{1}{3}\left|\int_{t_k}^{\tk}h(s)ds\right|.
\end{equation}
Then
\begin{equation}
\left|\frac{1}{1-t_k}\int_{t_k}^1h(s)ds\right|>\left| \frac{2}{3(1-t_k)}\int_{t_k}^{\tk}h(s)ds\right|\ge\frac{1}{3}\left|\frac{1}{\tk-t_k}\int_{t_k}^{\tk}h(s)ds\right|,
\end{equation}
a contradiction. Thus the claim follows. Since $\lim\limits_{t\to 1-}\frac{1}{1-t}\int_t^1h(s)ds=0$, combining bounds (\ref{eq:controlaverage},\ref{eq:controlofaverageforg},\ref{eq:controlaveragesforh}), we obtain
\begin{equation}
\limsup_{k\to\infty}\sup_{\sigma\in(0,1-\tilde{t}_k)}\frac{1}{\sigma}\int_{\tilde{t}_k}^{\tk+\sigma}(\tg(\tau)+h(\tau))d\tau\leq 0.
\end{equation}
{\it Case 2.} We show we can always reduce to Case 1. It is clear that we can choose a subsequence of $\{t_k\}$ (which we still denote as $\{t_k\}$ for ease of notation), such that
\begin{eqnarray}
&&1-t_{k+1}\leq \frac{1-t_k}{4},\\
&&\epsilon_k:=\left|\frac{1}{1-t_k}\int_{t_k}^1(\tg(s)+h(s))ds\right|\leq \frac{1}{4^k}.
\end{eqnarray}
Now set 
\begin{equation}
f(t)=\sum_{k=1}^{\infty}(4\,\epsilon_k+2^{-k})\,\chi_{[t_k,t_{k+1})}(t),\,\,{\rm for}\,\,\,t\in(0,1).
\end{equation}
It's straightforward to check that $\lim\limits_{t\to 1-}f(t)=0$. Note also
\begin{equation}
\frac{1}{1-t_j}\int_{t_j}^1(\tg(s)+h(s)+f(s))ds\,\ge\, -\epsilon_j+\frac{1}{1-t_j}\int_{t_j}^{t_{j+1}}(4\,\epsilon_j+2^{-j})ds\ge \epsilon_j+2^{-j-1}>0.
\end{equation}
Now we have $\tg+f\ge0$, $\lim\limits_{k\to\infty}\frac{1}{1-t_k}\int_{t_k}^1\tg(s)+f(s)ds=0$ and $\lim\limits_{t\to 1-}\frac{1}{1-t}\int_t^1h(s)ds=0$. Thus we can apply the method in {\it Case 1}, and conclude that there exists a sequence of times $\tk$ approaching $1$ as $k\to\infty$, with 
\begin{equation}
\limsup_{k\to\infty}\sup_{\sigma\in(0,1-\tk)}\frac{1}{\sigma}\int_{\tk}^{\tk+\sigma}(\tg(s)+h(s)+f(s))ds\leq 0.
\end{equation}
By the nonnegativity of $f$, we then obtain
\begin{equation}
\limsup_{k\to\infty}\sup_{\tau\in(0,1-\tk)}\frac{1}{\tau}\int_{\tk}^{\tk+\tau}(\tg(s)+h(s))ds\leq 0.
\end{equation}
The lemma is proved.\\

\bigskip

For the application to the global existence case, we also need a corresponding version of the above lemma in the global case. 
Let $u$ be a type II global solution to (\ref{eq:mainwaveequation}) and let $u^L$ be the free radiation term. Set $a=u-u^L$, then we have
\begin{eqnarray*}
&&\lim_{t\to\infty}\frac{1}{t}\int_0^t\int_{R^3}(\partial_sa)^2(x,s)\,dxds=0;\\
&&\lim_{t\to\infty}\frac{1}{t}\int_0^t\int_{R^3}|\nabla a|^2-|a|^3(x,s)\,dxds=0.
\end{eqnarray*}
Similarly let $\psi$ be a finite energy global solution to (\ref{eq:WP}), and let $\phi^L$ be the free radiation term. Set $a=\psi-\phi^L$, then we have
\begin{eqnarray*}
&&\lim_{t\to\infty}\frac{1}{t}\int_0^t\int_0^{\infty}(\partial_sa)^2(r,s)\,rdrds=0;\\
&&\lim_{t\to\infty}\frac{1}{t}\int_0^t\int_0^{\infty}\left(f'(a)(\partial_ra)^2+\frac{f^2(a)}{r^2}\right)\,rdrds=0.
\end{eqnarray*}
Set
\begin{eqnarray*}
&&g(t)=\int_{R^3}(\partial_ta)^2(x,t)\,dx,\\
&&h(t)=\int_{R^3}|\nabla a|^2-|a|^3(x,t)\,dx,
\end{eqnarray*}
in the case of (\ref{eq:mainwaveequation}) and set
\begin{eqnarray*}
&&g(t)=\int_0^{\infty}(\partial_ta)^2(r,t)\,rdr,\\
&&h(t)=\int_0^{\infty}\left(f'(a)(\partial_ra)^2+\frac{f^2(a)}{r^2}\right)\,rdr,
\end{eqnarray*}
in the case of (\ref{eq:WP}). Then $g\ge 0$, $g,\,h$ are bounded on $(0,\infty)$, and satisfy
\begin{eqnarray}
&&\lim_{t\to\infty}\frac{1}{t}\int_0^tg(s)\,ds=0,\label{eq:controlg}\\
&&\lim_{t\to\infty}\frac{1}{t}\int_0^th(s)\,ds=0.\label{eq:controlh}
\end{eqnarray}
Our goal is to prove the following statement.
\begin{lemma}\label{lm:keylemmaappendix}
Let $g\ge 0$, $g,\,h$ be bounded on $(0,\infty)$, and satisfy (\ref{eq:controlg}), (\ref{eq:controlh}). 
Then there exists a positive increasing function $\lambda(t)$ on $(0,\infty)$ with $\lim\limits_{t\to\infty}\lambda(t)=\infty$, and a sequence of times $\tilde{t}_k\to\infty$ as $k\to\infty$, such that
\begin{equation}\label{eq:vanishingaverageIIappendix}
\limsup_{k\to\infty}\sup_{\tau\in (0,\frac{\tk}{4})}\frac{1}{\tau}\int_{\tk}^{\tk+\tau}\left(\lambda(s)g(s)+h(s)\right)ds\leq 0.
\end{equation}
\end{lemma}

\smallskip
\noindent
{\it Proof.} Firstly we claim that there exists a positive increasing $\lambda(t)$ with $\lambda(t)\to\infty$ as $t\to\infty$, such that
\begin{equation}
\lim_{k\to\infty}\frac{1}{t_k}\int_0^{t_k}\lambda(s)g(s)\,ds=0,
\end{equation}
for some sequence $t_k\to\infty$ with additional properties (\ref{eq:no1sec5}), (\ref{eq:no2sec5}) below. To prove this claim, we can assume $g$ is not identically zero on $(T,\infty)$ for any $T>0$, otherwise the claim holds trivially. Let us choose a sequence of times $t_k$ with $t_1>t_0=0$ and $\int_0^{t_1}g(t)\,dt>0$, such that
\begin{eqnarray}
&&t_{k+1}>2t_k>0,\,\,{\rm for}\,\,k\ge 1;\label{eq:no1sec5}\\
&&\frac{1}{t_{k+1}-t_k}\int_{t_k}^{t_{k+1}}g(s)\,ds>0,\,\,\,{\rm for}\,\,k\ge 1;\\
&&\frac{1}{t_{k+1}-t_k}\int_{t_k}^{t_{k+1}}g(s)\,ds\leq\frac{1}{4}\frac{1}{t_k-t_{k-1}}\int_{t_{k-1}}^{t_k}g(s)\,ds\,\,{\rm for}\,\,k\ge 1.\label{eq:no2sec5}
\end{eqnarray}
The existence of such a sequence follows easily from our assumption that $g$ is not identically zero on $(T,\infty)$ for any $T>0$, and the fact that
\begin{equation}
\lim_{t\to\infty}\frac{1}{t-T}\int_T^tg(s)\,ds=0
\end{equation}
for any $T>0$, which is an immediate consequence of (\ref{eq:controlg}). Define
\begin{equation}
\lambda(t)=\sum_{k=1}^{\infty}2^k\chi_{[t_k,t_{k+1})}.
\end{equation}
It is clear that $\lambda(t)$ is an increasing function and $\lim\limits_{t\to\infty}\lambda(t)=\infty$. Now
\begin{eqnarray*}
&&\frac{1}{t_m}\int_{t_1}^{t_m}\lambda(s)g(s)\,ds\\
&&=\sum_{k=1}^{m-1}\frac{1}{t_m}2^k\int_{t_k}^{t_{k+1}}g(s)\,ds\\
&&\leq\sum_{k=1}^{m-1}\frac{2^k}{t_m}(t_{k+1}-t_k)\frac{1}{t_{k+1}-t_k}\int_{t_k}^{t_{k+1}}g(s)\,ds\\
&&\leq\sum_{k=1}^{m-1}\frac{2^k}{t_m}(t_{k+1}-t_k)4^{-k}\frac{1}{t_1}\int_{0}^{t_{1}}g(s)\,ds\\
&&\leq C\sum_{k=1}^{m-1}\frac{1}{t_m}2^{-k}(t_{k+1}-t_k).
\end{eqnarray*}
We note that for any integer $k_0>0$,
\begin{eqnarray*}
&&\limsup_{m\to\infty}\,\sum_{k=1}^{m-1}\frac{1}{t_m}2^{-k}(t_{k+1}-t_k)\\
&&=\limsup_{m\to\infty}\,\frac{1}{t_m}\sum_{k=k_0}^{m-1}2^{-k}(t_{k+1}-t_k)\\
&&\leq 2^{-k_0}.
\end{eqnarray*}
Hence we conclude that
\begin{equation}
\limsup_{m\to\infty}\,\frac{1}{t_m}\int_0^{t_m}\lambda(s)g(s)\,ds=0,
\end{equation}
which finishes the proof of our claim. \\
Now set 
\begin{equation}
\tg(t)=\lambda(t)g(t).
\end{equation}
We know that 
\begin{equation}
\tg(t)\ge 0,\,\,\,{\rm and}\,\,\,\lim_{k\to\infty}\frac{1}{t_k}\int_0^{t_k}\tg(s)ds=0,\,\,\,{\rm with}\,\,\,t_{k+1}>2t_k>0.
\end{equation}
Also,
\begin{equation}
\sup_{t\in(0,\infty)}\,|h(t)|\leq CM,\,\,\,{\rm and}\,\,\,\lim_{t\to\infty}\frac{1}{t}\int_0^th(s)ds=0.
\end{equation}
A moment's reflection shows
\begin{equation}\label{eq:convergeto0}
\lim_{k\to\infty}\frac{1}{t_{k+1}-t_k}\int_{t_k}^{t_{k+1}}\tg(t)\,dt=0.
\end{equation}
We now claim that we can find $\tk$ with 
\begin{equation}
\lim_{k\to\infty}\tk=\infty,
\end{equation}
such that
\begin{equation}
\limsup_{k\to\infty}\sup_{\tau\in (0,\frac{\tk}{4})}\frac{1}{\tau}\int_{\tk}^{\tk+\tau}\tg(t)+h(t)\,dt\leq 0.
\end{equation}
The above claim implies our lemma. To prove the claim, we fix any $\epsilon>0$, by (\ref{eq:convergeto0}) we can choose $k$ sufficiently large, such that
\begin{eqnarray}
&&\frac{1}{t_{k+1}-t_k}\int_{t_k}^{t_{k+1}}\tg(t)dt<\frac{\epsilon}{1000},\\
&&{\rm and}\,\,\,\left|\frac{1}{t}\int_0^th(s)\,ds\right|<\frac{\epsilon}{1000},\,\,\,{\rm for}\,\,t\ge t_k.
\end{eqnarray}
Suppose for all $t\in[t_k,t_k+\frac{t_{k+1}-t_k}{2}]$, we can always find some $\tau\in(0,t_{k+1}-t)$ with 
\begin{equation}\label{eq:ineqforcontradiction}
\frac{1}{\tau}\int_t^{t+\tau}\tg(s)+h(s)ds>\epsilon.
\end{equation}
We will get a contradiction. Let us first prove the following useful inequality.\\

Fix $c_1>0$. For any $t>0$ and $\tau>c_1t$, we have
\begin{equation}\label{eq:controlofaverage}
\left|\frac{1}{\tau}\int_t^{t+\tau}h(s)\,ds\right|\leq\frac{2+c_1}{c_1}\left(\left|\frac{1}{t}\int_0^th(s)\,ds\right|+\left|\frac{1}{t+\tau}\int_0^{t+\tau}h(s)\,ds\right|\right).
\end{equation}
The proof follows easily from a contradiction argument, we sketch some of the details for completeness. Suppose (\ref{eq:controlofaverage}) does not hold. Then
\begin{equation*}
\left|\frac{1}{\tau}\int_t^{t+\tau}h(s)\,ds\right|>\frac{2+c_1}{c_1}\left(\left|\frac{1}{t}\int_0^th(s)\,ds\right|+\left|\frac{1}{t+\tau}\int_0^{t+\tau}h(s)\,ds\right|\right).
\end{equation*}
It follows that
\begin{equation*}
\frac{1}{2+c_1}\left|\int_t^{t+\tau}h(s)\,ds\right|>\left|\int_0^th(s)\,ds\right|.
\end{equation*}
Hence
\begin{equation*}
\left|\int_0^{t+\tau}h(s)\,ds\right|\ge\left|\int_t^{t+\tau}h(s)\,ds\right|-\left|\int_0^{t}h(s)\,ds\right|>\frac{1+c_1}{2+c_1}\left|\int_t^{t+\tau}h(s)\,ds\right|.
\end{equation*}
Combining this with $t<\frac{\tau}{c_1}$, we get
\begin{eqnarray*}
&&\left|\frac{1}{t+\tau}\int_0^{t+\tau}h(s)\,ds\right|\\
&&>\frac{1}{(\frac{1}{c_1}+1)\tau}\frac{1+c_1}{2+c_1}\left|\int_t^{t+\tau}h(s)\,ds\right|\\
&&=\frac{c_1}{2+c_1}\left|\frac{1}{\tau}\int_t^{t+\tau}h(s)\,ds\right|,
\end{eqnarray*}
a contradition. Thus the inequality (\ref{eq:controlofaverage}) is proved.\\

Now let us continue to prove the claim. Take 
\begin{equation}
A=\left\{\tau\in(0,t_{k+1}-t_k]: \,\frac{1}{\tau}\int_{t_k}^{t_k+\tau}\tg(s)+h(s)\,ds\ge\epsilon\right\}.
\end{equation}
Then by assumption (\ref{eq:ineqforcontradiction}) with $t=t_k$, we see $A\neq \emptyset$. Note that for $\tau\in[\frac{t_{k+1}-t_k}{4},t_{k+1}-t_k]$, by inequality (\ref{eq:controlofaverage}), we have
\begin{eqnarray*}
&&\frac{1}{\tau}\int_{t_k}^{t_k+\tau}\tg(s)+h(s)ds\\
&&\leq\frac{4}{t_{k+1}-t_k}\int_{t_k}^{t_{k+1}}\tg(s)\,ds+24\left(\left|\frac{1}{t_k}\int_0^th(s)\,ds\right|+\left|\frac{1}{t_{k+1}}\int_0^{t_{k+1}}h(s)\,ds\right|\right)\\
&&\leq4\cdot\frac{\epsilon}{1000}+24\cdot\frac{\epsilon}{500}<\frac{\epsilon}{2}.
\end{eqnarray*}
Thus
\begin{equation}
\theta_k:=\sup A\leq \frac{t_{k+1}-t_k}{4}.
\end{equation}
Moreover, by continuity of $\frac{1}{\tau}\int_t^{t+\tau}\tg(s)+h(s)ds$ in $\tau$, we see $\theta_k\in A$. Thus (again by continuity)
\begin{equation}\label{eq:equalepsilon}
\frac{1}{\theta_k}\int_{t_k}^{t_k+\theta_k}\tg(s)+h(s)ds=\epsilon.
\end{equation}
Since $t_k+\theta_k<t_k+\frac{t_{k+1}-t_k}{2}$, by our assumption (\ref{eq:ineqforcontradiction}) with $t=t_k+\theta_k$, we can find $\tau\in (0,t_{k+1}-t_k-\theta_k)$, such that
\begin{equation}\label{eq:greaterthanepsilon}
\frac{1}{\tau}\int_{t_k+\theta_k}^{t_k+\theta_k+\tau}\tg(s)+h(s)ds>\epsilon.
\end{equation}
By a similar argument as before, we know $\tau<\frac{t_{k+1}-t_k}{4}$. Hence $t_k+\theta_k+\tau<t_k+\frac{t_{k+1}-t_k}{2}$. By (\ref{eq:equalepsilon}) and (\ref{eq:greaterthanepsilon}), we have
\begin{equation}
\frac{1}{\theta_k+\tau}\int_{t_k}^{t_k+\theta_k+\tau}\tg(s)+h(s)ds>\epsilon,\,\,\,{\rm and}\,\,\,\theta_k+\tau\in (0,t_{k+1}-t_k).
\end{equation}
This is a contradiction with the fact $\theta_k=\sup \,A$. Hence our assumption is false, and we can find $\tk\in[t_k,t_k+\frac{t_{k+1}-t_k}{2}]$, such that
\begin{equation}
\sup_{\tau\in(0,t_{k+1}-\tk)}\frac{1}{\tau}\int_{\tk}^{\tk+\tau}\tg(s)+h(s)ds\leq\epsilon.
\end{equation}
Note that $t_{k+1}>2t_k$ implies
\begin{equation}
t_{k+1}-\tk\ge\frac{t_{k+1}-t_k}{2}>\frac{t_{k+1}}{4}>\frac{\tk}{4}.
\end{equation}
Thus the lemma is proved with the time sequence $\tk\to\infty$.



\begin{section}{Appendix: Preliminaries on the Cauchy theory and the profile decomposition for the radial energy critical wave and equivariant wave map equations}
\begin{subsection}{Energy critical focusing semilinear wave equation in $R^6$}
In this section we recall some facts on the equation (\ref{eq:mainwaveequation}) that are used in the paper.
The following result from \cite{RKS} says that essentially the free waves travel at unit speed.
\begin{lemma}\label{lm:energypartition}
Suppose $u^L$ is a radial solution to the linear wave equation in $d$ dimensions, with $\OR{u}^L(0)\in \HL(R^d)$. Then we have
\begin{equation}
\lim_{R\to\infty}\limsup_{t\to\infty}\int_{||x|-t|\ge R}\,|\nabla_{t,x}u^L|^2(x,t)\,dx=0.
\end{equation}
\end{lemma}

For initial data $(u_0,u_1)\in\HL(R^6)$, denote $u^L(\cdot,t)=S(t)(u_0,u_1)$ as the solution to the six dimensional linear wave equation. By Strichartz estimates (see Ginibre-Velo \cite{GV}), we have $u^L\in L^2_tL^4_x(R^6\times(0,\infty))$. We can re-write equation (\ref{eq:mainwaveequation}) in the integral form as
\begin{equation}
u(\cdot,t)=S(t)(u_0,u_1)+\int_0^tS(t-s)|u|u(\cdot,s)\,ds.
\end{equation}
Then using the Strichartz estimates, see e.g. \cite{GV}, we obtain
\begin{equation*}
\left\|\int_0^tS(t-s)f(\cdot,s)\,ds\right\|_{C([0,\infty),\HL(R^6))\cap L^2_tL^4_x(R^6\times (0,\infty))}\leq C\|f\|_{L^1_tL^2_x(R^6\times (0,\infty))},
\end{equation*}
we can apply the Contraction Mapping Theorem to obtain a unique solution $u\in C([0,T),\HL(R^6))\cap L^2_tL^4_x(R^6\times (0,T))$ on $R^6\times [0,T)$, provided $\|u^L\|_{L^2_tL^4_x(R^6\times (0,T))}$ is sufficiently small (which can be achieved by making $T$ sufficiently small). This solution can then be continued to a maximal interval of existence $[0,T_+)$, and one also obtain a ``blow-up" criteria: $T_+<\infty$ if and only if $\|u\|_{L^2_tL^4_x(R^6\times (0,T_+))}=\infty$. We refer the reader to \cite{KenigYang} for full details. If $\|(u_0,u_1)\|_{\HL(R^6)}$ is sufficiently small the solution is global and scatters. Moreover, the following energy is preserved for $t\in [0,T_+)$
\begin{equation}
\E(\OR{u}(t)):=\int_{R^6}\left(\frac{|\partial_{t,x}u|^2}{2}-\frac{|u|^3}{3}\right)(x,t)\,dx\equiv \E(u_0,u_1).
\end{equation}
When the initial data is not assumed to be small, more interesting behavior may occur. On the one hand, by finite speed of propagation and study of the ODE $y''=|y|y$, it's clear that there exists a solution $\OR{u}$ to equation (\ref{eq:mainwaveequation}) which blows up in finite time. On the other hand, the ground state $W=\frac{1}{(1/24+|x|^2)^2}$ exists globally and does not  scatter. The ground state is charaterized as the minimizer for the best constant in the Sobolev inequality
\begin{equation}
\|f\|_{L^3(R^6)}\leq C \|\nabla f\|_{L^2(R^6)},
\end{equation}
and plays an important role in the global dynamics of equation (\ref{eq:mainwaveequation}). From \cite{kenigmerle} in the case of $R^3$, we know if $\E(u_0,u_1)<\E(W,0)$, then we have the following dichotomy:
\begin{itemize}
\item $\|\nabla u_0\|_{L^2(R^3)}<\|\nabla W\|_{L^2(R^3)}$, then the solution $u$ is global and scatters as $t\to\infty$;
\item $\|\nabla u_0\|_{L^2(R^3)}>\|\nabla W\|_{L^2(R^3)}$, the solution $u$ blows up in finite time.
\end{itemize}
The case $\|\nabla u_0\|_{L^2(R^3)}=\|\nabla W\|_{L^2(R^3)}$ can not happen. Similar results hold in $R^6$ (see \cite{KenigYang}).\\

An important ingredient in \cite{kenigmerle} (and in our analysis) is the Bahouri-Gerard profile decompositions, introduced in the context of wave equations in \cite{BaGe}, (see also \cite{Bulut} in the case of $R^{1+6}$). We present some of the details needed below in the simpler radial case.\\

Suppose $(u_{0,n},u_{1,n})\in \HL(R^6)$ is a sequence of radial initial data with 
\begin{equation*}
\sup\limits_n \|(u_0,u_1)\|_{\HL(R^6)}<\infty.
\end{equation*}
By passing to a subsequence, we can assume that $(u_{0,n},u_{1,n})$ has the following profile decomposition
\begin{equation}\label{eq:linearprofiledecomposition}
(u_{0,n},u_{1,n})=\sum_{j=1}^J\left(\frac{1}{\lambda_{jn}^2}U_L^j(\frac{x}{\lambda_{jn}},\frac{-t_{jn}}{\lambda_{jn}}),\frac{1}{\lambda_{jn}^3}\partial_tU_L^j(\frac{x}{\lambda_{jn}},\frac{-t_{jn}}{\lambda_{jn}})\right)+(w^J_{0,n},w^J_{1,n}),\,\,\,\forall J\ge 1,
\end{equation}
where $U_L^j(x,t)$ are radial finite energy solutions to the six dimensional linear wave equation. In addition, let $w^J_{n}(x,t)$ be the radial  solution to linear wave equation with initial data $(w^J_{0,n},w^J_{1,n})$, then $w^J_{n}$, $\lambda_{jn}>0$ and $t_{jn}\in R$ satisfy the following pseudo-orthogonality properties
\begin{eqnarray*}
&&t_{jn}\equiv 0\,\,{\rm for\,\,all\,\,}n,\,\,{\rm or}\,\,\lim_{n\to\infty}\frac{t_{jn}}{\lambda_{jn}}\in\{\pm\infty\};\\
&&\lim_{n\to\infty}\left(\frac{\lambda_{jn}}{\lambda_{j'n}}+\frac{\lambda_{j'n}}{\lambda_{jn}}+\frac{|t_{jn}-t_{j'n}|}{\lambda_{jn}}\right)=\infty, \,\,\,{\rm for\,\,all}\,\,j\neq j';\\
&&{\rm write}\,\,w^J_{n}(x,t)=\frac{1}{\lambda_{jn}^2}\widetilde{w}^j_{Jn}(\frac{x}{\lambda_{jn}},\,\frac{t-t_{jn}}{\lambda_{jn}}),\,\,{\rm then}\,\,\widetilde{w}^j_{Jn}\rightharpoonup 0,\,\,{\rm as}\,\,n\to\infty, \,\,{\rm for}\,\,1\leq j\leq J.
\end{eqnarray*}
Moreover, the dispersive part $w^J_{n}$ satisfies the bound
\begin{equation}
\lim_{J\to\infty}\limsup_{n\to\infty}\|w^J_{n}\|_{L^{\infty}_tL^3_x(R^6\times[0,\infty))\cap L^2_tL^4_x(R^3\times [0,\infty))}=0.
\end{equation}
Notice that by the pseudo-orthogonality of the parameters, the different profiles are ``supported" in different length scales or at different times (for example, when one looks at where the Strichartz norm of each profile is concentrated), hence one can expect that different profiles should not interact with each other very much, even in the nonlinear evolution. 
To make this idea more rigorous, let us introduce the ``nonlinear profiles" $U^j$ defined as follows
\begin{itemize}
\item If $\lim\limits_{n\to\infty}\frac{-t_{jn}}{\lambda_{jn}}=L\in\{\pm\infty\}$, then define $U^j$ as the solution to equation (\ref{eq:mainwaveequation}), with
\begin{equation}
\lim_{t\to L }\|\OR{U}^j(\cdot,t)-\OR{U}^j_L(\cdot,t)\|_{\HL(R^6)}=0.
\end{equation}
The existence of $U^j$ in a neighborhood of $L$ satisfying these conditions can be proved by standard perturbation arguments. Then 
\begin{equation}
U^j_n:=\frac{1}{\lambda^2_{jn}}U^j(\frac{x}{\lambda_{jn}},\frac{t-t_{jn}}{\lambda_{jn}})
\end{equation}
is a solution that approximates the solution to equation (\ref{eq:mainwaveequation}) with initial data $\frac{1}{\lambda_{jn}^2}\OR{U}_L^j(\frac{x}{\lambda_{jn}},-\frac{t_{jn}}{\lambda_{jn}})$;
\item If $t_{jn}=0$, then define $U^j$ simply as the solution to equation (\ref{eq:mainwaveequation}) with initial data $\OR{U}_L^j(x,0)$. In this case 
\begin{equation}
U^j_n:=\frac{1}{\lambda^2_{jn}}U^j(\frac{x}{\lambda_{jn}},\frac{t}{\lambda_{jn}})
\end{equation}
is exactly the solution to equation with initial data $\frac{1}{\lambda_{jn}^2}U_L^j(\frac{x}{\lambda_{jn}},0)$.
\end{itemize}
By suitable perturbation arguments, one can then prove the following nonlinear version of the profile decomposition (see \cite{DKMsmall}):\\
Let $(u_{0,n},u_{1,n})$ be a sequence of  radial initial data with uniformly bounded $\HL(R^6)$ norm. Suppose $(u_{0,n},u_{1,n})$ has the profile decomposition (\ref{eq:linearprofiledecomposition}). Let $u_n$ be the solution to equation (\ref{eq:mainwaveequation}) with initial data $(u_{0,n},u_{1,n})$. Let $\theta_n>0$ be such that 
\begin{equation*}
\forall j\ge 1, \,\,\frac{\theta_n-t_{jn}}{\lambda_{jn}}<T_+(U^j), 
\end{equation*}
for sufficiently large $n$, where $T_+(U^j)$ is the maximal time of existence for $U^j$, \footnote{Note that $T_+(U^j)=\infty$ if $\lim\limits_{n\to\infty}\frac{-t_{jn}}{\lambda_{jn}}=\infty$, and if $\lim\limits_{n\to\infty}\frac{-t_{jn}}{\lambda_{jn}}=-\infty$ then $U^j$ exists on $(-\infty,T_+(U^j))$. In both cases $\OR{U}^j$ scatters at $\lim\limits_{n\to\infty}\frac{-t_{jn}}{\lambda_{jn}}$.} and 
\begin{equation*}
\forall j,\,\,\,\limsup_{n\to\infty}\|U^j\|_{L^2_tL^4_x\left(R^6\times (-\frac{t_{jn}}{\lambda_{jn}},\frac{\theta_n-t_{jn}}{\lambda_{jn}})\right)}<\infty,
\end{equation*}
then for all $n$ sufficiently large, $u_n$ is defined on $[0,\theta_n]$ and we have the following nonlinear profile decomposition
\begin{equation}
\OR{u}_n(x,t)=\sum_{j=1}^J\OR{U}^j_n(x,t)+\OR{w}_{Jn}(t,x)+\OR{r}_{Jn}(x,t),\,\,\,\forall J\ge 1,
\end{equation}
where the residue term $\OR{r}_{Jn}$ vanishes in the sense that
\begin{equation}
\lim_{J\to\infty}\limsup_{n\to\infty}\left(\|r_{Jn}\|_{L^2_tL^4_x(R^6\times [0,\theta_n])}+\sup_{t\in[0,\theta_n]}\|\OR{r}_{Jn}(t)\|_{\HL(R^6)}\right)=0.
\end{equation}

We shall also need the following elementary real analysis lemma.
\begin{lemma}\label{lm:restriction}
Let $\Sigma=\{(s,t): t\ge 0,\,s\in R\}$. Suppose $f\in L^1_{{\rm loc}}(\Sigma)$ and $\partial_tf\in L^1_{{\rm loc}}(\Sigma)$, then $f(\cdot,0)\in L^1_{{\rm loc}}(R)$ is well defined. Moreover, assume a sequence of functions $f_n$ verify $f_n\in L^1_{{\rm loc}}(\Sigma)$, $\partial_tf_n\in L^1_{{\rm loc}}(\Sigma)$, with $f_n\to f$ in $L^1_{{\rm loc}}(\Sigma)$ and $\partial_tf_n\to\partial_tf$ in $L^1_{{\rm loc}}(\Sigma)$. Then we have $f_n(\cdot,0)\to f(\cdot,0)$ in $L^1_{{\rm loc}}(R)$.
\end{lemma}

\smallskip
\noindent
{\it Proof.} Take a smooth function $\eta\in C_c^{\infty}([0,\infty))$ with $\eta\equiv 1$ on $[0,1]$. For any $K\Subset R$ and continous function $\phi$ with ${\rm supp}\,\phi\Subset K$, we can define $f(\cdot,0)$ as a distribution with
\begin{equation}\label{eq:boundaryvalue1}
(f(\cdot,0),\,\phi):=-\int_{\Sigma}f(s,t)\phi(s)\eta'(t)\,dtds-\int_{\Sigma}\partial_tf(t,s)\phi(s)\eta(t)\,dtds.
\end{equation}
This definition clearly coincides with classical boundary values if $f$ is smooth on $\Sigma$. It's easy to see $f(\cdot,0)$ thus defined is a locally finite measure on $R$. Note also that as $|K|\to 0$ while keeping $\|\phi\|_{C(K)}\leq 1$, then $(f(\cdot,0),\,\phi)\to 0$. Hence $f(\cdot,0)$ thus defined is actually in $ L^1_{{\rm loc}}(R)$. The second part of the lemma follows directly from the definition (\ref{eq:boundaryvalue1}).
\end{subsection}

\begin{subsection}{Two dimensional equivariant wave maps}
In this section we review some results for the two dimensional equivariant wave map equation (\ref{eq:WP}).
Let $(\psi_0,\psi_1)$ be of finite energy
\begin{equation*}
\mathcal{E}(\psi_0,\psi_1)=\int_0^{\infty}\left((\partial_r\psi_0)^2+\psi_1^2+\frac{g^2(\psi_0)}{r^2}\right)\,rdr<\infty.
\end{equation*}

Clearly $\psi_0$ is locally H\"older continuous in $(0,\infty)$. Moreover $\psi_0$ has limits both as $r\to\infty$ and $r\to 0$, which is a consequence of the next lemma.
\begin{lemma}\label{lm:boundonpsi}
Let $\psi_0$ be as above, then $\psi_0$ is bounded. Moreover $\psi_0(0):=\lim\limits_{r\to 0+}\psi_0(r)$ and $\psi_0(\infty):=\lim\limits_{r\to\infty}\psi_0(r)$ exist and are zeros of $g$. 
\end{lemma}

\smallskip
\noindent
{\it Proof.} By H\"older inequality, we obtain 
\begin{equation}\label{eq:boundaryvalue}
|G(\psi_0(r_1))-G(\psi_0(r_2))|\leq \int_{r_1}^{r_2}|g(\psi_0(r))||\partial_r\psi_0|\,dr\leq \E(\psi_0,r_1,r_2).
\end{equation}
By $\E(\psi_0,r_1,r_2)<\infty$ and assumption (A1), we see that $\psi_0$ is bounded. Moreover, since both as $r_1\to 0,\,r_2\to 0$ and as $r_1\to\infty,\,r_2\to\infty$, the right hand side of $(\ref{eq:boundaryvalue})$ tends to zero, we can conclude that $\psi_0(0)$ and $\psi_0(\infty)$ are well defined. The fact that $\psi_0(0)$ and $\psi_0(\infty)$ are zeros of $g$ follows from the finiteness of the energy $\E(\psi_0,\psi_1)<\infty$.\\

The following lemma shows that whenever $\psi_0(r)$ is away from the zeros of $g(\psi)$, there is some concentration of energy near $r$.
\begin{lemma}\label{lm:concentrationofenergy}
Suppose $M:=\E(\psi_0,\psi_1)<\infty$, $\psi_0(\infty)=m\in\mathcal{V}$ and $r_0>0$. Assume
\begin{equation}
\delta=\inf_{l\in \mathcal{V}}\,|\psi_0(r_0)-l|>0.
\end{equation}
Then for any $\gamma>1$, there exists $\epsilon=\epsilon(\delta,\gamma,m,M,g)>0$, such that
\begin{equation}\label{eq:concentrationofenergy}
\int_{\frac{r_0}{\gamma}}^{\gamma r_0}\left((\partial_r\psi_0)^2+\frac{g^2(\psi_0)}{r^2}\right)\,rdr\ge\epsilon.
\end{equation}
\end{lemma}

\smallskip
\noindent
{\it Proof.} By H\"{o}lder inequality, for $\frac{r_0}{\gamma}<r_1<r_2<\gamma r_0$, we have 
\begin{eqnarray*}
|\psi_0(r_1)-\psi_0(r_2)|&=&\left|\int_{r_1}^{r_2}\partial_r\psi_0\,dr\right|\\
&\leq&C\left(\int_{r_1}^{r_2}(\partial_r\psi_0)^2\,rdr\right)^{\frac{1}{2}} \log^{\frac{1}{2}}{\frac{r_2}{r_1}}\leq C(\gamma)\left(\int_{r_1}^{r_2}(\partial_r\psi_0)^2\,rdr\right)^{\frac{1}{2}}.
\end{eqnarray*}
Hence, if 
\begin{equation*}
|\psi_0(r_1)-\psi_0(r_2)|>\frac{\delta}{2}
\end{equation*}
for some $r_1,\,r_2\in (\frac{r_1}{\gamma},\gamma r_0)$, we would have
\begin{equation*}
 \int_{\frac{r_0}{\gamma}}^{\gamma r_0}(\partial_r\psi_0)^2\,rdr\ge 
\frac{\delta^2}{4C(\gamma)}
\end{equation*}
and this finishes proves our claim in this case. Otherwise we have
\begin{equation*}
\sup_{r\in (\frac{r_0}{\gamma},\gamma r_0)}|\psi_0(r)-\psi_0(r_0)|\leq \frac{\delta}{2}.
\end{equation*}
Hence 
\begin{equation*}
\inf_{l\in \mathcal{V}}|\psi_0(r)-l|\ge\frac{\delta}{2},\,\,\,{\rm for\,\,any\,\,}r\in(\frac{r_0}{\gamma},\gamma r_0).
\end{equation*}
By Lemma \ref{lm:boundonpsi}, $\psi_0$ is pointwisely bounded by $C(M,m,g)$, hence using the local smoothness of $g$, we get
\begin{equation*}
\int_{r\in (\frac{r_0}{\gamma},\gamma r_0)}\frac{g^2(\psi_0)}{r^2}\,rdr\ge \epsilon(\gamma,\delta,M,m,g)>0.
\end{equation*}
The lemma is proved.\\

Let us recall the profile decomposition adapted to the equivariant wave map setting, introduced by C\^ote in \cite{cotesoliton}. We will give an alternative treatment at the end of this Appendix below.

\begin{lemma}\label{lm:decompositionforwavemap}
Let $(\psi_{0,n},\psi_{1,n})$ be a sequence with $M:=\sup\limits_{n\ge 1}\,\E(\psi_{0,n},\psi_{1,n})<\infty$ and $\psi_{0,n}(\infty)=m\in\mathcal{V}$. Passing to a subsequence (which we still denote as $(\psi_{0n},\psi_{1n})$ for ease of notations) if necessary, we can assume $\psi_{0,n}(0)$ is constant in $n$ with $k=|g'(\psi_{0,n}(0))|$ and there exist an integer $K\ge 0$,  sequences $r_{ln}>0$ for $1\leq l\leq K$, $\lambda_{jn}>0$, $t_{jn}\in R$, and functions $(\psi^l_0,0)$ with $\E(\psi^l_0,0)<\infty$, finite energy radial solutions $U^j_L$ to the $2k+2$ dimensional linear wave equation, such that we have the following profile decomposition for each $J\ge 1$
\begin{equation}\label{eq:wavemapdecomposition11}
\def\arraystretch{2.2}
\begin{array}{rl}
(\psi_{0,n},\psi_{1,n})=&\left(\sum_{l=1}^{K}(\psi^l_0(\frac{\cdot}{r_{ln}})-\psi^l_0(0))+\psi_0^K(0),0\right)\\
&+r^k\sum_{j=1}^J\left(\frac{1}{\lambda^k_{jn}}U^j_L(\frac{r}{\lambda_{jn}},\frac{-t_{jn}}{\lambda_{jn}}),\frac{1}{\lambda_{jn}^{k+1}}\partial_tU^j_L(\frac{r}{\lambda_{jn}},\frac{-t_{jn}}{\lambda_{jn}})\right)+r^k(w^J_{0,n},w^J_{1,n}).
               \end{array}
\end{equation}
If we let $w^J_n$ be the solution to the $2k+2$ dimensional wave equation with initial data $(w^J_{0,n},w^J_{1,n})$, then the parameters satisfy ``pseudo-orthogonality" conditions 
\begin{eqnarray*}
&&r_{Kn}\ll r_{K-1,n}\ll\cdots\ll r_{1n};\\
&&t_{jn}\equiv 0\,\,{\rm for\,\,all\,\,}n,\,\,{\rm or}\,\,\lim_{n\to\infty}\frac{t_{jn}}{\lambda_{jn}}\in\{\pm\infty\};\\
&&{\rm if}\,\,t_{jn}\equiv 0\,\,{\rm for\,\,all\,\,}n, \,\,{\rm then}\,\,\frac{\lambda_{jn}}{r_{ln}}+\frac{r_{ln}}{\lambda_{jn}}\to\infty,\,\,{\rm as}\,\,n\to\infty \,\,{\rm for\,\,each\,\,}j\,\,{\rm and}\,\,1\leq l\leq K;\\
&&\lim_{n\to\infty}\left(\frac{\lambda_{jn}}{\lambda_{j'n}}+\frac{\lambda_{j'n}}{\lambda_{jn}}+\frac{|t_{jn}-t_{j'n}|}{\lambda_{jn}}\right)=\infty, \,\,\,{\rm for\,\,all}\,\,j\neq j'.
\end{eqnarray*}
Furthermore,
\begin{eqnarray*}
&&\psi^1_0(\infty)=m,\,\,\psi^{l+1}_0(\infty)=\psi^l_0(0),\,\,{\rm for\,\,each\,\,} 1\leq j\leq K-1,\\
&&{\rm write}\,\,w^J_{n}(x,t)=\frac{1}{\lambda_{jn}^k}\widetilde{w}^j_{Jn}(\frac{x}{\lambda_{jn}},\,\frac{t-t_{jn}}{\lambda_{jn}}),\,\,{\rm then\,\,we\,\,have}\,\,\\
&&\widetilde{w}^j_{Jn}(\cdot,t)\rightharpoonup 0,\,\,{\rm in}\,\,\HL(R^{2k+2})\,\,{\rm as}\,\,n\to\infty,\,\,\forall t\in R \,\,{\rm and\,\, for}\,\,1\leq j\leq J;\\
&&{\rm write}\,\,\widetilde{w}^l_{Jn}=r_{ln}^kw^J_{0,n}(r_{ln}\cdot),\,\,{\rm then\,\,we\,\,have}\,\,\\
&&\widetilde{w}^l_{Jn}\rightharpoonup 0,\,\,{\rm in}\,\,\dot{H}^1(R^{2k+2})\,\,{\rm as}\,\,n\to\infty,\,\,\forall t\in R,\,\,1\leq l\leq K\,\,{\rm and\,\,each\,\,}J;\\
&&\lim_{J\to\infty}\limsup_{n\to\infty}\left(\|w^J_{n}\|_{L^{\infty}_tL^{\frac{2k+2}{k}}_x\cap L^{\frac{k+2}{k}}_tL^{\frac{2k+4}{k}}_x(R^{2k+2}\times [0,\infty))}+\|r^kw^J_{n}\|_{L^{\infty}}\right)=0.
\end{eqnarray*}
\end{lemma}

\smallskip
\noindent
{\it Remark.} If we know in addition that
\begin{equation}\label{vanishing111}
\int_0^{\infty}(\psi_{1,n})^2\,rdr\to 0,\,\,\,{\rm as}\,\,n\to\infty,
\end{equation}
which implies $\|u_{1,n}\|_{L^2(R^{2k+2})}=o_n(1)$ if we set $u_{1,n}=r^k\psi_{1,n}$, then we can get \footnote{See Lemma 5.9 in \cite{Cotemap1} and the arXiv version, 1209.3682v2 for a corrected proof.} the following simplified profile decomposition for the sequence $(\psi_{0,n},\psi_{1,n})$
\begin{eqnarray}
(\psi_{0,n},\psi_{1,n})&=&\left(\sum_{l=1}^{K}(\psi^l_0(\frac{\cdot}{r_{lh}})-\psi^l_0(0))+\psi_0^K(0),0\right)+\nonumber\\
                               &&\quad +r^k\sum_{j=1}^J\left(\frac{1}{\lambda^k_{jn}}U^j_L(\frac{r}{\lambda_{jn}},0),0\right)+(\zeta^J_{0,n},\zeta^J_{1,n}),\label{eq:improveddecomppre}
\end{eqnarray}
with
\begin{equation}
\lim_{J\to\infty}\limsup_{n\to\infty}\|\zeta^J_{0,n}\|_{L^{\infty}}+\lim_{n\to\infty}\|\zeta^J_{1,n}\|_{L^2(R_+,\,rdr)}= 0.
\end{equation}
In addition,
\begin{equation}
\left(\zeta^J_{0,n}(r_{ln}\cdot,),\,r_{ln}\zeta^J_{1,n}(r_{ln}\cdot)\right)\rightharpoonup 0,\,\,\forall 1\leq l\leq K,\,\,{\rm and}\,\,\left(\zeta^J_{0,n}(\lambda_{jn}\cdot),\,\lambda_{jn}\zeta_{1,n}^J(\lambda_{jn}\cdot)\right)\rightharpoonup 0,\,\,\forall 1\leq j\leq J,
\end{equation}
as $n\to\infty$.
The decomposition (\ref{eq:improveddecomppre}) is obtained from (\ref{eq:wavemapdecomposition11}) by absorbing the profiles with $1\leq j\leq J$ and $\lim_{n\to\infty}\frac{t_{jn}}{\lambda_{jn}}\in\{\pm \infty\}$ into the error term $(\zeta^J_{0,n},\,\zeta^J_{1,n})$. The profiles with $t_{jn}\equiv 0$ must satisfy $\partial_tU^j_L(\cdot,0)\equiv 0$ by (\ref{vanishing111}) and Lemma 5.9 in \cite{Cotemap1} (see also the arXiv version, 1209.3682v2 for a correction). 
One important point of the decomposition (\ref{eq:improveddecomppre}) for us is that $\psi_{0,n}(r_{ln}\cdot)\rightharpoonup \psi^l_0$, and $\psi_{0,n}(\lambda_{jn}\cdot)\rightharpoonup U^j_L(\cdot,0)$ as $n\to\infty$, which allows us to capture the correct length scales.

\bigskip

Next we give a version of Struwe's convergence result for wave maps\cite{Struwe2}. On the one hand this convergence is slightly weaker than that given by \cite{Struwe2}, (already sufficient for our purpose though) in the sense that the convergence is only for $r>0$, not across origin, while on the other hand it does not require a ``smallness" condition on the energy near $r=0$, hence is more robust.\\
\begin{lemma}\label{lm:struwe}
Suppose $\alpha_n,\,\beta_n,\in(0,\infty)$ satisfy $\lim\limits_{n\to\infty}\alpha_n=0$ and $\lim\limits_{n\to\infty}\beta_n=\infty$. Let $M:=\sup_n\E(\psi_{0,n},\psi_{1,n})<\infty$. Let $\psi_n$ be a sequence of solutions to (\ref{eq:WP}) on $(0,\infty)\times [0,T)$, and $\psi_n(\infty,t)=m\in\mathcal{V}$ for all $t\in [0,T)$. Suppose in addition
\begin{equation}\label{eq:timeaveragemap}
\lim_{n\to\infty}\int_0^T\int_{\alpha_n}^{\beta_n}(\partial_t\psi_n)^2\,rdrds=0.
\end{equation}
Then we can extract a subsequence of $\psi_n$ and a harmonic map $Q$, such that $\partial_{r,t}\psi_n\to \partial_{r,t}Q$ in $L^2_{{\rm loc}}((0,\infty)\times [0,T))$, and $\psi_n\to Q$  locally uniformly in $(0,\infty)\times[0,T)$. Moreover, we have $\partial_{t,r}(\psi_{0,n},\psi_{1,n})\to \partial_{t,r}(Q,0)$ in $L^2_{{\rm loc}}(0,\infty)$.
\end{lemma}

\smallskip
\noindent
{\it Proof.} By uniform boundedness of energy of $\psi_n$, we conclude $\sup_n\|\psi_n\|_{L^{\infty}}<C(m,M)$. We can select a subsequence of $\psi_n$, such that $\psi_n\rightharpoonup \psi$ in the sense of distributions in $R^2\times [0,T)$. By (\ref{eq:timeaveragemap}), we see $\psi(r,t)=\psi(r)$. Since $\psi_n$ is locally uniformly H\"older continuous in $(0,\infty)\times[0,T)$ by finiteness of energy, we can also assume $\psi_n\to \psi$ locally uniformly in $(0,\infty)\times[0,T)$. Hence 
\begin{equation}
\E(\psi,\partial_t\psi)\leq \liminf_{n\to\infty}\E(\psi_n,\partial_t\psi_n).
\end{equation}
Since $\psi_n$ satisfies equation (\ref{eq:WP}) and $\partial_t\psi_n\to 0$, we see $\psi$ verifies for $r>0$
\begin{equation}
-\partial_{rr}\psi-\frac{\partial_r\psi}{r}+\frac{f(\psi)}{r^2}=0.
\end{equation}
Hence $\psi=Q$ is a harmonic map (see \cite{Cote}). $\psi-\psi_n$ verifies
\begin{equation}\label{eq:differenceequation}
\partial_{tt}\psi_n-\partial_{rr}(\psi_n-\psi)-\frac{\partial_r(\psi_n-\psi)}{r}+\frac{f(\psi_n)-f(\psi)}{r^2}=0.
\end{equation}
Take $\eta\ge 0$, ${\rm supp}\,\eta \Subset (0,\infty)\times [0,T)$, and multiply equation (\ref{eq:differenceequation}) by $(\psi_n-\psi)\eta$ and integrate, we get 
\begin{eqnarray*}
&&\int_0^T\int_0^{\infty}-\partial_t\psi_n\,(\psi_n-\psi)\,\partial_t\eta-(\partial_t\psi_n)^2\eta\,rdrdt-\int_0^{\infty}\psi_{0,n}(\psi_n-\psi)(r,0)\,rdr+\\
&&+\int_0^T\int_0^{\infty}\partial_r(\psi_n-\psi)\cdot(\psi_n-\psi)\partial_r\eta+
(\partial_r(\psi_n-\psi))^2\eta\,rdrdt\\
&&-\int_0^T\int_0^{\infty}\frac{(\partial_r(\psi_n-\psi))\cdot(\psi_n-\psi)}{r}\eta\,rdrdt+\int_0^T\int_0^{\infty}\frac{f(\psi_n)-f(\psi)}{r^2}(\psi_n-\psi)\eta\,rdrdt=0.
\end{eqnarray*}
Letting $n\to\infty$ in the above equation, we get
\begin{equation*}
\lim_{n\to\infty}\int_0^T\int_0^{\infty}\left(\partial_r(\psi_n-\psi)\right)^2\eta \,rdrdt=0,
\end{equation*}
which implies the first claim of the lemma. Now let us prove the last claim of the lemma which is a consequence of the decay of outer energy. More precisely, fix $0<r_1<r_2$,  and any $0<S<r_1$. Multiplying equation by $\partial_t(\psi-\psi_n)$ and integrating in the region $\Sigma:=\{(r,t):\,r_1-t<r<r_2+t,\,t<S\}$, by the convergence 
\begin{equation*}
\left|\frac{f(\psi)-f(\psi_n)}{r^2}\right|\leq C(r_1,f)|\psi-\psi_n|\to 0\,\,\,{\rm as}\,\,n\to\infty, 
\end{equation*}
we obtain
\begin{eqnarray*}
&&\int_{r_1}^{r_2}(\partial_t\psi-\partial_t\psi_n)^2+(\partial_r\psi-\partial_r\psi_n)^2(r,0)\,rdr\\
&&\leq \int_{r_1-S}^{r_2+S}\left((\partial_t\psi-\partial_t\psi_n)^2+(\partial_r\psi-\partial_r\psi_n)^2\right)(r,S)\,rdr\\
&&\quad\quad+C\int_{\Sigma}\frac{|\psi-\psi_n|}{r^2}|\partial_t(\psi-\psi_n)|\,rdrdt\\
&&\leq \int_{r_1-S}^{r_2+S}\left((\partial_t\psi-\partial_t\psi_n)^2+(\partial_r\psi-\partial_r\psi_n)^2\right)(r,S)\,rdr+o_n(1),
\end{eqnarray*}
uniformly in $S<r_1$. Suppose for some positive $\epsilon$, 
\begin{equation}
\limsup_{n\to\infty}\int_{r_1}^{r_2}(\partial_t\psi-\partial_t\psi_n)^2+(\partial_r\psi-\partial_r\psi_n)(r,0)^2\,rdr\ge \epsilon.
\end{equation}
Then for $n$ sufficiently large (again independent of $S$)
\begin{equation*}
\int_{r_1-S}^{r_2+S}(\partial_t\psi-\partial_t\psi_n)^2+(\partial_r\psi-\partial_r\psi_n)^2(r,S)\,rdr\ge \frac{\epsilon}{2}.
\end{equation*}
This holds for any positive $S<r_1$, which is a contradiction to the fact that $\partial_{r,t}\psi_n\to\partial_{r,t}\psi$ in $L^2_{{\rm loc}}((0,T)\times (0,\infty))$ as $n\to\infty$. Hence the lemma is proved.\\

\smallskip
\noindent
{\it Remark.} $Q$ in the above lemma may be trivial. To obtain a nontrivial harmonic map in the limit, one needs to combine the Struwe convergence result with the profile decomposition, which allows us to capture the correct ``length scales". \\

We can use Lemma \ref{lm:decompositionforwavemap} and Lemma \ref{lm:struwe}, together with the asymptotic vanishing of kinetic energy of the wave map after subtracting the regular part in the finite time blow up case or the free radiation part in the global existence case in Section 3, to give an alternative proof of the following Theorem, which is an important first step in establishing the full soliton resolution decomposition for a sequence of times in \cite{cotesoliton}.
\begin{theorem}\label{th:cotesec2}
Let $\OR{\psi}(t)$ be a finite energy solution to equation (\ref{eq:WP}) with initial data $(\psi_0,\psi_1)$. Then there exists a sequence of times $t_n\uparrow T_+$, an integer $J\ge 0$, $J$ sequences of scales $0<r_{Jn}\ll \cdots\ll r_{2n}\ll r_{1n}$ and $J$ harmonic maps $Q_1,\cdots,Q_J$ such that
\begin{equation*}
Q_J(0)=\psi_0(0),\,\,\,Q_{j+1}(\infty)=Q_j(0),\,\,\,{\rm for}\,\,j=1,\cdots,J-1,
\end{equation*}
and such that the following holds.\\
\quad\quad (1)\,\, (Global case) If $T_+=\infty$, denote $\ell=\psi_0(\infty)$. Then $Q_1(\infty)=\ell$, $r_{1n}\ll t_n$ and there exists a radial finite energy solution $\OR{\phi}^L$ to the $2|g'(\ell)|+2$ dimensional linear wave equation such that 
\begin{equation}
\OR{\psi}(t_n)=\sum_{j=1}^J(Q_j(\cdot/r_{jn})-Q_j(\infty),0)+(\ell,0)+\OR{\phi}^L(t_n)+\OR{b}_n.
\end{equation}
(2)\,\,(Blow up case) If $T_+<\infty$, denote $\ell=\lim\limits_{t\to T_+}\psi(T_+-t,t)$ (it is well defined). Then $J\ge 1$, $\lambda_{1,n}\ll T_+-t_n$ and there exists a function $\OR{\phi}$ with $\E(\OR{\phi})<\infty$ such that $Q_1(\infty)=\phi(0)=\ell$ and 
\begin{equation}
\OR{\psi}(t_n)=\sum_{j=1}^J(Q_j(\cdot/r_{jn})-Q_j(\infty),0)+\OR{\phi}+\OR{b}_n.
\end{equation}
In both cases, $\OR{b}_n=(b_{0,n},b_{1,n})$ satisfies $\|b_n\|_{H\times L^2}=O(1)$ and vanishes in the following sense
\begin{equation}\label{eq:heheinequalitysec2}
\|b_{0,n}\|_{L^{\infty}}+\|b_{1,n}\|_{L^2(R^+,\,rdr)}\to 0\,\,{\rm as}\,\,n\to\infty,
\end{equation}
and moreover, for any sequence $\lambda_n>0$ and $A>1$
\begin{equation}\label{eq:dyadicsec2}
\|b_{0,n}\|_{H(\frac{\lambda_n}{A}\leq r\leq A\lambda_n)}\to 0,\,\,{\rm as}\,\,n\to\infty.
\end{equation}
\end{theorem}

\noindent
{\it Proof.} For simplicity, let us only consider the finite time blow up case. (The global existence case is similar, we only need to use Lemma \ref{lm:thetrickgeneralglobal} instead of Lemma \ref{lm:thetrickgeneral}.) Assume without loss of generality $T_+=1$. Take the sequence $t_n\uparrow 1$ from Lemma \ref{lm:thetrickgeneral}. Then we have 
\begin{equation}\label{eq:averagezerosec2}
\sup_{\sigma\in(0,1-t_n)}\frac{1}{\sigma}\int_{t_n}^{t_n+\sigma}\int_0^{\infty}(\partial_ta)^2\,rdrds=o_n(1),
\end{equation}
where $\OR{a}(r,t)=\OR{\psi}(r,t)-\OR{\phi}(r,t)+\OR{\phi}(0,t)$ and $\OR{\phi}$ is the regular part of $\OR{\psi}$. We note that $\OR{a}(r,t)=\OR{\phi}(r,t)$ for $r\ge 1-t$. (\ref{eq:averagezerosec2}) is a consequence of the fact that asymptotically there is no energy in the self similar region $r\sim 1-t$ and a virial type identity (see Lemma 2.2 and Corollary 2.3 of \cite{Shatah1}), together with a real analysis lemma (see Corollary 5.3 of \cite{KenigYang}).  See also Section 3 and Section 4 for the proof of Lemma \ref{lm:thetrickgeneral}. Set $k=|g'(\psi(0,t))|$. By (\ref{eq:averagezerosec2}) and the remark below Lemma \ref{lm:decompositionforwavemap}, passing to a subsequence of $t_n$ if necessary, we can assume $\OR{a}(t_n)=\OR{\psi}(t_n)-\OR{\phi}(t_n)+\OR{\phi}(0,t_n)$ has the following profile decomposition for each $K$
\begin{equation}\label{eq:decfirst}
\def\arraystretch{2.2}
\begin{array}{rl}
\OR{\psi}(t_n)-\OR{\phi}(t_n)+\OR{\phi}(0,t_n)=&(\sum_{j=1}^{J}(\psi^j_0(\frac{\cdot}{r_{jn}})-\psi^j_0(0))+\psi_0^J(0),0)\\
&+r^k\sum_{l=1}^K\left(\frac{1}{\lambda^k_{ln}}U^l_L(\frac{r}{\lambda_{ln}},0),0\right)+(\zeta^K_{0,n},\zeta^K_{1,n})
               \end{array}
\end{equation}
where $\psi^1_0(\infty)=\ell=\phi(0,t),\,\,\psi^{j+1}_0(\infty)=\psi^j_0(0),\,\,{\rm for\,\,each\,\,} 1\leq j\leq J-1,\,\psi^J_0(0)=\psi(0,t_n)$ and the parameters satisfy ``pseudo-orthogonality" conditions
\begin{eqnarray*}
&&0<r_{Jn}\ll r_{J-1,n}\ll\cdots\ll r_{1n};\\
&&\frac{\lambda_{ln}}{r_{jn}}+\frac{r_{jn}}{\lambda_{ln}}\to\infty,\,\,{\rm as}\,\,n\to\infty \,\,{\rm for\,\,each\,\,}l\,\,{\rm and}\,\,1\leq j\leq J;\\
&&\lim_{n\to\infty}\left(\frac{\lambda_{ln}}{\lambda_{l'n}}+\frac{\lambda_{l'n}}{\lambda_{ln}}\right)=\infty, \,\,\,{\rm for\,\,all}\,\,l\neq l'.
\end{eqnarray*}
$(\zeta^K_{0,n},\zeta^K_{1,n})$ vanishes asymptotically in the following senses:
\begin{equation}
\lim_{K\to\infty}\limsup_{n\to\infty}\|\zeta^K_{0,n}\|_{L^{\infty}}+\lim_{n\to\infty}\|\zeta^K_{1,n}\|_{L^2(R_+,\,rdr)}= 0,
\end{equation}
and
\begin{equation}
\left(\zeta^K_{0,n}(r_{jn}\cdot,),\,r_{jn}\zeta^K_{1,n}(r_{jn}\cdot)\right)\rightharpoonup 0,\,\,\forall 1\leq j\leq J,\,\,{\rm and}\,\,\left(\zeta^J_{0,n}(\lambda_{ln}\cdot),\,\lambda_{ln}\zeta_{1,n}^J(\lambda_{ln}\cdot)\right)\rightharpoonup 0,\,\,\forall 1\leq l\leq K,
\end{equation}
as $n\to\infty$.
By Lemma 2.2 in \cite{Shatah1} (see also Lemma \ref{lm:noselfsimilarenergymap}) there is asymptotically no energy in the self similar region $r\sim 1-t$ as $t\to 1-$, thus we can assume that $r_{1n}\ll 1-t_n$ and $\lambda_{jn}\ll 1-t_n$ as $n\to\infty$.
An important implication of the decomposition (\ref{eq:decfirst}) is that for each $j,\,l$
\begin{equation}\label{eq:convergenceapp}
\psi(r_{jn}\cdot,t_n)\rightharpoonup \psi^j_0,\,\,\,{\rm and}\,\,\,\psi(\lambda_{ln}\cdot,t_n)\rightharpoonup U^l_L(0).
\end{equation}
Since $\OR{\phi}(t)$ is continuous in the energy space, (\ref{eq:averagezerosec2}) implies
\begin{equation}\label{eq:zerosec2}
\sup_{\sigma\in(0,1-t_n)}\frac{1}{\sigma}\int_{t_n}^{t_n+\sigma}\int_0^{1-s}(\partial_t\psi)^2\,rdrds=o_n(1).
\end{equation}
By appropriate rescaling and time translations, (\ref{eq:zerosec2}) implies for each $T>0$ and sufficiently large $n$ that
\begin{eqnarray*}
&&\frac{1}{T}\int_0^T\int_0^{\frac{1-t_n}{2r_{jn}}}\left(r_{jn}\partial_t\psi(r_{jn}r,r_{jn}t+t_n)\right)^2\,rdrdt=o_n(1);\\
&&\frac{1}{T}\int_0^T\int_0^{\frac{1-t_n}{2\lambda_{ln}}}\left(\lambda_{ln}\partial_t\psi(\lambda_{ln}r,\lambda_{ln}t+t_n)\right)^2\,rdrdt=o_n(1).
\end{eqnarray*}
Hence the fact that $\lambda_{ln}\ll 1-t_n$ and $r_{jn}\ll 1-t_n$ for each $l,\,j$, the convergence (\ref{eq:convergenceapp}) and Lemma \ref{lm:struwe} imply that $\psi^j_0$ and $r^kU^l_L(\cdot,0)$ must be harmonic maps. Thus $\psi^j_0=Q_j$, $U^l_L\equiv 0$ for each $l$ and consequently $(\zeta^K_{0,n},\zeta^K_{1,n})$ does not depend on $K$, and there are no additional profiles with $t_{ln}\equiv 0$. Recalling that $\psi^j_0(0)=\psi^{j+1}_0(\infty)$ for $1\leq j\leq J-1$ and $\psi^1_0(\infty)=\ell=\phi(0,t)$, we can then simplify the decomposition (\ref{eq:decfirst}) as follows
\begin{equation*}
\OR{\psi}(t_n)=\sum_{j=1}^J(Q_j(\cdot/r_{jn})-Q_j(\infty),0)+\OR{\phi}(t_n)+(\zeta_{0,n},\zeta_{1,n}).
\end{equation*}
The bound (\ref{eq:heheinequalitysec2}) on the dispersive part $(b_{0,n},b_{1,n}):=(\zeta_{0,n},\zeta_{1,n})$ follows from the property of the profile decomposition. To prove bound (\ref{eq:dyadicsec2}), by the fact that there is no self similar energy and the regularity property of $\psi$ for $r\ge 1-t$, it is clear we only need to consider the case $\lambda_n\ll 1-t_n$ and $\lambda_n\psi(\lambda_n\cdot,t_n)\rightharpoonup 0$ as $n\to \infty$. Then (\ref{eq:dyadicsec2}) follows from the same argument as above. By the continuity of $\OR{\phi}(t)$ in $H\times L^2$, the Theorem is proved with $\OR{\phi}:=\OR{\phi}(\cdot,1)$.\\

\noindent
{\it Remark.} Our estimate on the dispersive part $(b_{0,n},b_{1,n})$ is slightly weaker than the following bound from \cite{cotesoliton}
\begin{equation*}
\|\partial_rb_{0,n}\|_{L^2(r\leq r_{Jn},\,rdr)}+\|b_{0,n}\|_{L^{\infty}}+\|b_{1,n}\|_{L^2(R^+,\,rdr)}\to 0\,\,{\rm as}\,\,n\to\infty,
\end{equation*}
 however it is already be sufficient for our arguments in the paper. One of the main achievements in this paper is to show that the dispersive part $\OR{b}_n$ actually vanishes in the energy space as $n\to\infty$.  \\

Below we prove Lemma \ref{lm:decompositionforwavemap}.
Let $(\psi_{0,n},\psi_{1,n})$ be a sequence with $M:=\sup\limits_{n\ge 1}\,\E(\psi_{0,n},\psi_{1,n})<\infty$ and $\psi_{0,n}(\infty)=m\in\mathcal{V}$. We distinguish between two cases.\\ 

\smallskip
\noindent
{\it The special case:\,  $\sup_n\,\sup_{r>0}|\psi_{0,n}(r)-m|<\frac{1}{2}\inf\limits_{l\in\mathcal{V},\,l\neq m}|l-m|$.} \\
Then by the non-vanishing of $g'(m)$, we see 
\begin{equation}\label{eq:equivalentbound}
\E(\psi_{0,n},\psi_{1,n})\sim \int_0^{\infty}\left((\partial_r\psi_{0,n})^2+\psi_{1,n}^2+\frac{(\psi_{0,n}-m)^2}{r^2}\right)\,rdr.
\end{equation}
Let $|g'(m)|=k$. If we set $(u_{0,n},u_{1,n}):=\left(\frac{\psi_{0,n}-m}{r^k},\frac{\psi_{1,n}}{r^k}\right)$, by (\ref{eq:equivalentbound})  
\begin{equation}
\sup_{n\ge 1}\int_{R^{2k+2}}|\nabla u_{0,n}|^2+u_{1,n}^2\,dx<\infty.
\end{equation}
By the profile decomposition of Bahouri-Gerard \cite{BaGe}, we conclude there exists a subsequence of $(u_{0,n},u_{1,n})$ (which we still denote as $(u_{0,n},u_{1,n})$ for ease of notation), finite energy radial solutions $U^j_L(r,t)$ to the $2k+2$ dimensional linear wave equation, scale sequences $\lambda_{jn}>0$, and time shifts $t_{jn}$, such that
\begin{equation}
(u_{0,n},u_{1,n})=\sum_{j=1}^J\left(\frac{1}{\lambda^k_{jn}}U^j_L(\frac{r}{\lambda_{jn}},\frac{-t_{jn}}{\lambda_{jn}}),\frac{1}{\lambda_{jn}^{k+1}}\partial_tU^j_L(\frac{r}{\lambda_{jn}},\frac{-t_{jn}}{\lambda_{jn}})\right)+(w^J_{0,n},w^J_{1,n}).
\end{equation}
Moreover, let $w^J_n(r,t)$ be the solution to the $2k+2$ dimensional linear wave equation with initial data $(w^J_{0,n},w^J_{1,n})$, then $w^J_n,\,\lambda_{jn}$ and $t_{jn}$ satisfies the following ``pseudo-orthogonality" condition:
\begin{eqnarray*}
&&t_{jn}=0\,\,{\rm for\,\,all\,\,}n,\,\,{\rm or}\,\,\lim_{n\to\infty}\frac{t_{jn}}{\lambda_{jn}}\in\{\pm\infty\};\\
&&\lim_{n\to\infty}\left(\frac{\lambda_{jn}}{\lambda_{j'n}}+\frac{\lambda_{j'n}}{\lambda_{jn}}+\frac{|t_{jn}-t_{j'n}|}{\lambda_{jn}}\right)=\infty, \,\,\,{\rm for\,\,all}\,\,j\neq j';\\
&&{\rm write}\,\,w^J_{n}(x,t)=\frac{1}{\lambda_{jn}^2}\widetilde{w}^j_{Jn}(\frac{x}{\lambda_{jn}},\,\frac{t-t_{jn}}{\lambda_{jn}}),\,\,{\rm then}\,\,\widetilde{w}^j_{Jn}\rightharpoonup 0,\,\,{\rm as}\,\,n\to\infty, \,\,{\rm for}\,\,1\leq j\leq J.
\end{eqnarray*}
In addition, the dispersive part $w^J_{n}$ satisfies
\begin{equation}
\lim_{J\to\infty}\limsup_{n\to\infty}\|w^J_{n}\|_{L^{\infty}_tL^{\frac{2k+2}{k}}_x\cap L^{\frac{k+2}{k}}_tL^{\frac{2k+4}{k}}_x(R^{2k+2}\times [0,\infty))}=0.
\end{equation}
Let us also note
\begin{equation}
\lim_{J\to\infty}\limsup_{n\to\infty}\|r^kw^J_{n}\|_{L^{\infty}}=0.
\end{equation}
Hence we get the following profile decomposition for $(\psi_{0,n},\psi_{1,n})$
\begin{equation}
(\psi_{0,n},\psi_{1,n})=r^k\sum_{j=1}^J\left(\frac{1}{\lambda^k_{jn}}U^j_L(\frac{r}{\lambda_{jn}},\frac{-t_{jn}}{\lambda_{jn}}),\frac{1}{\lambda_{jn}^{k+1}}\partial_tU^j_L(\frac{r}{\lambda_{jn}},\frac{-t_{jn}}{\lambda_{jn}})\right)+r^k(w^J_{0,n},w^J_{1,n})+(m,0).
\end{equation}

\bigskip
\noindent
{\it The general case.}\\
Now let us deal with the general case, only in this case the special features of wave maps appear. We can assume
\begin{equation*}
\sup_{n\ge 1}\sup_{r>0}|\psi_{0,n}(r)|\leq C(M,m).
\end{equation*}
Set
\begin{equation}
\delta=\frac{1}{4}\inf\left\{|l_1-l_2|:\,l_1,\,l_2\in\mathcal{V},\,\,l_1\neq l_2,\,\,{\rm with}\,\,l_1,\,l_2\in[-2C(M,m),2C(M,m)]\right\}.
\end{equation}
Choose a sequence of positive numbers 
\begin{equation}
0<\delta_1<\cdots<\delta_j<\delta_{j+1}<\frac{\delta}{2},\,\,\,\forall j\ge 1.
\end{equation}
Let us define 
\begin{equation}
r_{1n}=\sup\,\{r>0:\,|\psi_{0,n}(r)-m|\ge \delta_1\}.
\end{equation}
If the set on the right hand side is empty for a subsequence, then passing to the subsequence we will be in the special case and we can obtain profile decompositions using techniques there. If otherwise, set $\widetilde{\psi}_{0,n}=\psi_{0,n}(r_{1n}\cdot)$, then we see by scaling invariance of the energy and definition of $r_{1n}$
\begin{equation}
\sup_n\E(\tpsi_{0,n},0)\leq M, \,\,\,{\rm and}\,\,\,|\tpsi_{0,n}(1)-m|=\delta_1.
\end{equation}
Passing to a subsequence if necessary we can assume $\tpsi_{0,n}\to \psi^1_0$ locally uniformly in $R^+$.
Then $|\psi^1_0(1)-m|=\delta_1$ and $\E(\psi^1_0,0)\leq M$. By Lemma \ref{lm:concentrationofenergy} we also have $\E(\psi^1_0,0)>\epsilon>0$ with some constant $\epsilon$ depending only on $g$, $m$ and $M$ (and the choice of the sequence $\delta_j$). Now let us consider the new sequence 
\begin{equation}
h_{0,n}=\psi_{0,n}-\psi^1_0(\frac{\cdot}{r_{1n}})+\psi^1_0(0). 
\end{equation}
For any $L>2$ we have
\begin{equation}\label{eq:asy}
\lim_{n\to\infty}\sup_{\frac{r_{1n}}{L}\leq r\leq Lr_{1n}}\left|\psi_{0,n}(r)-\psi^1_0(\frac{r}{r_{1n}})\right|=0.
\end{equation}
If we choose $L>2$ sufficiently large, such that
\begin{equation*}
\sup_{r\ge \frac{1}{2}Lr_{1n}}\left|\psi^1_0(\frac{r}{r_{1n}})-\psi^1_0(\infty)\right|<\frac{\delta_1}{2}.
\end{equation*}
By (\ref{eq:asy}), for sufficiently large $n$, we then have
\begin{equation*}
\sup_{\frac{Lr_{1n}}{2}\leq r\leq Lr_{1n}}|\psi_{0,n}(r)-\psi^1_0(\infty)|<\delta_1.
\end{equation*}
Note that by definition of $r_{1n}$ we have
\begin{equation*}
\sup_{r\ge \frac{1}{2}Lr_{1n}}|\psi_{0,n}(r)-m|<\delta_1.
\end{equation*}
Thus $|m-\psi^1_0(\infty)|<2\delta_1<\delta$ and by definition of $\delta$, we must have $m=\psi^1_0(\infty)$. Therefore the new sequence $h_{0,n}$ satisfies for sufficiently large $n$,
\begin{equation*}
h_{0,n}(\infty)=\psi^1_0(0).
\end{equation*}
We observe, for any $L>1$ and sufficiently large $n$
\begin{equation}\label{eq:goodapproximate}
\sup_{r\ge \frac{r_{1n}}{L}}\,|h_{0,n}(r)-h_{0,n}(\infty)|=\sup_{r\ge \frac{r_{1n}}{L}}\,\left|\psi_{0,n}(r)-\psi^1_0(\frac{r}{r_{1n}})\right|< \delta_2.
\end{equation}
Note $\psi_{0,n}(\infty)=\psi^1_0(\infty)$, (\ref{eq:goodapproximate}) follows from
\begin{eqnarray*}
&&\sup_{r\ge \frac{r_{1n}}{L}}\,\left|\psi_{0,n}(r)-\psi^1_0(\frac{r}{r_{1n}})\right|\\
&&\leq \sup_{\frac{r_{1n}}{L}\leq r\leq Mr_{1n}}\,\left|\psi_{0,n}(r)-\psi^1_0(\frac{r}{r_{1n}})\right|+\\
&&\quad\quad+\sup_{r\ge Mr_{1n}}\,\left(|\psi_{0,n}(r)-\psi_{0,n}(\infty)|+|\psi^1_0(\frac{r}{r_{1n}})-\psi^1_0(\infty)|\right)\\
&&\leq o_n(1)+\delta_1+\sup_{r\ge Mr_{1n}}\,|\psi^1_0(\frac{r}{r_{1n}})-\psi^1_0(\infty)|\\
&&<\delta_2,\,\,{\rm if\,\,}M{\rm\,\,is\,\,taken\,\,large,\,\,for\,\,sufficiently\,\,large\,\,}n.
\end{eqnarray*}
We claim that the extraction of the bubble $\psi^1_0$ decreases energy up to the scale $O(r_{1n})$, more precisely we have the following lemma.
\begin{lemma}\label{lm:extractiondecreaseenergy}
We have for large $n$
\begin{equation}\label{eq:energyfinite}
\E(h_{0,n},0)\leq C \E(\psi_{0,n},0).
\end{equation}
Furthermore, we can find positive $\beta_{1n}\to \infty$ as $n\to\infty$, such that
\begin{equation}\label{eq:energydecomposition}
\int_{r\leq \beta_{1n}r_{1n}}(\partial_rh_{0,n})^2+\frac{g^2(h_{0,n})}{r^2}\,rdr\leq \E(\psi_{0,n},0)-\E(\psi^1_0,0)+o_n(1).
\end{equation}
\end{lemma}

\smallskip
\noindent
{\it Remark.} Due to the term $\psi^1_0(0)$ which changes the asymptotic behavior of $\psi_{0,n}$ for $r\ge \beta_{1n}r_{1n}$, the energy may in fact be increased slightly after the extraction of the bubble.\\

\smallskip
\noindent
{\it Proof.} The proof is not completely trivial due the term $\int_0^{\infty}\frac{g^2(h_{0,n})}{r^2}\,rdr$, which does not admit explicit expansion. Recall $h_{0,n}=\psi_{0,n}-\psi^1_0(\frac{r}{r_{1n}})+\psi^1_0(0)$. Since 
\begin{equation*}
\sup_{\frac{r_{1n}}{L}\leq r\leq Lr_{1n}}\,\left|\psi_{0,n}-\psi^1_0(\frac{r}{r_{1n}})\right|=o_n(1)
\end{equation*}
for all $L>1$, we can choose $\beta_{1n}\to\infty$ as $n\to+\infty$, such that
\begin{eqnarray*}
&&\sup_{\frac{r_{1n}}{\beta_{1n}}\leq r\leq \beta_{1n}r_{1n}}|\psi_{0,n}-\psi^1_0|=o_n(1);\\
&&\int_{\frac{r_{1n}}{\beta_{1n}}}^{\beta_{1n}r_{1n}}\frac{|\psi_{0,n}-\psi^1_0(\frac{r}{r_{1n}})|^2}{r^2}\,rdr=o_n(1).
\end{eqnarray*}
Note by definition of $\psi^1_0$, 
\begin{equation*}
r_{1n}\partial_rh_{0,n}(r_{1n}r)\rightharpoonup 0, 
\end{equation*}
weakly in $L^2((0,\infty),rdr)$. Hence
\begin{eqnarray}
\int_0^{\infty}(\partial_r\psi_{0,n})^2\,rdr&=&\int_0^{\infty}\left((\partial_r\psi^1_0)^2+2\partial_r\psi^1_0\cdot r_{1n}\partial_r h_{0,n}(r_{1,n}r)+(\partial_rh_{0,n})^2\right)\,\,rdr\nonumber\\
&=&\int_0^{\infty}\left((\partial_r\psi^1_0)^2+(\partial_rh_{0,n})^2\right)\,rdr+o_n(1).\label{eq:firstorthogonality}
\end{eqnarray}
Let us consider the term $\int_0^{\infty}\frac{g^2(h_{0,n})}{r^2}\,rdr$. We have
\begin{eqnarray*}
&&\int_0^{\infty}\frac{g^2(h_{0,n})}{r^2}\,rdr\\
&&=\int_{r\leq \frac{r_{1n}}{\beta_{1n}}}\frac{g^2(h_{0,n})}{r^2}\,rdr+\int_{\frac{r_{1n}}{\beta_{1n}}\leq r\leq \beta_{1n}r_{1n}}\frac{g^2(h_{0,n})}{r^2}\,rdr\\
&&\quad\quad\quad\quad+\int_{r\ge \beta_{1n}r_{1n}}\frac{g^2(h_{0,n})}{r^2}\,rdr\\
&&=I_1+I_2+I_3.
\end{eqnarray*}
For $r\leq \frac{r_{1n}}{\beta_{1n}}$, 
\begin{equation*}
h_{0,n}(r)=\psi_{0,n}(r)-\psi^1_0(\frac{r}{r_{1n}})+\psi^1_0(0)=\psi_{0,n}(r)+o_n(1).
\end{equation*}
Noting that $g\in C^3$ and $\psi_{0,n}$ are uniformly bounded, by the elementary inequality
\begin{equation*}
|g^2(x+y)-g^2(x)-2g(x)g'(x)y|\lesssim |y|^2,
\end{equation*}
the non-vanishing of $g'(\psi^1_0(0))$ and H\"older inequality, we obtain
\begin{eqnarray*}
I_1&=&\int_{r\leq \frac{r_{1n}}{\beta_{1n}}}\frac{g^2(\psi_{0,n})}{r^2}\,rdr\\
&&\quad\quad+\int_{r\leq \frac{r_{1n}}{\beta_{1n}}}\frac{g(\psi_{0,n})g'(\psi_{0,n})(\psi^1_0(0)-\psi^1_0(\frac{r}{r_{1n}}))}{r^2}\,rdr\\
&&\quad\quad\quad\quad+O\left(\int_{r\leq \frac{r_{1n}}{\beta_{1n}}}\frac{|\psi^1_0(\frac{r}{r_{1n}})-\psi^1_0(0)|^2}{r^2}\,rdr\right)\\
&=& \int_{r\leq \frac{r_{1n}}{\beta_{1n}}}\frac{g^2(\psi_{0,n})}{r^2}\,rdr+O\left(\int_{r\leq \frac{r_{1n}}{\beta_{1n}}}\frac{g^2(\psi^1_0(\frac{r}{r_{1n}}))}{r^2}\,rdr\right)\\
&&\quad\quad\quad\quad+O\left(\int_{r\leq \frac{r_{1n}}{\beta_{1n}}}\frac{g^2(\psi^1_0(\frac{r}{r_{1n}}))}{r^2}\,rdr\right)\\
&=&\int_{r\leq \frac{r_{1n}}{\beta_{1n}}}\frac{g^2(\psi_{0,n})}{r^2}\,rdr+o_n(1).
\end{eqnarray*}
For $\frac{r_{1n}}{\beta_{1n}}\leq r\leq \beta_{1n}r_{1n}$,
\begin{equation*}
h_{0,n}=\psi_{0,n}-\psi^1_0(\frac{r}{r_{1n}})+\psi^1_0(0)=\psi^1_0(0)+o_n(1).
\end{equation*}
Hence,
\begin{eqnarray*}
I_2&=&\int_{\frac{r_{1n}}{\beta_{1n}}\leq r\leq \beta_{1n}r_{1n}}\frac{g^2(h_{0,n})}{r^2}\,rdr\\
    &\lesssim& \int_{\frac{r_{1n}}{\beta_{1n}}\leq r\leq \beta_{1n}r_{1n}}\frac{|\psi_{0,n}-\psi^1_0(\frac{r}{r_{1n}})|^2}{r^2}\,rdr\\
&=&o_n(1).
\end{eqnarray*}
Note that the above arguments also imply
\begin{eqnarray*}
&&\int_{\frac{r_{1n}}{\beta_{1n}}\leq r\leq \beta_{1n}r_{1n}}\frac{g^2(\psi_{0,n})}{r^2}\,rdr\\
&&=\int_{\frac{r_{1n}}{\beta_{1n}}\leq r\leq \beta_{1n}r_{1n}}\frac{g^2(\psi^1_0(\frac{r}{r_{1n}})+\psi_{0,n}-\psi^1_0(\frac{r}{r_{1n}}))}{r^2}\,rdr\\
&&= \int_{\frac{r_{1n}}{\beta_{1n}}\leq r\leq \beta_{1n}r_{1n}}\frac{g^2(\psi^1_0(\frac{r}{r_{1n}}))}{r^2}\,rdr+o_n(1)\\
&&=\int_0^{\infty}\frac{g^2(\psi^1_0)}{r^2}\,rdr+o_n(1).
\end{eqnarray*}
For $r\ge \beta_{1n}r_{1n}$ and sufficiently large $n$,
\begin{eqnarray*}
&&|h_{0,n}-\psi^1_0(0)|\\
&&=\left|\psi_{0,n}-\psi^1_0(\frac{r}{r_{1n}})\right|\\
&&\leq |\psi_{0,n}(r)-\psi_{0,n}(\infty)|+\left|\psi^1_0(\frac{r}{r_{1n}})-\psi^1_0(\infty)\right|\\
&&<\delta_1+o_n(1)\\
&&< \delta_2.
\end{eqnarray*}
Hence,
\begin{eqnarray*}
I_3&=&\int_{r\ge \beta_{1n}r_{1n}}\frac{g^2(h_{0,n})}{r^2}\,rdr\\
&\lesssim& \int_{r\ge\beta_{1n}r_{1n}}\frac{\left(\psi_{0,n}-\psi^1_0(\frac{r}{r_{1n}})\right)^2}{r^2}\,rdr\\
&\lesssim&  \int_{r\ge\beta_{1n}r_{1n}}\left(\frac{\left(\psi_{0,n}-\psi_{0,n}(\infty)\right)^2}{r^2}+\frac{\left(\psi^1_0(\frac{r}{r_{1n}})-\psi^1_0(\infty)\right)^2}{r^2}\right)\,rdr\\
&\lesssim& \int_{r\ge\beta_{1n}r_{1n}}\left(\frac{g^2(\psi_{0,n})}{r^2}+\frac{g^2(\psi^1_0(\frac{r}{r_{1n}}))}{r^2}\right)\,rdr\\
&\lesssim& \int_{r\ge\beta_{1n}r_{1n}}\frac{g^2(\psi_{0,n})}{r^2}\,rdr+o_n(1),
\end{eqnarray*}
by the choice of $\beta_{1n}$. The estimates on $I_1,\,I_2,\,I_3$ and (\ref{eq:firstorthogonality}) finish the proof of (\ref{eq:energyfinite}).
Combining $I_1,\,I_2$, we obtain
\begin{eqnarray*}
&&\int_{r\leq \beta_{1n}r_{1n}}\frac{g^2(h_{0,n})}{r^2}\,rdr\\
&&=\int_{r\leq \frac{r_{1n}}{\beta_{1n}}}\frac{g^2(\psi_{0,n})}{r^2}\,rdr+o_n(1)\\
&&\leq \int_0^{\infty}\frac{g^2(\psi_{0,n})}{r^2}\,rdr-\int_0^{\infty}\frac{g^2(\psi^1_0)}{r^2}+o_n(1).
\end{eqnarray*}
This inequality, combined with (\ref{eq:firstorthogonality}) finishes the proof of the lemma.\\

Let us turn to the extraction of the second bubble. Define
\begin{equation*}
r_{2n}:=\sup\,\left\{r:\,|h_{0,n}-h_{0,n}(\infty)|>\delta_2\right\}.
\end{equation*}
By (\ref{eq:goodapproximate}), we see $r_{2n}=o(r_{1n})$ as $n\to\infty$. Passing to a subsequence, we can assume that $h_{0,n}(r_{2n}\cdot)\to \psi^2_0$ locally uniformly in $(0,\infty)$, for some $\psi^2_0$ satisfying 
\begin{equation*}
|\psi^2_0(1)-\psi^1_0(0)|=\delta_2,\,\,\,{\rm and}\,\,\E(\psi^2_0,0)\in(\epsilon,\,M).
\end{equation*}
In fact by $r_{2n}=o(r_{1n})$ and the formula for $h_{0,n}$, we actually have $\psi_{0,n}(r_{2n}\cdot)\to \psi^2_0$ locally uniformly in $(0,\infty)$. By arguments similar to those in the extraction of the first bubble, we conclude for any $L>1$ and sufficiently large $n$
\begin{eqnarray*}
&&h_{0,n}(\infty)=\psi^2_0(\infty)=\psi^1_0(0), \,\,\,{\rm and}\\
&&\sup_{r\ge \frac{r_{2n}}{L}}\left |\psi_{0,n}(r)-
\left(\psi^1_0(\frac{r}{r_{1,n}})-\psi^1_0(0)\right)-\psi^2_0(\frac{r}{r_{2n}})\right|< \delta_3.
\end{eqnarray*}
Re-define 
\begin{equation}
h_{0,n}=\psi_{0,n}-\left(\psi^1_0(\frac{r}{r_{1n}})-\psi^1_0(0)\right)-\left(\psi^2_0(\frac{r}{r_{2n}})-\psi^2_0(0)\right).
\end{equation}
We have for any $L>1$, sufficiently large $n$, and some positive $\beta_{2n}\to\infty$ as $n\to\infty$ with 
\begin{equation*}
\beta_{2n}r_{2n}\ll r_{1n},
\end{equation*}
such that
\begin{eqnarray*}
&&h_{0,n}(\infty)=\psi^2_0(0);\\
&&\sup_{r\ge \frac{r_{2n}}{L}}|h_{0,n}-h_{0,n}(\infty)|<\delta_3;\\
&&\E(h_{0,n},0)\leq C\E(\psi_{0,n},0),\,\,\,{\rm and}\,\,\\
&&\int_{r\leq \beta_{2n}r_{2n}}(\partial_rh_{0,n})^2+\frac{g^2(h_{0,n})}{r^2}\,rdr\leq \E(\psi_{0,n},0)-\E(\psi^1_0,0)-\E(\psi^2_0,0)+o_n(1).
\end{eqnarray*}

Then we can repeat this process to get $r_{3,n},\cdots,r_{Kn}$ and $\psi^3_0,\cdots,\psi^K_0$ and $\beta_{ln}\to\infty$ with the following properties
\begin{eqnarray*}
&&r_{1n}\gg r_{2n}\gg\cdots\gg r_{Kn},\,\,\,{\rm and}\,\,\,r_{l+1,n}\beta_{l+1,n}\ll r_{ln};\\
&&|\psi^l_0(1)-\psi^l_0(\infty)|=\delta_l;\\
&&\psi^l_0(0)=\psi^{l+1}_0(\infty),\,\,\forall 1\leq l\leq K-1;\\
&&\psi_{0,n}(r_{in}\cdot)-\sum_{l=1}^{i-1}\left(\psi^l_0\left(\frac{r_{in}\cdot}{r_{lh}}\right)-\psi^l_0(0)\right)\to \psi^{i}_0,\,\,\,{\rm locally\,\,uniformly\,\,in\,\,}(0,\infty);\\
&& {\rm or}\,\, {\rm equivalently},\,\, \psi_{0,n}(r_{in}\cdot)\to \psi^i_0, {\rm \,\,locally\,\, uniformly\,\, in}\,\, R^+;\\
&&\limsup_{n\to\infty}\sup_{r\ge \frac{r_{in}}{L}}\left|\psi_{0,n}-\sum_{l=1}^{i-1}(\psi^l_0(\frac{\cdot}{r_{lh}})-\psi^l_0(0))-\psi^i(\frac{\cdot}{r_{in}})\right|< \delta_{i+1}, \,\,\,\forall L>1.
\end{eqnarray*}
Morover, re-define 
\begin{equation*}
h_{0,n}=\psi_{0,n}-\sum_{l=1}^{K}(\psi^l_0(\frac{\cdot}{r_{lh}})-\psi^l_0(0)),
\end{equation*}
then
\begin{equation*}
\int_{r\leq \beta_{Kn}r_{Kn}}(\partial_rh_{0,n})^2+\frac{g^2(h_{0,n})}{r^2}\,rdr\leq \E(\psi_{0,n},0)-\sum_{l=1}^K\E(\psi^l_0,0)+o_n(1).
\end{equation*}
We note that each $\psi^l_0$ contains fixed amound of energy (independent of $l$), hence this process can not continue indefinitely, thus after extracting $K$ bubbles, we can no longer extract any more bubbles with this process, which implies
\begin{equation}
\limsup_{n\to\infty}\sup_{r>0}\,\left|\psi_{0,n}-\sum_{l=1}^{K}(\psi^l_0(\frac{\cdot}{r_{lh}})-\psi^l_0(0))-\psi_0^K(0)\right|\leq \delta_{K+1}<\frac{\delta}{2}.
\end{equation}
Hence the sequence $(\psi_{0,n}-\sum_{l=1}^{K}(\psi^l_0(\frac{\cdot}{r_{lh}})-\psi^l_0(0)),\psi_{1,n})$ verifies conditions in the special case, and by passing to a subsequence if necessary, (with $k=g'(\psi_0^K(0))$) we obtain the following profile decomposition
\begin{eqnarray*}
&&\left(\psi_{0,n}-\sum_{l=1}^{K}(\psi^l_0(\frac{\cdot}{r_{lh}})-\psi^l_0(0))-\psi^K(0),\psi_{1,n}\right)\\
&&=r^k\sum_{j=1}^J\left(\frac{1}{\lambda^k_{jn}}U^j_L(\frac{r}{\lambda_{jn}},\frac{-t_{jn}}{\lambda_{jn}}),\frac{1}{\lambda_{jn}^{k+1}}\partial_tU^j_L(\frac{r}{\lambda_{jn}},\frac{-t_{jn}}{\lambda_{jn}})\right)+r^k(w^J_{0,n},w^J_{1,n});
\end{eqnarray*} 
or equivalently
\begin{equation}
\def\arraystretch{2.2}
\begin{array}{rl}
(\psi_{0,n},\psi_{1,n})=&\left(\sum_{l=1}^{K}(\psi^l_0(\frac{\cdot}{r_{lh}})-\psi^l_0(0))+\psi_0^K(0),0\right)\\
&+r^k\sum_{j=1}^J\left(\frac{1}{\lambda^k_{jn}}U^j_L(\frac{r}{\lambda_{jn}},\frac{-t_{jn}}{\lambda_{jn}}),\frac{1}{\lambda_{jn}^{k+1}}\partial_tU^j_L(\frac{r}{\lambda_{jn}},\frac{-t_{jn}}{\lambda_{jn}})\right)+r^k(w^J_{0,n},w^J_{1,n})
               \end{array}
\end{equation}
where the parameters satisfy ``pseudo-orthogonality" conditions as in the special case and furthermore
\begin{eqnarray*}
&&r_{1n}\gg\cdots\gg r_{Kn}\,\,\,{\rm and\,\,if\,\,}t_{jn}\equiv 0,\,\,{\rm then}\,\,\frac{r_{ln}}{\lambda_{jn}}+\frac{\lambda_{jn}}{r_{ln}}\to\infty\,\,\,{\rm for\,\,each\,\,}j,\,l;\\
&&\psi^1_0(\infty)=m,\,\,\psi^{j+1}_0(\infty)=\psi^j_0(0),\,\,{\rm for\,\,each\,\,} 1\leq j\leq K-1;\\
&&{\rm write}\,\,\widetilde{w}^l_{Jn}:=r_{ln}^kw^J_{0,n}(r_{ln}\cdot),\,\,{\rm then}\,\,\widetilde{w}^l_{Jn}\rightharpoonup 0\,\,{\rm in}\,\,\dot{H}^1(R^{2k+2})\,\,{\rm as}\,\,n\to\infty,\,\,{\rm for}\,\,\forall l,\,J;\\
&&\lim_{J\to\infty}\limsup_{n\to\infty}\left(\|w^J_{n}\|_{L^{\infty}_tL^{\frac{2k+2}{k}}_x\cap L^{\frac{k+2}{k}}_tL^{\frac{2k+4}{k}}_x(R^{2k+2}\times [0,\infty))}+\|r^kw^J_{n}\|_{L^{\infty}}\right)=0.
\end{eqnarray*}

\end{subsection}

\end{section}

\end{document}